\documentclass[10pt,a4paper]{amsart} 
\usepackage{multicol}

\usepackage[utf8]{inputenc}
\usepackage[T1]{fontenc}

\usepackage[colorlinks=true, citecolor=purple, urlcolor=blue, linkcolor=purple%, pagebackref
]{hyperref}
\usepackage[noabbrev]{cleveref}
\usepackage{amsfonts}
\usepackage{amssymb}
\usepackage{amsthm} 

\usepackage[all]{xy}

\usepackage{mathtools} % coloneqq

\usepackage{mathdots} % iddots

\usepackage{mathrsfs}

\usepackage{dutchcal} % new calligraphic font with small letters

\usepackage{tikz}
\usetikzlibrary{cd}

\numberwithin{equation}{subsection}

\newtheorem{main}{Theorem}
\crefname{main}{theorem}{theorems}
\Crefname{main}{Theorem}{Theorems}

\newtheorem{conjecture}[main]{Conjecture}
\newtheorem*{conjecture*}{Conjecture}

\newtheorem{theorem}[equation]{Theorem}
\newtheorem{definition}[equation]{Definition}
\newtheorem{proposition}[equation]{Proposition}
\newtheorem*{proposition*}{Proposition}
\newtheorem{corollary}[equation]{Corollary}
\newtheorem*{corollary*}{Corollary}
\newtheorem{lemma}[equation]{Lemma}

\theoremstyle{definition}

\theoremstyle{remark}
\newtheorem{remark}[equation]{Remark}
\newtheorem*{remark*}{Remark}
\newtheorem{example}[equation]{Example}
\newtheorem*{example*}{Example}

\newtheorem*{examples*}{Examples}

\newcommand{\defeq}{\coloneqq}

\newcommand{\hooklongrightarrow}{\lhook\joinrel\longrightarrow}
\newcommand{\twoheadlongrightarrow}{\relbar\joinrel\twoheadrightarrow}

\newcommand{\C}{\mathbf{C}}

\newcommand{\Z}{\mathbf{Z}}

\DeclareMathOperator{\gr}{gr}

\DeclareMathOperator{\colim}{colim}

\DeclareMathOperator{\Ker}{Ker}
\DeclareMathOperator{\Fil}{F}
\newcommand{\Li}{\mathrm{Li}}
\newcommand{\FilLi}{\Fil_{\Li}}
\DeclareMathOperator{\Span}{Span}
\DeclareMathOperator{\Spec}{Spec}

\DeclareMathOperator{\Alt}{\textstyle\bigwedge}

\newcommand{\g}{\mathfrak{g}}
\newcommand{\h}{\mathfrak{h}}

\renewcommand{\sl}{\mathfrak{sl}}
\newcommand{\m}{\mathfrak m}
\newcommand{\n}{\mathfrak n}

\renewcommand{\l}{\mathfrak l}

\DeclareMathOperator{\ad}{ad}
\DeclareMathOperator{\Ad}{Ad}

\newcommand{\Orb}{\mathcal{O}}

\newcommand{\Jet}{\mathrm{J}}
\newcommand{\Hgy}{\mathrm{H}}

\newcommand{\A}{\mathcal{A}}
\newcommand{\Cpx}{\mathcal{C}}
\newcommand{\F}{\mathcal{F}}
\newcommand{\V}{\mathcal{V}}
\newcommand{\W}{\mathcal{W}}
\newcommand{\Ham}{\mathcal{H}}
\newcommand{\I}{\mathcal{I}}
\newcommand{\E}{\mathcal{E}}

\newcommand{\interm}{\mathrm{int}}
\newcommand{\new}{\mathrm{new}}
\newcommand{\old}{\mathrm{old}}
\newcommand{\HS}{\mathrm{HS}}

\newcommand{\One}{\mathbf 1}
\newcommand\NO[1]{\mathopen{\vcentcolon}#1\mathclose{\vcentcolon}}

\usepackage[backend=biber,style=alphabetic,autocite=plain, maxbibnames=10]{biblatex}
\title{Reduction by stages for affine W-algebras}

\author[N. Genra]{Naoki Genra}
\address{Faculty of Science, Academic Assembly, University of Toyama 3190 Gofuku, Toyama 930--8555, Japan}
\email{genra@sci.u-toyama.ac.jp}

\author[T. Juillard]{Thibault Juillard}
\address{Universität Hamburg, Hamburg, Germany}
\email{thibault.juillard@uni-hamburg.de}

\begin{document}

\begin{abstract}
    Given a pair of nilpotent orbits in a simple Lie algebra, one can associate a pair of vertex algebras called affine W-algebras. Under some compatibility conditions on these orbits, we prove that one of these W-algebras can be obtained as the quantum Hamiltonian reduction of the other. This property is called \emph{reduction by stages}. We provide several examples in classical and exceptional types.
    
    To prove reduction by stages for affine W-algebras, we use our previous work on reduction by stages for the Slodowy slices associated with these nilpotent orbits, these slices being the associated varieties of the W-algebras. We also prove and use the fact that each W-algebra can be defined using several equivalent BRST cohomology constructions: choosing the right BRST complexes allows us to connect the two W-algebras in a natural way.
\end{abstract}
	
\maketitle

\tableofcontents

\section{Introduction}

Vertex algebras are algebraic structures that were defined by Borcherds \cite{borcherds1986vertex}. They provide a rigorous axiomatisation of two-dimensional conformal field theories and play a major role in several areas of mathematics: the representation theory of infinite-dimensional Lie algebras, the Monstrous Moonshine Conjecture, the Langlands program, etc. Affine W-algebras form an important family of vertex algebras obtained as quantum Hamiltonian reductions of universal affine vertex algebras.

In this paper, we study a sufficient condition on a pair of affine W-algebras to ensure that one of them can be constructed as the quantum Hamiltonian reduction of the other one. This property is called \emph{reduction by stages}, by analogy with symplectic geometry \cite{marsden2007hamiltonian, morgan2014quantum, genra2024reduction}. Our main motivations for such result are: establishing connections between categories of modules of various W-algebras, constructing embeddings of W-algebras or isomorphisms between their simple quotients. Affine W-algebras are related to Slodowy slices by their natural Li filtration. By means of this relation, we will use reduction by stages for Slodowy slices, which was established in our previous work~\cite{genra2024reduction}, to prove reduction by stages for affine W-algebras.

\subsection{Context}

The \emph{affine W-algebra} $\W^k(\g, f)$ is a vertex algebra defined from the data of a finite-dimensional simple complex Lie algebra $\g$, a nilpotent element $f$ in~$\g$ and a complex number $k$ \cite{feigin1990quantization, kac2003quantum, kac2004quantum}. It is constructed by applying a BRST cohomology functor, denoted by $\Hgy^0_f$, to the universal affine vertex algebra~$\V^k(\g)$ associated with $\g$ and $k$ (\Cref{subsection:brst-vertex,subsection:definition-w-algebras}):
$$ \W^k(\g, f) \defeq \Hgy^0_f(\V^k(\g)). $$

This definition is the quantum affine analogue of the construction of the Slodowy slice~$S_f$ by Hamiltonian reduction of the dual vector space $\g^*$ equipped with its Kirillov--Kostant Poisson structure \cite{kostant1978whittaker, premet2002special, gan2002quantization}. The affine W-algebra is equipped with the Li filtration and the corresponding graded algebra is a Poisson vertex algebra \cite{li2005abelianizing}. Since $S_f$ is a Poisson variety, the coordinate ring of its arc space $\Jet_\infty S_f$ is also a Poisson vertex algebra~\cite{arakawa2012remark}. These two Poisson vertex algebras are isomorphic \cite{desole2006finite, arakawa2015associated}, hence geometric properties of Slodowy slices can be used to study affine W-algebras, see~\cite{arakawa2024arc} and the references in it.

The BRST cohomology functor can also be applied to the universal enveloping algebra $\mathcal{U}(\g)$, and the resulting noncommutative algebra, denoted by
$\mathcal{U}(\g, f)$, is called a \emph{finite W-algebra}. Finite W-algebras were studied in \cite{kostant1978whittaker,lynch1979generalized,deboer1994relation,ragoucy1999yangian,premet2002special, gan2002quantization,losev2010quantized}. They are natural quantisations of Slodowy slices \cite{premet2002special, gan2002quantization}. Affine and finite W-algebras are related by the Zhu algebra construction \cite{arakawa2007representation,desole2006finite}.

\subsection{Reduction by stages}

Let $f_1, f_2$ be two nilpotent elements in $\g$ and denote by $\mathbf{O}_1, \mathbf{O}_2$ their adjoint orbits in $\g$. Assume the inclusion $\overline{\mathbf{O}_1} \subseteq \overline{\mathbf{O}_2}$ of their Zariski closures and assume that the element~$f_0 \defeq f_2 - f_1$ is nilpotent. We say that \emph{reduction by stages} holds for the corresponding pair of affine W-algebras if there exists a BRST cohomology functor $\Hgy^0_{f_0}$ that can be applied to the first W-algebra to get the second one, up to a natural isomorphism of vertex algebras:
$$ \Hgy^0_{f_0}(\W^k(\g, f_1) ) \cong \W^k(\g, f_2).  $$
In other words, we have the commutative triangle:
$$ \begin{tikzcd}[, column sep=tiny]
	& \V^k(\g) \arrow[ld, mapsto, "\Hgy_{f_1}^0"'] \arrow[rd, mapsto,"\Hgy_{f_2}^0"] & \\
	\W^k(\g, f_1) \arrow[rr, mapsto, "\Hgy_{f_0}^0"'] & & \W^k(\g, f_2).
\end{tikzcd} $$

The analogous problem for the corresponding Slodowy slices was first studied by Morgan in his PhD thesis \cite{morgan2014quantum} and by the authors of the present paper in \cite{genra2024reduction}. In the latter, sufficient conditions are given on a pair of nilpotent elements~$f_1, f_2$ in $\g$ to get reduction by stages for the slices $S_{f_1}, S_{f_2}$, and also for the corresponding finite W-algebras $\mathcal{U}(\g, f_1), \mathcal{U}(\g, f_2)$. Let us mention that reduction by stages also appears in the context of slices in the affine Grassmaniann~\cite{kamnitzer2022hamiltonian,butson2025inverse}.

For affine W-algebras, reduction by stages was studied by Madsen and Ragoucy in the theoretical physics paper \cite{madsen1997secondary}, for the case when the simple Lie algebra $\g$ is of type $\mathrm A$ and the nilpotent elements $f_1, f_2$ correspond to hook-type partitions. Their approach consists in proving the isomorphism
$$ \Hgy^0_{f_0}(\Hgy^0_{f_1}(\V^k(\g))) \cong \Hgy^0_{f_2}(\V^k(\g))  $$
by noticing that the right-hand side cohomology is the total cohomology of a double cochain complex, and the left hand-side is the second page of the associated spectral sequence. We will generalise this approach.

Recently, reductions by stages of affine W-algebras (also called the \emph{partial reductions}) were intensively studied by using the free-field realizations of affine W-algebras developed by the first-named author in \cite{genra2020screening}. In general, this approach allows to study the reduction by stages with its associated \emph{inverse reduction}. Inverse reduction originally consists in ``inverting'' the BRST cohomology functor $\Hgy^0_f$ in the sense of constructing a vertex algebra embedding 
$$ \V^k(\g) \hooklongrightarrow \W^k(\g, f) \otimes_\C \mathcal{D}_{\mathrm{ch}}, $$
where $\mathcal{D}_{\mathrm{ch}}$ denotes the vertex algebra of chiral differential operators on an affine variety of the form $\C^{n_1} \times (\C^\times)^{n_2}$. These embeddings are useful tools for the representation theory of affine W-algebras \cite{adamovic2019realizations, adamovic2021realisation, adamovic2024relaxed}. This can be generalised to reduction by stages. In this context, an inverse reduction is an embedding of the form 
$$ \W^k(\g, f_1) \hooklongrightarrow \W^k(\g, f_2) \otimes_\C \mathcal{D}_{\mathrm{ch}}. $$

(Inverse) reductions by stages were studied in \cite{fehily2023subregular, fehily2024inverse} for type $\mathrm A$ and hook-type nilpotent elements; in \cite{fasquel2025orthosymplectic} for type $\mathrm B$ and subregular, regular nilpotent elements; in \cite{creutzig2025structure} for type $\mathrm A$ in small ranks; in \cite{fasquel2026connecting} for $\g = \sl_4$ and any ordered pair of nilpotent orbits; and in~\cite{fasquel2025virasoro} for $\g$ being of classical type and $\Hgy_{f_0}$ being a Virasoro-type reduction.

\subsection{Main results}

Let $\g$ be a simple Lie algebra and $\h$ be a Cartan subalgebra of~$\g$. For $i=1, 2$, let $H_i$ be in $\h$ such that the associated grading
$$ \g = \bigoplus_{\delta \in \Z} \g^{(i)}_\delta, \quad \text{where} \quad \g^{(i)}_\delta \defeq \{x \in \g ~ | ~ [H_i,x] = \delta x\}, $$
is a \emph{good grading} for the nilpotent element $f_i$ (\Cref{subsection:good-grading}). Since $[H_1, H_2]=0$, one gets a bigrading on $\g$:
$$ \g = \bigoplus_{\delta_1, \delta_2 \in \Z} \g_{\delta_1, \delta_2}, \quad \text{where} \quad \g_{\delta_1, \delta_2} \defeq \g^{(1)}_{\delta_1} \cap \g^{(2)}_{\delta_2}. $$
Set $f_0 \defeq f_2 - f_1$. Our sufficient conditions for reduction by stages are defined by:
\begin{equation}
    \begin{split} 
   & \g^{(1)}_{\geqslant 2}  \subseteq \g^{(2)}_{\geqslant 1} \subseteq \g^{(1)}_{\geqslant 0}, \quad  \g^{(1)}_1 \subseteq  \bigoplus_{\delta=0}^2 \g_{1, \delta}, \quad \g^{(2)}_1 \subseteq \bigoplus_{\delta=0}^2 \g_{\delta, 1}, \\
    & f_0 \in \g_{0, -2}.  \end{split}\tag{\ref{conditions}}
\end{equation}
The conditions~\eqref{conditions} are a refinement of the conditions used to prove reduction by stages for Slodowy slices in \cite[Main Theorems 1 and 2]{genra2024reduction}. These conditions imply the inclusion $\overline{\mathbf{O}_1} \subseteq \overline{\mathbf{O}_2}$ of nilpotent orbits closures. Our main result is the following.

\begin{main}[\Cref{theorem:reduction-by-stages-w-algebras}] \label{main:reduction-by-stages-w-algebras}
    If the conditions~\eqref{conditions} hold, there is a BRST cochain complex $\Cpx^\bullet_{f_0}(\W^k(\g, f_1))$ whose cohomology is isomorphic to $\W^k(\g, f_2)$ as vertex algebras: 
    $$ \Hgy^\bullet \big(\Cpx^\bullet_{f_0}(\W^k(\g, f_1))\big) \cong \mathsf{\delta}_{\bullet = 0} \, \W^k(\g, f_2). $$
\end{main}

\begin{examples*}   
    \Cref{main:reduction-by-stages-w-algebras} holds in the cases described in \Cref{table:examples}. By convention, the partitions indexing the nilpotent orbits in classical types are represented by nonincreasing sequences. A hook-type partition of the integer $n$ is a partition of the form $(a, 1^{n-a})$ for $1 \leqslant a \leqslant n$. See \Cref{subsection:examples} for details.
\end{examples*}

\begin{table}[h]
    \caption{Examples of reductions by stages}
    \centering 
     \begin{tabular}{|c||c||c||c|}
     \hline
         $\g$ & $f_1$ & $f_2$ & Reference \\\hline\hline
         type A & hook-type &  hook-type & \cite{madsen1997secondary, fehily2024inverse} \\\hline
         type $\mathrm{A}_{3}$ & \begin{tabular}{@{}c@{}}
            partition of $4$: \\
            $(2,1^2)$
         \end{tabular} & \begin{tabular}{@{}c@{}}
            partition of $4$: \\
            $(2,2)$
         \end{tabular} & \cite{creutzig2025structure} \\\hline
         type $\mathrm{A}_{n-1}$ & \begin{tabular}{@{}c@{}}
            partition of $n$: \\
            $(a_1, \dots, a_{r-1}, a_r, 1^p)$
         \end{tabular} & \begin{tabular}{@{}c@{}}
            partition of $n$: \\
            $(a_1, \dots, a_{r-1}, a_r + 1, 1^{p - 1})$
         \end{tabular} & \begin{tabular}{@{}c@{}}
            new \\
            for $n > 3$
         \end{tabular} \\\hline
         type B & subregular & regular & \cite{fasquel2025orthosymplectic} \\\hline 
         type $\mathrm{C}_r$ & \begin{tabular}{@{}c@{}}
            partition of $r$: \\
            $(2^2,1^{2r-4})$
         \end{tabular} & regular & new \\\hline
         type $\mathrm{G}_2$ & \begin{tabular}{@{}c@{}}
            Bala-Carter \\
            label $\widetilde{\mathrm{A}}_1$
         \end{tabular} & subregular & new \\\hline
    \end{tabular}  
    \label{table:examples}
\end{table}

To establish \Cref{main:reduction-by-stages-w-algebras}, we prove two important results which are also of independent interest. The first one is the definition of a new BRST complex, denoted by $\widetilde\Cpx^\bullet_{f_2}(\V^k(\g))$, to reconstruct the W-algebra $\W^k(\g, f_2)$, with an embedding 
$$ \Cpx^\bullet_{f_0}\big(\W^k(\g, f_1))\big) \hooklongrightarrow \widetilde\Cpx^\bullet_{f_2}(\V^k(\g)). $$
Therefore, one gets an induced vertex algebra homomorphism
$$ \Theta : \Hgy^0\big(\Cpx^\bullet_{f_0}(\W^k(\g, f_1))\big) \longrightarrow \Hgy^0\big(\widetilde\Cpx^\bullet_{f_2}(\V^k(\g))\big). $$

\begin{main}[{Theorem \ref{theorem:new-construction-w-algebra}}] \label{main:new-brst-complex}
    The cohomology of the complex $\widetilde\Cpx^\bullet_{f_2}(\V^k(\g))$ is isomorphic to the affine W-algebra $\W^k(\g, f_2)$:
    $$ \Hgy^\bullet\big(\widetilde\Cpx^\bullet_{f_2}(\V^k(\g)))\big) \cong \delta_{\bullet = 0} \, \W^k(\g, f_2). $$
\end{main}

Then there is a natural vertex algebra homomorphism
$$ \Theta : \Hgy^0\Big(\Cpx^\bullet_{f_0}(\W^k(\g, f_1))\Big) \longrightarrow \W^k(\g, f_2). $$
To show \Cref{main:new-brst-complex} and to prove that the vertex algebra homomorphism $\Theta$ is an isomorphism, we need to use the Li filtration and results about Slodowy slices. The associated graded map $\gr \Theta$ is an isomorphism as a consequence of the reduction by stages for the associated Slodowy slices $S_{f_1}, S_{f_2}$ and their arc spaces (\Cref{theorem:geometric-reduction-stages}). The Li filtration is not a complete filtration on the BRST cochain complex, so one has convergence issues. To overcome them, we generalise~\cite[Theorem~9.7]{arakawa2024arc}.

\begin{main}[Theorem \ref{theorem:vanishing-vertex}] \label{main:vanishing-vertex}
    Let $(\Cpx^\bullet, \mathcal{d})$ be a vertex algebra BRST cochain complex. Assume technical conditions which are roughly:
    \begin{itemize}
    	\item the existence of vertex algebra gradings on the complex to control the convergence of spectral sequence associated to the Li filtration,
    	
    	\item some good geometric conditions to get the vanishing of the cohomology of the associated graded  cochain complex $(\gr_{\Li} \Cpx^\bullet, \gr_{\Li} \mathcal{d})$.
    \end{itemize} 

    Then the cohomology vanishes in degrees other than $0$:
    $$ \Hgy^n(\Cpx^\bullet, \mathcal{d}) = 0 \quad \text{for} \quad n \neq 0. $$ 
    In degree 0 there is a natural isomorphism 
    $$ \gr_{\Fil}\Hgy^0(\Cpx^\bullet, \mathcal{d}) \overset{\sim}{\longrightarrow} \Hgy^0(\gr_{\Li} \Cpx^\bullet, \gr_{\Li} \mathcal{d}), $$
    where the filtration $\Fil$ on the cohomology $\Hgy^0(\Cpx^\bullet, \mathcal{d})$ is induced by the Li filtration on the complex $(\Cpx^\bullet, \mathcal{d})$.
\end{main}

\Cref{main:vanishing-vertex} holds for all the complexes mentioned above, allowing us to show \Cref{main:new-brst-complex}, and then \Cref{main:reduction-by-stages-w-algebras}. Another application of this theorem is to prove that all the possible BRST cohomology constructions of~$\W^k(\g, f)$ are equivalent, see \Cref{theorem:equivalent-constructions} or the remark before Theorem~9.7 in \cite{arakawa2024arc}. This proof is inspired by \cite{arakawa2015localization}, where equivalence of definitions is proved for the $\hbar$-adic W-algebra~$\W^k(\g, f)^\wedge_{\hbar}$.

\subsection{Future work and applications} \label{subsection:future}

In an upcoming work, we will extend \Cref{main:reduction-by-stages-w-algebras,main:new-brst-complex} to the~$\V^k(\g)$-modules in the Kazhdan--Lusztig category~$\mathsf{KL}^k(\g)$. 

Using this functorial version of reduction by stages, we aim to prove the fact that any W-algebra in type A can be reconstructed only using hook-type reduction~\cite[Conjecture A]{creutzig2025structure}. More precisely, let $\g$ be the simple Lie algebra $\sl_n$ and let $f$ be a nilpotent element corresponding to a partition~$(a_1, \dots, a_r)$ of $n$. Then the conjecture says that, for any module $V$ in~$\mathsf{KL}^k(\g)$, there is an isomorphism of~$\W^k(\sl_n, f)$-modules
$$ \Hgy^0_{f_r} \cdots \Hgy^0_{f_1}(V) \cong \Hgy^0_{f}(V), $$
where each $f_i$ is a nilpotent element in $\sl_{n_i}$ associated with the partition~$(a_i, 1^{n_i - a_i})$ of $n_i \defeq n - a_1 - \cdots a_{i-1}$.

This conjecture relates to a generalisation of the Kac--Roan--Wakimoto embedding. Let~$(a_1, \dots, a_s, a_{s+1}, \dots, a_r)$ be a partition of $n$ and let $f_2$ be a nilpotent element of $\sl_n$ corresponding to this partition. Let $f_1$ be a nilpotent element corresponding to $(a_1, \dots, a_s, 1^p)$, where $p \defeq  a_{s+1} + \cdots + a_r$. In general, $f_1, f_2$ do not satisfy the conditions \eqref{conditions}. Let $f_0$ be a nilpotent element of $\sl_p$ corresponding to the partition $(a_{s+1}, \dots, a_r)$ of $p$. By~\cite{kac2003quantum} (recalled in \Cref{proposition:krw-embedding}), there is a level $\mathcal{l}$ in $\C$ such that there is a vertex algebra embedding $\V^{\mathcal{l}}(\sl_p) \hookrightarrow \W^k(\sl_n, f_1)$.

\begin{conjecture} \label{conjecture:embedding}
    There is an isomorphism $\Hgy^0_{f_0}(\W^k(\g, f_1)) \cong \W^k(\g, f_2)$ of vertex algebras and, for any module $V$ in $\mathsf{KL}^k(\g)$, there is a $\W^k(\g, f_2)$-module isomorphism
    $$ \Hgy^0_{f_0}(\Hgy^0_{f_1}(V)) \cong \Hgy^0_{f_2}(V). $$
    Therefore, there is a vertex algebra embedding
    \begin{equation} \W^{\mathcal{l}}(\sl_p, f_0) \hooklongrightarrow \W^k(\sl_n, f_2). \label{equation:conjecture-embedding} \end{equation}
\end{conjecture}

Note that the analogue of \eqref{equation:conjecture-embedding} for finite W-algebras is known and plays an important role for their representation theory, see \cite{futorny2020gelfand}. It follows from the description of finite W-algebras as truncated shifted Yangians \cite{brundan2006shifted}. Let us mention that an additional motivation for \Cref{conjecture:embedding} is to prove isomorphisms of simple quotient of affine W-algebras at admissible level \cite[Conjecture 8.11]{arakawa2024singularities}, see~\cite[Section~1.4]{genra2024reduction} for more details.

\subsection{Recent developments}

This subsection gives more recent developments on reduction by stages since the submission of the present paper, in September 2025. Two proofs for inverse reduction and reduction by stages have come out, under the assumption that $k$ is generic.

%, under the assumptions that $\g$ is of type A, $f_1$ and $f_2$ belong to adjacent nilpotent orbits, and $k$ is generic.
\begin{itemize}
	\item Fasquel and Nakatsuka provide a proof using free-field realisation technics \cite{fasquel2026quantum}. Their conditions on $f_1$ and $f_2$ are similar to \eqref{conditions} \cite[Remark 4.2]{fasquel2026quantum}, but they only deal with even good grading. However, their conditions are general enough to get reduction by stages for $f_1$ and $f_2$ belonging to any pair of adjacent orbits in type A \cite[Theorem 1.3]{fasquel2026quantum}, and some non-adjacent orbits \cite[Remark 6.3]{fasquel2026quantum}.
	
	\item In \cite{butson2025inverse-affine}, Butson and Nair use the geometry of the affine Grassmannian \cite{butson2025inverse} and the theory of chiral quantisation they have developed in \cite{butson2026deformation}. Their reductions by stages correspond to cases in type A when the W-algebras are related by a reduction with respect to the action of a commutative nilpotent Lie algebra $\n_0$ (see \Cref{proposition:construction-groups} for the definition of $\n_0$). 
\end{itemize}  

Let us mention that our result on reduction by stages for affine W-algebras has been used to simplify the study of the Griess algebra inside a given affine W-algebra: see \cite[Section 6.2]{adamovic2026walgebras}.

\subsection{Outline}

Notations, conventions and tools about vertex algebras and Poisson geometry are introduced in \Cref{section:vertex-poisson}. 

In \Cref{section:vanishing-brst}, we explain BRST cohomology for Poisson varieties, their arc spaces and vertex algebras. We state vanishing results and prove \Cref{main:vanishing-vertex}~(\ref{theorem:vanishing-vertex}).

In \Cref{section:w-algebras}, the different definitions of affine W-algebras are recalled. We provide an explicit proof for the equivalence of all these definitions (\Cref{theorem:equivalent-vertex}) using the associated Slodowy slices. Finally, we give a geometric interpretation to the Kac--Roan--Wakimoto embedding from the perspective of reduction by stages (\Cref{proposition:krw-embedding-geometry}).

In \Cref{section:reduction-stages}, we prove reduction by stages for affine W-algebras, that is to say \Cref{main:reduction-by-stages-w-algebras} (\ref{theorem:reduction-by-stages-w-algebras}), by using the analogous result for Slodowy slices. For this, we provide a new definition of the affine W-algebra $\W^k(\g, f_2)$ in \Cref{main:new-brst-complex} (\ref{theorem:new-construction-w-algebra}).

\subsection{Notation and convention} 

Unless otherwise stated, all objects (such as vector spaces, Lie algebras, algebras, schemes, varieties, algebraic groups, vertex algebras...) are defined over the complex number field $\C$. If $V$ is a vector space, $V(\!(z)\!)$ denotes the vector space of Laurent series with coefficients in $V$.

For any affine scheme $X = \Spec R$, the coordinate ring is denoted by $\C[X] \defeq R$. If $f : X \rightarrow Y$ is a homomorphism of affine schemes, the comorphism is denoted by~$f^* : \C[Y] \rightarrow \C[X]$. Any affine variety (reduced affine scheme of finite type) is identified with the set of its $\C$-points. A linear algebraic group is an affine variety with a group scheme structure.

If $G$ is a linear algebraic group, its Lie algebra is denoted by $\g$. Actions of an algebraic group on a scheme or a vector space are assumed to be algebraic, that is to say given by a co-action \cite{milne2017algebraic}. If both $G$ and $\g$ act on a vector space $V$, for any vector $v$ in $V$ set
$$ G^v \defeq \{g \in G ~ | ~ g \cdot v = v \} \quad \text{and} \quad \g^v \defeq \{x \in \g ~ | ~ x \cdot v = 0 \}. $$
Note that $\g^v$ is the Lie algebra of $G^v$.

\subsection{Acknowledgment} 

Both authors are grateful to Anne Moreau and Tomoyuki Arakawa for useful discussions and for sharing the preliminary versions of \cite{arakawa2024arc}. They are grateful to Anne Moreau and Sven M\"oller for the many relevant remarks they made about the present paper. They thank Dra\v{z}en Adamovi\'c, Jethro van Ekeren, David Ridout and Lewis Topley for organizing great conferences. Some part of this work was done while they were visiting the Instituto de Matem\'{a}tica Pura e Aplicada, Rio de Janeiro, Brazil, in March and April 2022, the Centre de Recherches Math\'{e}matiques, Universit\'{e} de Montr\'{e}al, Qu\'{e}bec, Canada in October and November 2022, the Inter-University Center, Dubrovnik, Croatia, in June 2023, and the International Center for Mathematical Sciences, Edinburgh, Scotland, in August 2023. They are grateful to those institutes for their hospitality.

The first named author is deeply grateful to Thomas Creutzig and Andrew Linshaw for valuable comments and suggestions that led him to find problems with the reduction by stages on affine W-algebras. He wishes to thank Jinwei Yang, David Ridout, Zachary Fehily, Vyacheslav Futorny, Ivan Loseu, Takashi Takebe and Hiraku Nakajima for valuable comments and discussions. Some part of this work was done while he was visiting the ICTP South American Institute for Fundamental Research, S\~{a}o Paulo, Brazil in April 2022, Laboratoire de Math\'{e}matiques d'Orsay, France in January 2023, the University of Melbourne, Australia in May 2023, Research Center for Theoretical Physics, Jagna, Bohol and Belmont Hotel, Mactan, Lapu-Lapu city, Cebu, Philippines in July and August 2023, the SUSTech International Center for Mathematics, Shenzhen, China in March 2024, Shanghai Jiao Tong University, China in June 2024, the Beijing Institute of Mathematical Sciences and Applications, Beijing, China in July 2024 and New Uzbekistan University in August 2024. He is grateful to those institutes, universities, schools, and hotels for their hospitality. He is supported by the World Premier International Research Center Initiative (WPI), MEXT, Japan, and JSPS KAKENHI Grant Number JP21K20317 and JP24K16888.

The second named author wants to thank Micha\"el Bulois, Dylan Butson, Maxime Fairon, Justine Fasquel, Reimundo Heluani, and Shigenori Nakatsuka for valuable comments and discussions. He is supported by a public grant from the Fondation Mathématique Jacques Hadamard. 

Both authors are grateful to the anonymous referees for their careful reading and their improvement suggestions, in particular for the proof of \Cref{lemma:poisson-reduction}.

\section{Vertex algebra and associated Poisson geometry} \label{section:vertex-poisson}

% In Subsection \ref{subsection:vertex-algebras}, we introduce what we need about vertex algebras and their gradings by Hamiltonian operators. In Subsection \ref{subsection:poisson-vertex}, we give two fundamental construction of Poisson vertex algebras: by using the Li filtration of a vertex algebra or by considering the arc space of a Poisson variety. In Subsection \ref{subsection:vertex-examples}, we define the vertex superalgebras that we will need in the rest of the paper.  

\subsection{Graded vertex superalgebras} \label{subsection:vertex-algebras} 

For the details about the general theory of vertex algebras, see~\cite{kac1998vertex,frenkel2004vertex,desole2006finite}.

Let $\V$ be a vertex superalgebra over the complex number field $\C$. Denote by $\One$ the vacuum vector, which is a non-zero even vector in $\V$. Denote by $\partial : \V \rightarrow \V$ the translation operator, which is an even operator. The state-field correspondence is an even $\C$-linear which maps $a \otimes b$ in $\V \otimes_\C \V$ to a Laurent series $a(z) b$ in~$\V(\!(z)\!)$. It is denoted by $a(z)b \defeq \sum_{n \in \Z} a_{(n)}b\, z^{-n-1}$ and the coefficient $a_{(n)}b$ is called the $n$-th product of $a$ and $b$.

For $a, b$ in $\V$, the {normally ordered product} is $\NO{a b} \, \defeq \, a_{(-1)} b$. For~$a, b, c$ in $\V$, we adopt the convention: $ \NO{a b c} \, \defeq \, \NO{a (\NO{b c})}$. For $n$ a nonnegative integer, there is the relation $a_{(-n-1)} b = \frac{1}{n!} \NO{\partial(a) b}$. The $\lambda$-bracket is the formal sum defined by~$[a_\lambda b] \defeq \sum_{n=0}^\infty (a_{(n)}b)\frac{\lambda^n}{n!}$ in $V \otimes_\C \C[\lambda]$, it controls the nonassociativity and the noncommutativity.

Let $\{a_i\}_{i=1}^n$ be a set of odd or even vectors in $\V$. Denote by $\leqslant$ the lexicographic order on $\{1, \ldots, n\} \times \Z_{\geqslant 0}$. The vertex algebra $\V$ is said \emph{freely generated} by the \emph{strong generators} $\{a_i\}_{i=1}^n$ if the normally ordered products
\begin{equation*}
    \NO{\partial^{n_1}(a_{i_1}) \cdots \partial^{n_k}(a_{i_k})} \quad \text{for} \quad 
    \begin{cases}
        & k \geqslant 0 \quad \text{is an integer}, \\
        & i_\bullet \in \{1, \ldots, n\}^k \quad \text{and} \quad n_\bullet \in (\Z_{\geqslant 0})^k, \\
        & (i_j, k_j) \leqslant (i_{j+1}, k_{j+1}) \quad \text{if} \quad a_{i_j} \quad \text{is even}, \\
        & (i_j, k_j) < (i_{j+1}, k_{j+1}) \quad \text{if} \quad a_{i_j} \quad \text{is odd},
    \end{cases}
\end{equation*}
form a basis of $\V$ as a vector space. 

A semisimple operator $\Ham$ on $\V$ is called Hamiltonian if it is diagonalizable and satisfies, for any $a$ in $\V$, the relation $ [\Ham, a(z)] = z\partial_z a(z) + (\Ham a)(z)$. If $a$ is an eigenvector for $\Ham$, the corresponding eigenvalue is denoted by $\Delta(a)$. The decomposition in eigenspaces,
$$ \V = \bigoplus_{\Delta \in \C} \V(\Delta), \quad \text{where} \quad \V(\Delta) \defeq \{a \in \V ~  | ~  \Ham a = \Delta a\},  $$
makes $\V$ a $\C$-graded vector space.

\begin{definition} \label{definition:graded-vertex}
    In this paper, we call a \emph{graded vertex algebra} a vertex algebra equipped with a Hamiltonian operator whose eigenvalues lie in the discrete set $\frac{1}{K}\Z$ for some integer $K \geqslant 1$.
\end{definition}

An even vector $L$ in $\V$, with associated field~$L(z) = \sum_{n \in \Z}L_n\, z^{-n-2}$, is called a \emph{conformal vector} if the operator $L_0$ is a Hamiltonian operator on $\V$, if the operator~$L_{-1}$ coincides with the translation operator $ \partial$ and if there is a complex number~$c$, called the \emph{central charge}, such that $$[L_\lambda L] = (\partial + 2\lambda)L + \frac{c}{12}\lambda^3.$$

\subsection{Graded Poisson vertex algebra associated with the Li filtration}

Given a vertex superalgebra~$\V$, the Li filtration $\FilLi^\bullet\V$ is defined in the following way~\cite{li2005abelianizing}. For an integer $p \leqslant 0$, $\FilLi^p \V$ is equal to $\V$. For $p \geqslant 0$, the linear subspace~$\FilLi^p\V$ is spanned by the elements of the form
$$
    \NO{\partial^{n_1}(a_1) \cdots \partial^{n_k}(a_k)} \quad \text{where} \quad 
    \begin{cases}
        & k \geqslant 0, \quad a_\bullet \in \V^k \quad \text{and} \quad n_\bullet \in (\Z_{\geqslant 0})^k, \\
        & \text{so that} \quad \sum_{i=1}^k n_i \geqslant p.
    \end{cases}
$$
It is a decreasing filtration of $\V$ satisfying the following properties:
\begin{align}
    & \notag \FilLi^0 V = V, \quad \partial(\FilLi^p\V)  \subseteq \FilLi^{p+1}\V, \\
    \label{equation:Li-filtration}  & (\FilLi^p\V)_{(n)}(\FilLi^q\V)\subseteq
    \begin{cases}
        \FilLi^{p+q-n-1}\V & \quad  \mathrm{for} \quad n\in\Z,\\
        \FilLi^{p+q-n}\V & \quad \mathrm{for} \quad n\in\Z_{\geqslant 0}.
    \end{cases}
\end{align}

Denote by $[\bullet]_p : \FilLi^p\V \twoheadrightarrow \FilLi^p\V / \FilLi^{p+1}\V$ the canonical projection. It is known that the associated graded space  $\gr^\bullet_\Li \V \defeq \bigoplus_{p \in \Z}\FilLi^p\V/\FilLi^{p+1}\V$ is a differential supercommutative superalgebra defined by
$$ [a]_p \cdot [b]_q \defeq [\NO{a b}]_{p+q},\quad
\partial [a]_p \defeq [\partial a]_{p+1}, $$
where $a$ is in $\FilLi^p\V$ and $b$ is in $\FilLi^q\V$. It has a structure of a Poisson vertex superalgebra whose Poisson $\lambda$-bracket si given by:
$$\{ [a]_p{}_\lambda [b]_q \} = \sum_{n=0}^\infty [a_{(n)}b]_{p+q-n} \frac{\lambda^n}{n!}.$$
That is, $\{\bullet_\lambda\bullet\}$ satisfies the sesquilinearity, left and right Leibniz rules, skewsymmetry, and Jacobi identity. See \cite{dandrea2006three} for these definitions.

\subsection{Arc spaces of Poisson varieties}
Let $X = \Spec R$ be an affine scheme of finite type, where $R \defeq \C[x_1, \dots, x_r] / I $ is the quotient of a polynomial algebra by an ideal $I$ generated by polynomials $F_1, \ldots, F_s$. 

Consider the polynomial algebra~$\C[x_{i, (-n-1)} ~ | ~1 \leqslant i \leqslant r, ~ n \geqslant 0]$ equipped with the differential operator $\partial$ that is defined on generators by 
$$\partial(x_{i, (-n-1)}) \defeq (n+1) x_{i, (-n-2)}. $$
Denote by $I_\infty$ the associative ideal generated by the derivatives $ \partial^n(F_i)$, for all integers~$1 \leqslant i \leqslant r$ and~$n \geqslant 0$. Then the quotient algebra 
$$ R_\infty \defeq \C[x_{i, (-n-1)} ~ | ~1 \leqslant i \leqslant r, ~ n \geqslant 0] / I_\infty $$
is a differential algebra called the \emph{universal differential algebra} generated by $R$ and the \emph{arc space} of $X$ is the affine scheme $\Jet_\infty X \defeq \Spec R_\infty$. Note that there is a natural inclusion of the algebra $R$ in $R_\infty$.

Moreover, by \cite{arakawa2012remark}, if $X$ is an affine Poisson variety, then the differential algebra $\C[\Jet_\infty X]$ is a Poisson vertex algebra whose Poisson $\lambda$-bracket is given for two functions~$F_1, F_2$ in $\C[X] \subseteq \C[\Jet_\infty X]$ by $\{F_1{}_\lambda F_2\} \defeq \{F_1, F_2\}$, where $\{\bullet , \bullet\}$ is the Poisson bracket on $\C[X]$. 

One can also define the algebras of $m$-jets for any nonnegative integer $m$,
$$ R_m \defeq \C[x_{i, (-n-1)} ~ | ~1 \leqslant i \leqslant r, ~ 0 \leqslant n \leqslant m] / I_m, $$
where $I_m$ is the ideal generated by the derivatives $\partial^n(F_i)$ for the integers $ 1 \leqslant i \leqslant r$ and~$0 \leqslant n \leqslant m$ (with the convention $x_{i, (-n-1)} = 0$ for $n > m$). The \emph{$m$-jet space} of~$X$ the affine scheme $ \Jet_m X \defeq \Spec R_m$. For a pair of integers $0 \leqslant m_1 \leqslant m_2 \leqslant \infty$, there is a natural inclusion of algebras~$ R_{m_1} \hookrightarrow R_{m_2}$,
and the colimit is given by the algebra isomorphism 
$$ R_\infty \cong \colim_{m  \geqslant 0} R_m. $$ 
Note that these results have obvious analogues in the superalgebra setting, we omit details.

\subsection{Examples of vertex superalgebras} \label{subsection:vertex-examples}

\subsubsection{Affine vertex algebras} \label{subsubsection:Kac-Moody}

Let $\mathfrak{a}$ be a finite dimensional Lie algebra and denote by $\kappa : \mathfrak{a} \otimes \mathfrak{a} \rightarrow \C$ an $\mathfrak{a}$-invariant symmetric bilinear form. The \emph{universal affine vertex algebra} $\V^\kappa(\mathfrak{a})$ associated with $\mathfrak{a}$ at level $\kappa$ is a vertex algebra freely generated by any basis of $\mathfrak{a}$. One has the $\lambda$-brackets
$$ [a_\lambda b] = [a, b] + \kappa(a,b)\lambda \quad \text{for} \quad a, b \in \mathfrak{a}. $$

Recall that the Lie brackets extends to a Poisson bracket on the symmetric algebra of $\mathfrak{a}$, which is the coordinate ring of the dual space $\mathfrak{a}^*$. There is a Poisson vertex isomorphism $\gr_\Li \V^\kappa(\mathfrak{a}) \cong \C[\Jet_\infty \mathfrak{a}^*].$ The Poisson $\lambda$-bracket is given for $a, b$ in $\mathfrak{a}$ by $\{a_\lambda b\} = [a, b]$. 

There is a Hamiltonian operator $\Ham^\mathfrak{a}$ such that the conformal degree of $a$ in~$\mathfrak{a}$ is~$\Delta(a)=1$. 

Let $\g$ be a simple finite dimensional complex Lie algebra. Denote by~$(\bullet|\bullet)$ the non-degenerate symmetric invariant bilinear form on $\g$ given by 
$$ (\bullet|\bullet) \defeq (2 \mathcal{h}^\vee)^{-1} \kappa_\g, $$ 
where $\mathcal{h}^\vee$ be the dual Coxeter number of $\g$ and $\kappa_\g$ is the Killing form of the Lie algebra $\g$. The \emph{universal affine vertex algebra} $\V^k(\g)$ associated with $\g$ at level $k$ in $\C$ is defined as the vertex algebra $\V^\kappa(\g)$ for $\kappa \defeq k (\bullet|\bullet)$.

According to \cite[Remark 1.23]{desole2006finite}, for any $H$ in a Cartan subalgebra $\h$ of $\g$, the operator $\Ham^\g + (\partial H)_{(1)}$ defines a Hamiltonian operator on $\V^k(\g)$. The conformal degree of $a$ in $\g$ is $\Delta(a)=1 - \delta$ if there is a complex number~$\delta$ such that $[H, a] = \delta a$. 

Let $\n$ be a finite dimensional nilpotent Lie algebra. The universal affine vertex algebra $\V(\n)$ associated with $\n$ is the vertex algebra $\V^\kappa(\n)$ for $\kappa \defeq 0$. 

\subsubsection{Clifford vertex superalgebras} \label{subsubsection:Clifford} Suppose that $\n$ is a finite-dimensional vector space with a basis $\{x_i\}_{i=1}^{n}$ and~$\n^*$ be the dual space with the dual basis $\{\xi_i\}_{i=1}^{n}$. Let~$\F(\n \oplus \n^*)$ be the \emph{Clifford vertex superalgebra} associated with~ $\n \oplus \n^*$, also called the vertex superalgebra of charged fermions in \cite{kac2003quantum, kac2004quantum}. It is freely generated by odd elements $\phi_i, \phi_i^*$, where $1 \leqslant i \leqslant d$. 

If $x = \sum_{i=1}^n c_i x_i$ and $\xi = \sum_{i=1}^n c_i^* \xi_i$ are elements respectively in $\n$ and~$\n^*$ for some coefficients $c_{\bullet}, c_{\bullet}^*$ in $\C^n$, then set $\phi_x \defeq \sum_{i=1}^n c_i \phi_i$ and $\phi^*_\xi \defeq \sum_{i=1}^n c_i^* \phi_i^*$ that are elements in $\F(\n \oplus \n^*)$. There $\lambda$-brackets are given for $x,y$ in $ \n$ and $\xi, \eta$ in $\n^*$ by
$$ [\phi_x{}_\lambda \phi_\xi^*]=\xi(x), \quad [\phi_x{}_\lambda {\phi_y}] = [\phi_\xi^*{}_\lambda \phi_\eta^*]=0. $$

The Clifford vertex superalgebra $\F(\n \oplus \n^*)$ is graded by the so-called \emph{charge},
$$\F^\bullet(\n \oplus \n^*) = \bigoplus_{n\in\Z} \F^n(\n \oplus \n^*), \quad \text{where} \quad \begin{cases}
    &\operatorname{charge}(\One) \defeq 0, \\
    &\operatorname{charge}(\partial^i(\phi_x)) \defeq -1 \quad \text{for} \quad x \in \n, \\
    &\operatorname{charge}(\partial^i(\phi^*_\xi)) \defeq 1 \quad \text{for} \quad \xi \in \n^*
\end{cases} $$    
and $i \geqslant 0$. The charge is additive for the normally ordered products.

The exterior algebra $\Alt(\n \oplus \n^*)$ is a Poisson superalgebra whose bracket is given for~$x,y$ in $\n$ and $\xi, \eta$ in $\n^*$ by 
$$ \{\phi_x, \phi_\xi^*\} = \xi(x) \quad \text{and} \quad \{\phi_x, \phi_y\} =\{\phi_\xi^*, \phi_\eta^*\} =0.$$ Its universal differential superalgebra $\Alt_\infty(\n \oplus \n^*)$ is a Poisson vertex superalgebra and it has an analogue charge grading. There is a Poisson vertex superalgebra isomorphism $\gr_\Li \F^\bullet(\n \oplus \n^*) \cong \Alt^\bullet_\infty(\n \oplus \n^*)$. This isomorphism respects the charge grading.

The element $ L^\F \defeq \sum_{i=1}^{n}\NO{(\partial\phi_i^*)\phi_i} $ defines a conformal vector of $\F(\n \oplus \n^*)$. The elements $H_i^\F(z) \defeq - \NO{ \phi_i^*(z) \phi_i}$
satisfy the relations $[{H_i^\F}_\lambda{H_j^\F}] = \delta_{i=j} \lambda$. More generally, fix a complex number $m$ for each $1 \leqslant i \leqslant n$. Then, according to~\cite[Remark 1.23]{desole2006finite}, the element $$ L^\F(m_\bullet)  \defeq  L^\F  + \sum_{i=1}^{n} m_i \partial H_i^\F $$ defines a conformal vector of $\F(\n \oplus \n^*)$, the corresponding central charge equals to~$c^\F(m_\bullet) = \sum_{i=1}^{n}\left(-12m_i^2+12m_i-2\right)$. The conformal degrees are $\Delta(\phi_i) = 1 - m_i$ and $\Delta(\phi_i^*) = m_i$ for $1 \leqslant i \leqslant d$.

\subsubsection{Weyl vertex algebras} 

Let $(V, \omega)$ be a finite-dimensional symplectic space with a basis $\{v_i\}_{i=1}^{2 s}$. Let $\A(V)$ be the \emph{Weyl vertex algebra} of $V$, it is also called the neutral fermion vertex algebra in \cite{kac2003quantum, kac2004quantum} or $\beta\gamma$-system. It is freely generated by fields $\psi_i$ for $1 \leqslant i \leqslant 2s$. If $v = \sum_{i=1}^{2s} c_i v_i$ is a generic element in $V$ for some coefficients $c_\bullet$ in $\C^d$, set $ \psi_v \defeq \sum_{i=1}^{2 s} c_i\psi_i$ in $\A(V)$.
One has the $\lambda$-brackets
$$ [\psi_{v}{}_\lambda \psi_{w}] = \omega(v, w) \quad \text{for} \quad v, w \in V. $$

There is a Poisson vertex algebra isomorphism $ \gr_\Li \A(V) \cong \C[\Jet_\infty V]$, where the Poisson $\lambda$-bracket on the right-hand side is induced by the symplectic form $\omega$ in the following sense. The symplectic form induces a linear isomorphism $V \cong V^*$ given by $v \mapsto \omega(\bullet, v)$, and the Poisson bracket is given for $v, w$ in $V \hookrightarrow \C[V]$  by~$ \{v, w\} = \omega(v, w)$.

Let $\{v^i\}_{i=1}^{2 s}$ be the dual basis of $V$ with respect to $\omega$, entirely defined by the relations $\omega(v_i, v^j) = \delta_{i = j}$, and set $\psi^i \defeq \psi_{v^i}$. Then the element $ L^\A \defeq \frac{1}{2} \sum_{i=1}^{2 s}\NO{\partial(\psi^i)\psi_i} $ is a conformal vector of $\A(V)$. 

More generally, assume that there is a symplectic isomorphism of $V$ with the cotangent bundle $\mathrm{T}^*W$ of a Lagrangian subspace $W$ of $V$. Assume that $\{v_i\}_{i=1}^{s}$ is a basis of $W$, so $\{v^i\}_{i=1}^{s}$ is a basis of $W^*$. Set $ H^\A_i = -\NO{\psi^i \psi_i}$ for $1 \leqslant i \leqslant s$. Their $\lambda$-brackets are $[H_i^\A{}_\lambda{H_j^\A}] = - \delta_{i=j}\lambda$. For each $1 \leqslant i \leqslant s$, fix a complex number $a_i$. Then, the element 
$$ L^\A(a_\bullet) \defeq L^\A + \sum_{i=1}^{s} a_i \partial(H_i^\A) $$ 
defines a conformal vector of $\A(V)$, and~$c^\A(a_\bullet) = -\frac{1}{2}\dim V + 12\sum_{i=1}^{s} {a_i}^2$ is the corresponding central charge \cite[Remark 1.23]{desole2006finite}. Then the conformal degree of~$\psi_i$ is $\Delta(\psi_i) = \frac{1}{2} - a_i$ and the one of $\psi^i$ is $\Delta(\psi^i) = \frac{1}{2} + a_i$ for $1 \leqslant i \leqslant s$.

\section{Vanishing theorems for BRST cohomology} \label{section:vanishing-brst}

We recall the definition of the BRST (Becchi--Rouet--Stora--Tyutin) cochain complex for three types of objects: for Poisson varieties in \Cref{subsection:brst-poisson}, for their arc spaces in \Cref{subsection:brst-arc} and for vertex algebras in \Cref{subsection:brst-vertex}. \Cref{subsection:filtration-brst} contains the proof of the cohomology vanishing result of \Cref{main:vanishing-vertex}~(\ref{theorem:vanishing-vertex}), generalising ideas from \cite{arakawa2024arc}. For the classical construction of BRST cohomology in these three contexts, see \cite{arakawa2015associated,arakawa2015localization, arakawa2024arc}.

\subsection{Hamiltonian reduction with two groups}

Let $X$ be an affine scheme endowed with an action of an affine algebraic group $N$. There is an action of the group~$N$ on the coordinate ring $\C[X]$ and one can consider the subalgebra of invariant functions 
$$ \C[X]^N = \{F \in \C[X] ~ | ~ \text{for all} ~ g \in N, ~  g \cdot \phi = \phi \}.$$
The \emph{affine GIT quotient} of $X$ modulo the action of $N$ is defined as the affine scheme $X /\!/ N \defeq \Spec \C[X]^N$.  The Lie algebra $\n$ of $N$ acts by derivations on the coordinate ring $\C[X]$ and if $N$ is connected, then
$$ \C[X]^N = \C[X]^\n \defeq \{F \in \C[X] ~ | ~ \text{for all} ~  x \in \n, ~  x \cdot F = 0 \}. $$

Assume that $M, N$ are two connected affine algebraic groups such that $M$ is a normal subgroup of $N$, and denote by $\m, \n$ their respective Lie algebras. Let~$X$ be an affine Poisson variety with an algebraic action of $N$ by Poisson automorphisms. Assume the existence of an~$N$-equivariant map~$\mu : X \rightarrow \m^*$ such that $\mu$ is a moment map for the $M$-action on $X$, that is to say the comorphism~$ \mu^* : \m \rightarrow (\C[X], \{\bullet, \bullet\}) $ is a Lie algebra homomorphism which fits in the following commutative triangle:
$$ \begin{tikzcd}[column sep = small]
    \m \arrow[rr, "\mathrm{action}"] \arrow[rd, "{\mu^*}"'] & & {\operatorname{Der} \C[X]} \\
    & {\C[X]}, \arrow[ru, "{F \mapsto \{F, \bullet\}}"'] &
\end{tikzcd} $$ 
where the Lie algebra of derivations on $\C[X]$ is denoted by $\operatorname{Der} \C[X]$. The quotient~$\mu^{-1}(0) /\!/ N$ is called~\emph{Hamiltonian reduction}.

\begin{lemma} \label{lemma:poisson-reduction}
    The Hamiltonian reduction $\mu^{-1}(0) /\!/ N$ is a Poisson scheme.
\end{lemma}

The statement is well-known when $M = N$, see \cite[Chapter 5]{laurent-gengoux2013poisson}. In this case, any action with a moment map is an action by Poisson automorphisms. 

\begin{proof} 
	By the classical result from \cite[Chapter 5]{laurent-gengoux2013poisson}, the Hamiltonian reduction
	$$ \mu^{-1}(0) /\!/ M = \Spec(\C[\mu^{-1}(0)]^M)  $$
	carries a Poisson structure induced by the one on $X$. Because $M$ is a normal subgroup of $N$, the action of $N$ descends to the Hamiltonian reduction $\mu^{-1}(0) /\!/ M$ and it preserves the Poisson structure because $N$ acts by Poisson automorphisms on $X$.
	
	Therefore the two stage quotient $(\mu^{-1}(0) /\!/ M) /\!/ N$ has an induced Poisson structure. But it is clear that there is an equality of algebras
	$$ (\C[\mu^{-1}(0)]^M)^N = \C[\mu^{-1}(0)]^N $$
	inside the ambient algebra .
	
	Therefore there is a natural scheme isomorphism
	$$  \mu^{-1}(0) /\!/ N \cong (\mu^{-1}(0) /\!/ M) /\!/ N, $$
	and the left-hand side algebra can be equipped with the Poisson structure coming from the right-hand side algebra.
\end{proof}

Denote by~$\Orb$ the $N$-orbit of an element $\chi$ in $\m^*$, and assume that this orbit is a closed subvariety. Denote by~$\mu^{-1}(\Orb)$ the scheme-theoretic fiber, which is an~$N$-invariant affine subscheme of~$X$. The affine scheme $\mu^{-1}(\Orb) /\!/ N$ is Poisson and is called the \emph{Hamiltonian reduction} of~$X$ with respect to the action of~$N$, the moment map~$\mu$ and the orbit $\Orb$. Equivalently, one can consider the product Poisson variety~$X \times \Orb^-$, where~$\Orb^-$ denotes the orbit of~$-\chi$. This Poisson variety is equipped with the diagonal $N$-action and the~$\Orb$-twisted moment map
$$ \mu_\Orb : X \times \Orb^- \longrightarrow \m^*, \quad (x, \xi) \longmapsto \mu(x) + \xi. $$
The natural projection $X \times \Orb^- \rightarrow X$ induces a Poisson isomorphism
$$ {\mu_\Orb}^{-1}(0) /\!/ N \cong \mu^{-1}(\Orb) /\!/ N. $$
So, up to this $\Orb$-twist, we can always assume that the Hamiltonian reduction is realized with respect to the trivial orbit $\{0\}$. This is a very classical construction, see \cite[Section 24.4]{dasilva2001lectures}.

\subsection{BRST cohomology for Poisson varieties} \label{subsection:brst-poisson}

Let $N$ be a connected affine algebraic group acting on a Poisson variety $X$ by Poisson automorphisms, $M$ be a connected closed normal subgroup of $N$ and $\mu : X \rightarrow \m^*$ be an $N$-equivariant moment map. Assume the data of a Poisson variety $\widetilde X$ endowed with an action of the group $N$ and a moment map $\widetilde\mu : \widetilde X \rightarrow \n^*$. In addition, assume the existence of a dominant $N$-equivariant Poisson map $\rho : \widetilde X \rightarrow X$ such that the following square commutes:
\begin{equation} \begin{tikzcd}
    \widetilde X \arrow[d, "\rho"'] \arrow[r, "\tilde\mu"] & \n^* \arrow[d, two heads] \\
    X \arrow[r, "\mu"'] & \m^*.
\end{tikzcd} \label{equation:geometric-square} \end{equation}

The following construction follows \cite{kostant1987symplectic}. Denote by $\{x_i\}_{i \in I(\n)}$ a basis of the vector space $\n$. Consider the graded Poisson superalgebra
$$ \widetilde C^\bullet \defeq \C\big[\widetilde X\big] \otimes_\C \Alt^\bullet(\n \oplus \n^*), $$
with the notations introduced in~\Cref{subsubsection:Clifford}. It contains the following element:
\begin{equation}
    \widetilde Q \defeq \sum_{i \in I(\n)} \widetilde\mu^*(x_i) \phi^i - \frac{1}{2} \sum_{i, j \in I(\n)} \phi_{[x_i, x_j]} \phi^i \phi^j, \label{equation:brst-element}
\end{equation}
which is of charge one and odd. Consider the Poisson adjoint action of $\widetilde Q$, denoted by $\widetilde d \defeq \big\{\widetilde Q, \bullet \big\}$. The pair $\big(\widetilde C^\bullet, \widetilde d \big)$ forms a cochain complex and the associated cohomology~$\Hgy^\bullet\big(\widetilde C^\bullet, \widetilde d \big)$ has a natural graded Poisson superalgebra structure.

By definition, the dominant map $\rho : \widetilde X \rightarrow X$ induces a Poisson algebra embedding $ \C[X] \hookrightarrow \C\big[\widetilde X\big]$. For simplicity, we identify the left-hand side algebra to its image in the right-hand side. Under this identification, we can define the following graded Poisson super-subalgebra of $\widetilde C$,
$$ C^\bullet \defeq \C[X] \otimes_\C \Alt^\bullet(\m \oplus \n^*), $$
where $\Alt^\bullet(\m \oplus \n^*)$ denotes the Poisson super-subalgebra generated by the elements~$ \phi_x$ and $\phi^*_\xi$, where $x$ is in $\m$ and $\xi$ is in $\n^*$. The super-subalgebra $C$ is stable by the differential~$\widetilde d$. 

Denote by $d$ the restriction of the differential $\widetilde d$ to the superalgebra $C$. Then the pair $(C^\bullet, d)$ forms a cochain complex, called the \emph{BRST cochain complex} associated with the Poisson variety $X$, the acting group $N$ and the moment map $\mu : X \rightarrow \m^*$. Its homology~$\Hgy^\bullet\big(\widetilde C^\bullet, \widetilde d \big)$ has a natural graded Poisson superalgebra structure. The following statement is classical, see~\cite[Theorem~7.1]{arakawa2024arc}.

\begin{theorem}
    Let $X$ be a Poisson variety with an action of the group $N$ by Poisson automorphisms and an $N$-equivariant moment map $\mu : X \rightarrow \m^*$ for the action of the normal subgroup $M$ of $N$, and define the corresponding BRST complex~$(C, d)$. Make the following assumptions:
    \begin{enumerate}    
        \item the moment map $\mu : X \rightarrow \m^*$ is smooth,

        \item there exists a closed subvariety $S$ of $\mu^{-1}(0)$ such that the action map
        $$ \alpha : N \times S \longrightarrow \mu^{-1}(0), \quad (g,x) \longmapsto g \cdot x $$
        is an isomorphism.
    \end{enumerate}
    Then the subvariety $S$ is isomorphic to $\mu^{-1}(0) /\!/ N$, so the variety $S$ inherits a Poisson structure. There is a natural Poisson isomorphism
    $$ \Hgy^\bullet(C, d) \cong \Hgy^\bullet(\n, \C[N]) \otimes_\C \C[S], $$
    where the right-hand side is the tensor product of the trivial Poisson superalgebra $\Hgy^\bullet(\n, \C[N])$ given by the Lie algebra cohomology of the $\n$-module $\C[N]$ (the action is by left-invariant derivations) and the Poisson algebra $\C[S]$.
\end{theorem}

\begin{remark} \label{remark:two-groups}
	In \cite{arakawa2024arc}, a proof is provided for $M = N$ by using a spectral sequence argument : the BRST cohomology is the total cohomology of a double complex mixing Koszul homology and Lie algebra cohomology. This argument also works in the case $M \neq N$.
\end{remark}

\subsection{Affine GIT quotient of an arc space} \label{subsection:git-arc} 

Let $X$ be an affine scheme with an action of a connected affine algebraic group $N$. For any nonnegative integer $m$, there is an induced action of $\Jet_m N$, the $m$-jet space of the group, on $\Jet_m X$, the $m$-jet space of the variety. The scheme $\Jet_m N$ is a connected variety because the group~$N$ is smooth and connected, its Lie algebra is~$\n[t] / (t^{m+1}) \defeq \n \otimes_\C \C[t]/(t^{m+1})$, with the~$\C[t]$-linear Lie bracket extending the one on $\n$. 

By imposing that the ideal $\n[t]t^{m+1}$ acts by zero, we get an action by derivation of the Lie algebra $\n[t]  \defeq \n \otimes_\C \C[t]$ on each coordinate ring $\C[\Jet_m X]$ and the invariant functions for the $\n[t]/(t^{m+1})$ and~$\n[t]$-actions are the same:
$$ \C[\Jet_m X]^{\n[t]/(t^{m+1})} = \C[\Jet_m X]^{\n[t]}. $$
There is also an induced action of $\n[t]$ on the colimit 
$$ \C[\Jet_\infty X] \cong \colim_{m \geqslant 0} \C[\Jet_m X] $$ 
and the subalgebra $\n[t]$-invariant functions of $\Jet_\infty X$ is given by
$$ \C[\Jet_\infty X]^{\n[t]} \cong \colim_{m \geqslant 0} \C[\Jet_m X]^{\n[t]}. $$

The infinite type group scheme $\Jet_\infty N$ acts on the arc space $\Jet_\infty X$ through a coaction $\Phi : \C[\Jet_\infty X] \rightarrow \C[\Jet_\infty N] \otimes_\C \C[\Jet_\infty X]$, and the subalgebra of $\Jet_\infty N$-invariant functions is defined by
$$ \C[\Jet_\infty X]^{\Jet_\infty N} \defeq \{F \in \C[\Jet_\infty X] ~ | ~ \Phi(F) = 1 \otimes F\}. $$
One can show that it is given by the colimit $$ \C[\Jet_\infty X]^{\Jet_\infty N} \cong \colim_{m \geqslant 0} \C[\Jet_m X]^{\Jet_m N}.$$ 
Whence, because the group $N$ is connected, the $\Jet_\infty N$ and $\n[t]$-invariants coincide:
$$ \C[\Jet_\infty X]^{\Jet_\infty N} = \C[\Jet_\infty X]^{\n[t]} . $$

For more details about group schemes and coactions, see for instance \cite{milne2017algebraic}.

\subsection{BRST cohomology for arc spaces} \label{subsection:brst-arc} 

Let $\rho : \widetilde X \rightarrow X$ be a $N$-equivariant dominant Poisson map between two affine Poisson varieties $\widetilde X, X$ equipped with $N$-actions by Poisson automorphisms. Let $\widetilde \mu : \widetilde X \rightarrow \n^*$ and $\mu : X \rightarrow \m^*$ be two $N$-equivariant moment maps such that \eqref{equation:geometric-square} commutes. 

The Poisson superalgebra $\widetilde C^\bullet = \C\big[\widetilde X\big] \otimes_\C \Alt^\bullet(\n \oplus \n^*) $
induces a Poisson vertex superalgebra
$$ \widetilde C^\bullet_\infty \defeq \C\big[\Jet_\infty\widetilde X\big] \otimes_\C \Alt^\bullet_\infty(\n \oplus \n^*), $$
with the notations introduced in \Cref{subsubsection:Clifford}. It contains the element $\widetilde Q$, introduced in \eqref{equation:brst-element}, which is of charge degree one and odd. Consider the adjoint action of~$\widetilde Q$, denoted by $ \widetilde d_\infty \defeq \big\{\widetilde Q_\lambda \bullet \big\}_{\lambda=0}$. The pair $\big(\widetilde C^\bullet_\infty, \widetilde d _\infty\big)$ forms a cochain complex and its associated cohomology~$\Hgy^\bullet\big(\widetilde C^\bullet_\infty, \widetilde d_\infty \big)$ has a natural graded Poisson vertex superalgebra structure.

By definition, the dominant map $\rho : \widetilde X \rightarrow X$ induces a Poisson vertex algebra embedding $\C[\Jet_\infty X] \hookrightarrow \C\big[\Jet_\infty \widetilde X\big]$, so to simplify, we identify the left-hand side algebra to its image in the right-hand side. Under this identification, we can define the following graded Poisson super-subalgebra of $\widetilde C_\infty$,
$$ C^\bullet_\infty\defeq \C[\Jet_\infty X] \otimes_\C \Alt^\bullet_\infty(\m \oplus \n^*), $$
where $\Alt^\bullet_\infty(\m \oplus \n^*)$ denotes the Poisson vertex super-subalgebra generated by elements~$ \partial^n \phi_x$ and $\partial^n \phi^*_\xi$ for $x$ in $\m$, $\xi$ in $\n^*$ and $n$ a nonnegative integer. The super-subalgebra $C_\infty$ is stable by the differential $\widetilde d_\infty$. 
 
Denote by $d_\infty$ the restriction of the differential $\widetilde d_\infty$ to the superalgebra $C_\infty$. Then the pair $(C_\infty, d_\infty)$ forms a cochain complex, called the \emph{Poisson vertex BRST cochain complex} associated with the arc space of the Poisson variety $X$ with its action of $N$ and the moment map $\mu : X \rightarrow \m^*$. Its cohomology~$\Hgy^\bullet(C_\infty, d_\infty)$ has a natural graded Poisson superalgebra structure. The following statement is classical, see~\cite[Theorem~9.2]{arakawa2024arc}.

\begin{theorem} \label{theorem:vanishing-arc}
    Let $X$ be a Poisson variety with an action of the group $N$ by Poisson automorphisms and an $N$-equivariant moment map $\mu : X \rightarrow \m^*$ for the action of the normal subgroup $M$ of $N$, and define the corresponding BRST complex $(C_\infty, d_\infty)$. Make the following assumptions:
    \begin{enumerate}    
        \item the moment map $\mu : X \rightarrow \m^*$ is smooth,

        \item there exists a closed subvariety $S$ of $\mu^{-1}(0)$ such that the action map
        $$ \alpha : N \times S \longrightarrow \mu^{-1}(0), \quad (g,x) \longmapsto g \cdot x $$
        is an isomorphism.
    \end{enumerate}
    Then there is a natural Poisson vertex isomorphism
    $$ \Hgy^\bullet(C_\infty, d_\infty) \cong \Hgy^\bullet(\n[t], \C[\Jet_\infty N]) \otimes_\C \C[\Jet_\infty S], $$
    where the right-hand side is the tensor product of the trivial Poisson vertex superalgebra~$\Hgy^\bullet(\n[t], \C[\Jet_\infty N])$ given by the cohomology of the $\n[t]$-module~$\C[\Jet_\infty N]$ (the action is by left-invariant derivations) and the Poisson vertex algebra~$\C[\Jet_\infty S]$.
\end{theorem}

\begin{remark}
	\Cref{theorem:vanishing-arc} is slightly different than the statement in~\cite{arakawa2024arc}. First, we deal with the case when $M$ and $N$ are not necessarily equal, see \Cref{remark:two-groups}.
	
	Moreover, in~\cite{arakawa2024arc}, the hypotheses of the theorem are less restrictive: the moment maps $\mu : X \rightarrow \m^*$ and $\Jet_\infty \mu : \Jet_\infty X \rightarrow \Jet_\infty \m^*$ are assumed to be flat. We claim that if $\mu$ is smooth, then $\Jet_\infty \mu$ is automatically flat.	Indeed, if~$\mu$ is smooth, then $\Jet_m \mu : \Jet_m X \rightarrow \Jet_m \m^*$ is flat for any integer $m \geqslant 0$ \cite[Proposition~3.7.4]{chambert2018motivic}. By using the equational criterion for flatness \cite[Corollary 6.5]{eisenbud1995commutative}, we deduce that~$\Jet_\infty \mu$ is flat too.
\end{remark}

\begin{remark} 
    For any nonnegative integer $m$, the Lie algebra~$\n[t] / (t^{m+1})$ is finite dimensional and the Lie algebra cohomology of the module~$\C[\Jet_m N]$ is defined as the cohomology of the associated Chevalley--Eilenberg complex $C_{\mathrm{CE}, m} ^\bullet$, see~\cite[Section 7.7]{weibel1994introduction} for details. But for $\n[t]$, which is infinite-dimensional, the cohomology~$\Hgy^\bullet(\n[t], \C[\Jet_\infty N])$ shall be understood as the cohomology of the complex $C_{\mathrm{CE}, \infty}^\bullet \defeq \colim_{m \geqslant 0} C_{\mathrm{CE}, m}^\bullet$.
\end{remark}

\subsection{BRST cohomology for vertex algebras} \label{subsection:brst-vertex}

Let $M, N$ be two connected affine algebraic groups such that $M$ is a normal subgroup of $N$, and denote by $\m, \n$ their respective Lie algebras. We assume from now that the groups $M$ and $N$ are \emph{unipotent}, in particular the Killing forms of their Lie algebras are zero. Let $\V, \widetilde \V$ be two vertex algebras such that $\V$ is a subalgebra of $\widetilde \V$. 

Let $\V(\m)$, $\V(\n)$ be the universal affine vertex algebras associated respectively with~$\m$, $\n$. Let $\widetilde\Upsilon : \V(\n) \rightarrow \widetilde \V$ be a vertex algebra homomorphism, it induces an $\n[t]$-action on $\widetilde \V$. We assume that the vertex subalgebra~$\V$ is a $\n[t]$-submodule of~$\widetilde\V$. Moreover, we assume the inclusion $\widetilde\Upsilon(\V(\m)) \subseteq \V$. Whence, by restriction of~$\widetilde\Upsilon$, we get a vertex algebra map $\Upsilon : \V(\m) \rightarrow \V$ such that the following square is commutative:
\begin{equation} \begin{tikzcd}[column sep = large]
    \V(\n) \arrow[r, " \widetilde\Upsilon"] &  \widetilde\V  \\
    \V(\m) \arrow[r, "\Upsilon"'] \arrow[u, hook] &  \V \arrow[u, hook].
\end{tikzcd} \label{equation:chiral-comoment-square} \end{equation}
In particular, the vertex algebra $\V$ is a module over~$\n[t] \oplus \m[t^{-1}]$. The homomorphism $\Upsilon : \V(\m) \rightarrow \V$ is called a \emph{chiral comoment map}.

Let $\{x_i\}_{i \in I(\n)}$ be a basis of $\n$. Let $\F^\bullet(\m \oplus \n^*)$ be the graded vertex subalgebra of the Clifford vertex algebra $\F^\bullet(\n \oplus \n^*)$ strongly generated by elements $\phi_x$ and~$\phi^*_\xi$, where $x$ is in $\m$ and $\xi$ is in $\n^*$. Consider the charge-graded vertex superalgebra
$$ \widetilde \Cpx^\bullet \defeq \widetilde\V \otimes_\C \F^\bullet(\n \oplus \n^*). $$
It contains the element
$$ \widetilde Q \defeq \sum_{i \in I(\n)} \NO{\widetilde\Upsilon(x_i) \phi^*_i} - \frac{1}{2}\sum_{i, j \in I(\n)} \NO{\phi_{[x_i, x_j]}\phi^*_i \phi^*_j} $$
which is of charge one and odd. Consider the operator given by 0-th mode of $\widetilde Q$, denoted by $\widetilde{\mathcal{d}} \defeq \widetilde Q_{(0)} = \big[\widetilde Q \,_\lambda \bullet \big]_{\lambda = 0}$. The pair $\big(\widetilde \Cpx^\bullet, \widetilde{\mathcal{d}}\big)$ forms a cochain complex and its associated cohomology~$\Hgy^\bullet\big(\widetilde \Cpx^\bullet, \widetilde{\mathcal{d}}\big)$ has a natural graded vertex superalgebra structure.

Define the following graded vertex super-subalgebra
$$ \Cpx^\bullet \defeq \V \otimes_\C \F^\bullet(\m \oplus \n^*) $$
of the cochain complex $\widetilde \Cpx$. 

\begin{lemma}
    The super-subalgebra $\Cpx$ is stable by the differential $\widetilde{\mathcal{d}}$. 
\end{lemma}

\begin{proof}
    Using the inclusions~$[\n,\m] \subseteq \m$ and $\widetilde\Upsilon(\V(\m)) \subseteq \V$, and the fact that $\V$ is a~$\n[t]$-submodule of $\widetilde\V$, the lemma follows from a computation on generators.
\end{proof}

Denote by $\mathcal{d}$ the restriction of the differential $\widetilde{\mathcal{d}}$ to the vertex superalgebra $\Cpx$. Then the pair $(\Cpx^\bullet, \mathcal{d})$ forms a cochain complex, called the \emph{BRST cochain complex} associated with the vertex algebra $\V$ equipped with an action of $\n[t]$ and the chiral comoment map $\Upsilon : \V(\m) \rightarrow \V$.

\subsection{Non-negatively graded quotient complex} \label{subsection:quotient}

Assume that the vertex algebra $\widetilde \V$ is graded in the sense of \Cref{definition:graded-vertex} and that $\V$ is a graded subalgebra. For any element $x \in \widetilde{\V}$, denote by $x_\Delta$ the image of $x$ through the natural projection~$\widetilde{\V} \twoheadrightarrow \widetilde{\V}(\Delta)$. 

\begin{lemma}\label{lemma:standard}
    Suppose that the grading of $\widetilde{\V}$ and $\V$ lies in $\frac{1}{K}\Z_{\geqslant 0}$. Assume that, for any element $x$ in $\n$, the image~$\widetilde\Upsilon(x)$ lies in~$\bigoplus_{\Delta \leqslant 1}\widetilde{\V}(\Delta)$. 
    
    Then we have the following $\lambda$-bracket for all $x,y$ in $\n$:
    $$ \big[\widetilde\Upsilon(x)_1{}_\lambda\widetilde\Upsilon(y)_1\big] = \widetilde\Upsilon([x,y])_1. $$
\end{lemma}

\begin{proof}
    Set $X \defeq \widetilde\Upsilon(x)$ and $Y \defeq \widetilde\Upsilon(y)$ in $\widetilde\V$. They decompose as $$ X = \sum_{i=0}^{\infty} X_{1 - \frac{i}{K}} \quad \text{and} \quad Y = \sum_{j=0}^{\infty} X_{1 - \frac{j}{K}}.$$ By usual properties of vertex algebra grading, one has 
    $$ \Delta\big(X_{1 - \frac{i}{K}}{}_{(n)}Y_{1 - \frac{j}{K}}\big) = 1 - \frac{i+j}{K} - n $$
    that equals $1$ if and only if $i=j=n=0$. 
    
    Set $Z \defeq \widetilde\Upsilon([x, y])$. Because $[x{}_\lambda y] = [x,y]$, one has
    $$ Z =  \widetilde\Upsilon([x{}_\lambda y]) = [X{}_\lambda Y] = X_{(0)} Y, $$
    the last equality follows since $Z$ does not contain terms of strictly positive degree in the formal variable $\lambda$. 
    
    Comparing the degrees, we get the desired equality: $Z_1 = X_1{}_{(0)} Y_1$. 
\end{proof}

When the hypotheses of \Cref{lemma:standard} hold, there is a vertex algebra homomorphism~$\widetilde\Upsilon_{\mathrm{st}} : \V(\n) \rightarrow \widetilde\V$ such that $\widetilde\Upsilon_{\mathrm{st}}(x) \defeq \widetilde\Upsilon(x)_1$ for $x$ in $\n$. This homomorphism is called the \emph{standard comoment map} associated with $\Upsilon$. Set
$$ \widetilde Q_{\mathrm{st}} \defeq \sum_{i \in I(\n)} \NO{\widetilde\Upsilon_{\mathrm{st}}(x_i) \phi^*_i} - \frac{1}{2}\sum_{i, j \in I(\n)} \NO{\phi_{[x_i, x_j]}\phi^*_i \phi^*_j}  \quad \text{and} \quad \widetilde{\mathcal{d}}_{\mathrm{st}} \defeq \widetilde{Q}_{\mathrm{st}}{}_{(0)}. $$
Denote by $\mathcal{d}_{\mathrm{st}}$ the restriction of the differential $\widetilde{\mathcal{d}}_{\mathrm{st}}$ to the vertex superalgebra $\Cpx$. Then, the pair $(\Cpx^\bullet, \mathcal{d}_{\mathrm{st}})$ forms a cochain complex, called the \emph{standard BRST cochain complex} by analogy with \cite{feigin1990quantization, kac2003quantum}.

Consider the linear subspace 
$$  \I \defeq  \Span_\C\{\phi_x{}_{(-n-1)}c,  ~ (\mathcal{d}\phi_x){}_{(-n-1)}c ~ | ~ n \geqslant 0, ~ x \in \m, ~ c \in \Cpx\} $$
of the BRST complex $\Cpx$. This subspace is closed by $\mathcal{d}$ and so it descends to the quotient $ \Cpx^\bullet_+ \defeq \Cpx^\bullet / \I^\bullet$. Denote by $\mathcal{d}_+$ the induced differential. Hence the pair~$(\Cpx^\bullet_+, \mathcal{d}_+)$ forms a cochain complex of vector spaces, \emph{with no vertex algebra structure}.

The following theorem is a generalisation of \cite[Proposition 9.3]{arakawa2024arc}.

\begin{theorem}\label{theorem:quotient-quasi-isomorphic}
    Let $N$ be a unipotent affine algebraic group and $M$ be a unipotent normal subgroup. Let $\widetilde\V$ be a nonnegatively graded vertex algebra and $\V$ be a graded vertex subalgebra. Let $\widetilde\Upsilon : \V(\n) \rightarrow \widetilde\V$ be a vertex algebra homomorphism which restricts to a homomorphism $\Upsilon : \V(\m) \rightarrow \V$. Assume that $\V$ is a $\n[t]$-submodule of~$\widetilde\V$. Define the associated BRST cochain complex $(\Cpx^\bullet, \mathcal{d})$ and the quotient complex~$(\Cpx^\bullet_+, \mathcal{d}_+)$ as above. Assume the following conditions: 
    \begin{enumerate}
        \item for all $\Delta$ in $\frac{1}{K}\Z_{\geqslant 0}$, the homogenous subspace $\widetilde\V(\Delta)$ is finite-dimensional,

        \item for any element $x$ in $\n$, the image $\widetilde\Upsilon(x)$ lies in~$\bigoplus_{\Delta \leqslant 1}\widetilde{\V}(\Delta)$,

        \item the standard chiral comoment map $\widetilde\Upsilon_{\mathrm{st}} : \V(\n) \rightarrow \widetilde\V$ defined by \Cref{lemma:standard} induces a free action of the enveloping algebra $\mathcal{U}(\m[t^{-1}]t^{-1})$ on $\V$. 
    \end{enumerate}

    Then there is an isomorphism $\Hgy^\bullet(\Cpx^\bullet, \mathcal{d}) \cong \Hgy^\bullet(\Cpx^\bullet_+, \mathcal{d}_+)$ of vector spaces induced by the canonical projection $\Cpx \twoheadrightarrow \Cpx_+$.
\end{theorem}

To prove this theorem, extend the Hamiltonian operator $\Ham$ of $\widetilde\V$ to the BRST complex $\widetilde\Cpx = \widetilde\V \otimes_\C \F(\n \oplus \n^*)$ by adding the conformal vector~$L^\F_0$. Then $\Ham(\phi_x) = 1$ and $\Ham(\phi^*_\xi) = 0$ for $x$ in $\n$ and $\xi$ in $\n^*$. This grading descends to the subalgebra~$\Cpx$ and we denote the homogeneous subspaces by
$$ \Cpx(\Delta) \defeq \{ c \in \Cpx ~ | ~ \Ham(c) = \Delta c\} \quad \text{for} \quad \Delta \in \frac{1}{K} \Z_{\geqslant 0}. $$ 
One gets an induced decreasing filtration on $\Cpx$ defined by
$$ \Fil_{\Ham}^p\Cpx \defeq \bigoplus_{\Delta \leqslant -p/K} \Cpx(\Delta) \quad \text{for} \quad p \in \Z. $$
The filtration is exhaustive and bounded from below, that is to say $$ \Cpx = \bigcup_{p\in\Z}\Fil_{\Ham}^p\Cpx \quad \text{and} \quad \Fil_{\Ham}^p\Cpx_\V = 0 \quad \text{for} \quad p > 0. $$

This filtration is preserved by the coboundary operator $\mathcal{d}$ by the second assumption of \Cref{theorem:quotient-quasi-isomorphic}. Therefore there is a spectral sequence $\{(\E_{r}, \mathcal{d}_{r})\}_{r=0}^\infty$ associated with the filtered cochain complex $(\Fil_{\Ham}^\bullet\Cpx^\bullet, \mathcal{d})$ \cite[Section 5.4]{weibel1994introduction}, and it converges to $\Hgy^\bullet(\Cpx, \mathbcal{d})$ \cite[Theorem 5.5.1]{weibel1994introduction}.

The $0$-th coincides with the standard BRST complex introduced after \Cref{lemma:standard}:
$$\E_0^{p, q} = \Cpx^{p+q}\Big(\frac{1}{K}p\Big), \quad \mathcal{d}_0 = \mathcal{d}_{\mathrm{st}}, $$
for $p, q \in \Z$. In particular, the standard differential respects the conformal grading. Whence, the first page is:
$$ \E_{1}^{p, q} = \Hgy^{p + q}\bigg(\Cpx^\bullet\Big(\frac{1}{K}p\Big), \mathcal{d}_{\mathrm{st}}\bigg). $$

The standard comoment map $\Upsilon_{\mathrm{st}} : \V(\m) \rightarrow \V$ induces an action of the universal enveloping algebra $\mathcal{U}(\m[t^{-1}]t^{-1})$ of the Lie algebra $\m[t^{-1}]t^{-1}$, denoted by $\cdot_{\mathrm{st}}$.

\begin{lemma} \label{lemma:Hochschild--Serre}
    The 1-st page $\E_{1} = \Hgy^\bullet(\Cpx^\bullet, \mathcal{d}_{\mathrm{st}})$ is isomorphic to Lie algebra cohomology~$\Hgy^\bullet\Big(\n[t], \V/\big(\mathcal{U}(\m[t^{-1}]t^{-1})\cdot_{\mathrm{st}}\V\big)\Big)$.
\end{lemma}

\begin{proof}
    Apply the co-analog of the Hochschild--Serre spectral sequence \cite[Theorem 2.3]{voronov1993semi} on $(\Cpx, \mathcal{d}_{\mathrm{st}})$. It corresponds to the filtration defined by
    $$ \HS^p\Cpx^n \defeq \Span_\C\Big\{\phi^*_{\xi_1}{}_{(-m_1 - 1)} \cdots \phi^*_{\xi_p}{}_{(-m_p - 1)} c \quad \Big| \quad
    \begin{aligned}
        & c \in \Cpx^{n-p} \\
        & \xi_\bullet \in (\n^*)^p, ~  m_\bullet \in (\Z_{\geqslant 0})^p
    \end{aligned}
    \Big.\Big\}, $$
    where $n, p$ are in $\Z$. This filtration is decreasing, bounded from above (it means that $\HS^0\Cpx^n = \Cpx^n$), and it is preserved by the standard differential. There is an induced spectral sequence denoted by $\{(\E_{\HS, r}, \mathcal{d}_{\HS, r})\}_{r=0}^\infty$.

    The standard BRST complex splits into a direct sum of finite-dimensional subcomplexes because the Hamiltonian grading is preserved by the standard differential and the dimension of each summand is finite as a consequence of the first assumption of Theorem \ref{theorem:quotient-quasi-isomorphic}. For any $\Delta$ in $\frac{1}{K} \Z_{\geqslant 0}$, the induced Hochschild--Serre filtration on $\Cpx(\Delta)$ is finite and the associated spectral sequence $\{(\E_{\HS, r}(\Delta), \mathcal{d}_{\HS, r})\}_{r=0}^\infty$ is convergent \cite[Theorem 5.5.1]{weibel1994introduction}.

    The $0$-th is given by the following vector space isomorphism
    $$ \E_{\HS, 0}^{p,q} = \gr_\HS^p\Cpx^{p+q}  \cong \V \otimes_\C \Alt^{-q}_\infty(\m) \otimes_\C \Alt^p_\infty(\n^*), $$
    for $p, q$ in $\Z$. The first page is given by
    $$ \E_{\HS,1}^{p,q} = \Hgy_{-q}(\m[t^{-1}]t^{-1},\V) \otimes_\C \Alt^p_\infty(\n^*), $$
    where $\Hgy_\bullet(\m[t^{-1}]t^{-1},\V)$ is the Lie algebra homology with coefficients in the vector space $\V$ equipped $\m[t^{-1}]t^{-1}$-module induced by $\Upsilon_{\mathrm{st}}$. The action of $\mathcal{U}(\m[t^{-1}]t^{-1})$ on $\V$ is free by assumption. Hence, this homology is given by
    $$ \Hgy_\bullet (\m[t^{-1}]t^{-1},\V) = \delta_{\bullet = 0} \, \Hgy_{0}(\m[t^{-1}]t^{-1},\V) = \V/\big(\m[t^{-1}]t^{-1}\cdot_{\mathrm{st}}\V\big). $$

    Therefore, the spectral sequence collapses at the 2-nd page:
    $$ \E_{\HS,2}^{p,q} = \delta_{q = 0} \, \E_{\HS,2}^{p,0} = \delta_{q = 0} \, \Hgy^p\Big(\n[t], \V/\big(\m[t^{-1}]t^{-1}\cdot_{\mathrm{st}}\V\big)\Big). $$
    The spectral sequence is convergent on each homogeneous component (for the Hamiltonian grading), so the infinity term of the spectral sequence is
    $$ \E_{\HS,\infty}^{p,q} = \delta_{q = 0} \, \gr_\HS^p\Hgy^{p}(\Cpx^\bullet, \mathcal{d}_{\mathrm{st}}). $$ 
    
    The convergence also implies that the filtration on $\Hgy^{\bullet}(\Cpx^\bullet, \mathcal{d}_{\mathrm{st}})$ is complete on each homogeneous component and then the collapsing implies the equality
    $$ \Hgy^{p}(\Cpx^\bullet, \mathcal{d}_{\mathrm{st}}) = \gr_\HS^p\Hgy^{p}(\Cpx^\bullet, \mathcal{d}_{\mathrm{st}}), $$
    and it also implies the isomorphism $\E_{\HS,\infty}^{p,q} \cong \E_{\HS,2}^{p,q}$.
    Finally, we get the desired isomorphism for all $p \geqslant 0$:
    $$ \Hgy^p(\Cpx^\bullet, \mathcal{d}_{\mathrm{st}}) \cong \Hgy^p\Big(\n[t], \V/\big(\m[t^{-1}]t^{-1}\cdot_{\mathrm{st}}\V\big)\Big). $$
\end{proof}

\begin{proof}[Proof of Theorem \ref{theorem:quotient-quasi-isomorphic}]
    Consider the filtration induced on the quotient complex:
    $$ \Fil_{\Ham}^p\Cpx_+ \defeq \Fil_{\Ham}^p\Cpx /(\Fil_{\Ham}^p\Cpx \cap \I) \quad \text{for} \quad p \in \Z. $$
    Then there is a spectral sequence $\{(\E_{+,r}, \mathcal{d}_{+,r})\}_{r=0}^\infty$ converging to $\Hgy^\bullet(\Cpx_{+}, \mathcal{d}_+)$.

    The graded ideal associated with $\I$ is its standard analogue:
    $$  \gr_{\Fil_{\Ham}} \I = \I_{\mathrm{st}} \defeq  \Span_\C\{\phi_x{}_{(-n-1)}c, ~ \mathcal{d}_{\mathrm{st}}(\phi_x){}_{(-n-1)}c ~ | ~ n \geqslant 0, ~ x \in \m, ~ c \in \Cpx\}, $$
    which is graded by the Hamiltonian grading. The 0-th page of the spectral sequence coincides with the corresponding quotient of the standard BRST complex:
    $$\E_{+,0}^{p, q} = \Cpx^{p+q}\Big(\frac{1}{K}p\Big) / \I_{\mathrm{st}} \Big(\frac{1}{K}p\Big), \quad \mathcal{d}_{+,0} = \mathcal{d}_{\mathrm{st}, +}, \quad \text{for} \quad p, q \in \Z. $$
    
    In particular, it coincides with the following Lie algebra cohomology complex:
    $$ \E_{+,0}^{p} =  \V/\big(\m[t^{-1}]t^{-1}\cdot_{\mathrm{st}}\V\big) \otimes_\C \Alt^p_\infty(\n^*). $$
    So the first page is the following Lie algebra cohomology:
    $$ \E_{+,1}^{p} = \Hgy^p\Big(\n[t], \V/\big(\m[t^{-1}]t^{-1}\cdot_{\mathrm{st}}\V\big)\Big). $$

    By \Cref{lemma:Hochschild--Serre}, the filtered complex map $\Fil_{\Ham}^\bullet \Cpx^\bullet \twoheadrightarrow \Fil_{\Ham}^\bullet \Cpx_+^\bullet$ induces a homomorphism of the corresponding convergent spectral sequences which is an isomorphism at the first page:
    $$ \E_1^{p} \cong \E_{+,1}^{p} \cong \Hgy^p\Big(\n[t], \V/\big(\m[t^{-1}]t^{-1}\cdot_{\mathrm{st}}\V\big)\Big). $$
    Hence, by the comparison theorem \cite[Theorem 5.2.12]{weibel1994introduction}, the cohomologies to which the spectral sequences converge are isomorphic:
    $$\Hgy^p(\Cpx^\bullet, \mathcal{d}) \cong \Hgy^p(\Cpx^\bullet_+, \mathcal{d}_+). $$
\end{proof}

\subsection{Induced Li filtration on BRST cohomology} \label{subsection:filtration-brst}

 Let $\V, \widetilde \V$ be two vertex algebras such that $\V$ is a subalgebra of $\widetilde \V$. Let $\widetilde\Upsilon : \V(\n) \rightarrow \widetilde \V$ be a vertex algebra homomorphism that the vertex subalgebra~$\V$ is a $\n[t]$-submodule of~$\widetilde\V$ and the inclusion $\widetilde\Upsilon(\V(\m)) \subseteq \V$ holds. In particular we get the commutative diagram \eqref{equation:chiral-comoment-square}.

Let $\rho : \widetilde X \rightarrow X$ be a Poisson $N$-equivariant dominant map between two affine Poisson varieties $\widetilde X, X$ equipped with $N$-actions by Poisson automorphisms. We assume that there is a commutative square of Poisson vertex algebra homomorphisms,
$$ \begin{tikzcd}
    \gr_\Li \widetilde \V \arrow[r, "\sim"] & \C[\Jet_\infty \widetilde X]  \\
    \gr_\Li \V \arrow[r, "\sim"'] \arrow[u] & \C[\Jet_\infty X], \arrow[u, hook]
\end{tikzcd} $$
where the horizontal arrows are isomorphisms.

Let $\widetilde \mu : \widetilde X \rightarrow \n^*$ and $\mu : X \rightarrow \m^*$ be two $N$-equivariant moment maps such that the diagram \eqref{equation:geometric-square} commutes.
One gets the commutative square
$$ \begin{tikzcd}[column sep = large]
    \C[\Jet_\infty \n^*] \arrow[r, "\Jet_\infty \widetilde\mu^*"] & \C[\Jet_\infty \widetilde X]  \\
    \C[\Jet_\infty \m^*] \arrow[r, "\Jet_\infty \mu^*"'] \arrow[u, hook] & \C[\Jet_\infty X], \arrow[u, "\Jet_\infty\rho^*"', hook]
\end{tikzcd} $$
and we assume that this diagram coincides with the square induced by \eqref{equation:chiral-comoment-square}:
$$ \begin{tikzcd}[column sep = large]
    \gr_\Li \V(\n) \arrow[r, "\gr_\Li \widetilde\Upsilon"] &  \gr_\Li \widetilde\V  \\
    \gr_\Li \V(\m) \arrow[r, "\gr_\Li \Upsilon"'] \arrow[u, hook] & \gr_\Li \V \arrow[u, hook].
\end{tikzcd} $$

As before, introduce BRST cochain complexes for the vertex algebras,
$$  \Cpx^\bullet \defeq \V \otimes_\C \F^\bullet(\m \oplus \n^*) \subseteq \widetilde \Cpx^\bullet \defeq \widetilde\V \otimes_\C \F^\bullet(\n \oplus \n^*), $$
and the BRST cochain complexes for the arc spaces,
$$ C^\bullet_\infty \defeq \C[\Jet_\infty X] \otimes_\C \Alt^\bullet_\infty(\m \oplus \n^*) \\ \subseteq \widetilde C^\bullet_\infty \defeq \C\big[\Jet_\infty\widetilde X\big] \otimes_\C \Alt^\bullet_\infty(\n \oplus \n^*). $$
Clearly, one gets the following Poisson vertex superalgebra isomorphisms, compatible with the coboundary operators:
$$ \gr_\Li \Cpx^\bullet \cong C^\bullet_\infty \quad \text{and} \quad \gr_\Li \widetilde\Cpx^\bullet \cong \widetilde C^\bullet_\infty. $$

The following theorem is a generalisation of \cite[Theorem 9.7]{arakawa2024arc}.

\begin{theorem} \label{theorem:vanishing-vertex}
    Consider the data introduced above. Make the following technical assumptions:
    \begin{enumerate}
        \item the vertex superalgebra $\widetilde\Cpx$ is graded, $\Cpx $ is a graded subalgebra and the element~$\widetilde Q$ defining the coboundary operator $\widetilde{\mathcal{d}}$ is homogeneous of degree $1$,

        \item the space $$  \I \defeq  \Span_\C\{\phi_x{}_{(-n-1)}c, ~ (\mathcal{d}\phi_x){}_{(-n-1)}c ~ | ~ n \geqslant 0, ~ x \in \m, ~ c \in \Cpx\}$$ is a graded subspace of $\Cpx$,

        \item \label{theorem:vanishing-vertex,item:nonnegative} the induced grading on the quotient space $\Cpx_+ \defeq \Cpx / \I$ is nonnegative,

        \item \label{theorem:vanishing-vertex,item:quotient} there is an isomorphism $\Hgy^\bullet(\Cpx^\bullet, \mathcal{d}) \cong \Hgy^\bullet(\Cpx^\bullet_+, \mathcal{d}_+)$ of vector spaces induced by the linear projection $\Cpx^\bullet \twoheadrightarrow \Cpx_+^\bullet$,

        \item \label{theorem:vanishing-vertex,item:section} the moment map $\mu : X \rightarrow \m^*$ is smooth, and there exists a closed subvariety $S$ of $\mu^{-1}(0)$ such that the action map
        $$ \alpha : N \times S \longrightarrow \mu^{-1}(0), \quad (g,x) \longmapsto g \cdot x $$
        is an isomorphism.
    \end{enumerate}

    Then the cohomology vanishes in degrees other than $0$:
    $$ \Hgy^n(\Cpx^\bullet, \mathcal{d}) = 0 \quad \text{for} \quad n \neq 0. $$ 
    In degree 0 there is a natural isomorphism 
    $$ \gr_{\Fil}\Hgy^0(\Cpx^\bullet, \mathcal{d}) \overset{\sim}{\longrightarrow} \Hgy^0(\gr_{\Li} \Cpx^\bullet, \gr_{\Li} \mathcal{d}), $$
    where the filtration $\Fil$ on the cohomology $\Hgy^0(\Cpx^\bullet, \mathcal{d})$ is induced by the Li filtration on the complex $(\Cpx^\bullet, \mathcal{d})$.
\end{theorem}

\begin{remark}
    Our assumption \eqref{theorem:vanishing-vertex,item:nonnegative} is slightly less restrictive than \cite[Theorem~9.7]{arakawa2024arc}, where it is assumed also that $\Cpx^\bullet_+(0) = \C \One$. It will be necessary to prove \Cref{proposition:vanishing-new-complex}.
\end{remark}

\begin{remark}
    In the rest of the paper, we will use \Cref{theorem:quotient-quasi-isomorphic} to get the condition \eqref{theorem:vanishing-vertex,item:quotient} of \Cref{theorem:vanishing-vertex}. In particular, this means that two Hamiltonian operators will be needed: we will usually denote by~$\Ham^\old$ the one used for \Cref{theorem:quotient-quasi-isomorphic} and by~$\Ham^\new$ the one used for \Cref{theorem:vanishing-vertex}.
\end{remark}

Assume the hypotheses of \Cref{theorem:vanishing-vertex}. The Li filtration on $\Cpx$ induces a filtration on the quotient $\Cpx_+$, also denoted by $\FilLi$. The Li filtration is preserved by the coboundary operators $\mathcal{d}$ and $\mathcal{d}_+$ because of the property \eqref{equation:Li-filtration} of this filtration applied to the element $\widetilde Q$. Denote by $\Fil$ the filtration induced on their cohomologies. There are natural maps
\begin{align*}
   \gr_{\Fil}\Hgy^\bullet(\Cpx^\bullet, \mathcal{d})  & \longrightarrow \Hgy^\bullet(\gr_{\Li} \Cpx^\bullet, \gr_{\Li} \mathcal{d}), \\ 
     \gr_{\Fil}\Hgy^\bullet(\Cpx_+^\bullet, \mathcal{d}_+) & \longrightarrow \Hgy^\bullet(\gr_{\Li} \Cpx_+^\bullet, \gr_{\Li} \mathcal{d}_+)
\end{align*}  
which make the following square commutative:
\begin{equation} \begin{tikzcd}
    \gr_{\Fil}\Hgy^\bullet(\Cpx^\bullet, \mathcal{d}) \arrow[d, "\sim"'{sloped}] \arrow[r] & \Hgy^\bullet(\gr_{\Li} \Cpx^\bullet, \gr_{\Li} \mathcal{d}) \arrow[d] \\
    \gr_{\Fil}\Hgy^\bullet(\Cpx_+^\bullet, \mathcal{d}_+) \arrow[r] & \Hgy^\bullet(\gr_{\Li} \Cpx_+^\bullet, \gr_{\Li} \mathcal{d}_+).
\end{tikzcd} \label{equation:Li-quotient} \end{equation}
The left vertical map is an isomorphism by assumption of \Cref{theorem:vanishing-vertex}.

\begin{lemma} \label{lemma:quotient-quasi-isomorphic}
    The natural projection $\gr_\Li \Cpx^\bullet \twoheadrightarrow \gr_\Li \Cpx_+^\bullet$ induces an isomorphism of vector spaces $\Hgy^\bullet(\gr_\Li \Cpx^\bullet, \gr_\Li \mathcal{d}) \cong \Hgy^\bullet(\gr_\Li \Cpx^\bullet_+, \gr_\Li \mathcal{d}_+)$.
\end{lemma}

\begin{proof}
    By assumption, there is an isomorphism of Poisson vertex algebras:
    $$ \gr_\Li \Cpx^\bullet \cong C_\infty = \C[\Jet_\infty X] \otimes_\C \Alt^\bullet_\infty(\m \oplus \n^*) $$
    which maps the boundary operator $\mathcal{d}$ to $d_\infty$. It induces an isomorphism of differential ideals 
    $$ \gr_\Li \I \cong I_{\mu, \infty} \otimes_\C I_{\m, \infty},$$ 
    where $I_{\mu, \infty}$ is the ideal of $\C[\Jet_\infty X]$ corresponding to the closed embedding of schemes $\Jet_\infty \mu^{-1}(0) \hookrightarrow \Jet_\infty X$, and $I_{\m, \infty}$ is the differential ideal of $\Alt^\bullet_\infty(\m \oplus \n^*)$ generated by the elements $\phi_x$, where $x$ is in $\m$.

    Then we get a cochain complex isomorphism 
    $$ \gr_\Li \Cpx_+^\bullet \cong \C[\Jet_\infty \mu^{-1}(0)] \otimes_\C \Alt^\bullet_\infty(\n^*), $$
    where the right-hand side corresponds to the Lie algebra cochain complex of the~$\n[t]$-module $\C[\Jet_\infty \mu^{-1}(0)]$. Then, it follows from the proof of \Cref{theorem:vanishing-arc} that the natural projection $\gr_\Li \Cpx^\bullet \twoheadrightarrow \gr_\Li \Cpx_+^\bullet$ induces an isomorphism of vector spaces between their cohomologies.
\end{proof}

For any nonnegative integer $m$, the Lie algebra cohomology of $\C[\Jet_m N]$ is the algebraic de Rham cohomology of~$\Jet_m N$, and because this group is unipotent, its cohomology is $\C$ in degree $0$ and zero otherwise. Taking the colimit, the Lie algebra cohomology of $\C[\Jet_\infty N]$  is $\C$ in degree $0$ and zero otherwise. Thanks to \Cref{theorem:vanishing-arc} and \Cref{lemma:quotient-quasi-isomorphic}, we deduce that the cohomologies
$$ \Hgy^n(\gr_\Li \Cpx^\bullet, \gr_\Li \mathcal{d}) = \Hgy^n(\gr_\Li \Cpx^\bullet_+, \gr_\Li \mathcal{d}_+) = 0 $$
for any nonzero integer $n$.

\begin{lemma} \label{lemma:li-quotient}
    The cohomology $ \Hgy^n(\Cpx_+^\bullet, \mathcal{d})$ is zero if the integer $n$ is nonzero, and the natural map $\gr_{\Fil}\Hgy^\bullet(\Cpx_+^\bullet, \mathcal{d}_+) \rightarrow \Hgy^\bullet(\gr_{\Li} \Cpx_+^\bullet, \gr_{\Li} \mathcal{d}_+)$ is an isomorphism.
\end{lemma}

\begin{proof}
    Recall that the complex $\Cpx_+$ is nonnegatively graded by the third assumption of Theorem \ref{theorem:vanishing-vertex}. Moreover, by the first assumption, each homogeneous component~$\Cpx_+(\Delta)$, for $\Delta$ in $\frac{1}{K}\Z_{\geqslant 0}$, is a subcomplex. According to \cite[Proposition~2.6.1]{arakawa2012remark}, $\FilLi^p \Cpx_+(\Delta)$ is zero when $p > \Delta$ because the Hamiltonian grading is nonnegative on the quotient complex. So the filtration is finite on each homogeneous component.

    To the filtered complex $(\FilLi^\bullet \Cpx^\bullet, \mathcal{d}_+)$ is associated a spectral sequence, denoted by~$\{(\E_{+, \Li, r}, \mathcal{d}_{+,\Li,r})\}_{r=0}^\infty$, which is convergent on each homogeneous component for the grading induced by the Hamiltonian operator. The first and infinity pages are:
    \begin{align*}
        \E_{+, \Li, 1}^{p, q} & = \Hgy^{p + q}(\gr_{\Li}^p \Cpx_{+}^\bullet, \gr_{\Li} \mathcal{d}_+),\\
        \E_{+, \Li, \infty}^{p, q} & = \gr_{\Fil}^p \Hgy^{p + q}(\Cpx_{+}^\bullet, \mathcal{d}_+),
    \end{align*}
    for $p, q$ in $\Z$. Because of the vanishing implied by the previous lemma and the convergence on homogeneous subcomplexes, we get the desired isomorphism:
    $$ \gr_{\Fil} \Hgy^{n}(\Cpx_{+}^\bullet, \mathcal{d}_+) \cong \Hgy^{n}(\gr_{\Li} \Cpx_{+}^\bullet, \gr_{\Li} \mathcal{d}_+) $$
    for any integer $n$, and both are zero if $n \neq 0$. So we deduce the equality
    $$ \Hgy^n(\Cpx_+^\bullet, \mathcal{d}) = 0 \quad \text{for} \quad n \neq 0. $$
\end{proof}

\begin{proof}[Proof of \Cref{theorem:vanishing-vertex}]
    We go back to the commutative diagram~\eqref{equation:Li-quotient}. The vertical maps and the bottom horizontal map are isomorphisms as a consequence of the previous lemmas. So we can conclude that the top map is the desired isomorphism too. The second lemma and the isomorphism $\Hgy^\bullet(\Cpx^\bullet, \mathcal{d}) \cong \Hgy^\bullet(\Cpx^\bullet_+, \mathcal{d}_+)$ given by the fourth condition implies that $\Hgy^n(\Cpx^\bullet, \mathcal{d})$ is zero when $n$ is nonzero.
\end{proof}

\section{Affine W-algebras} \label{section:w-algebras}

In \Cref{subsection:slodowy-poisson-structure}, we recall the various constructions of the Poisson structure of Slodowy slices by Hamiltonian reduction. By analogy, we define various BRST complexes to construct affine W-algebras in \Cref{subsection:definition-w-algebras},  and prove the equivalence of all these constructions in \Cref{subsection:equivalence-definitions}, see \Cref{theorem:equivalent-constructions}. In \Cref{subsection:structure-w-algebra}, we remind additional facts about the structure of affine W-algebras. In \Cref{subsection:geometric-krw}, we recall the Kac--Roan--Wakimoto embedding and give a geometric interpretation in \Cref{proposition:krw-embedding-geometry}.

\subsection{Good grading} \label{subsection:good-grading}

Let $\g$ be a simple finite-dimensional complex Lie algebra and~$G$ be a connected algebraic group whose Lie algebra is $\g$. Let $\h$ be a Cartan subalgebra of $\g$. Denote by $(\bullet|\bullet)$ the non-degenerate symmetric invariant bilinear form on $\g$ given by
$$ (\bullet|\bullet) \defeq \frac{1}{2 \mathcal{h}^\vee} \kappa_\g, $$
where $\mathcal{h}^\vee$ be the dual Coxeter number of $\g$ and $\kappa_\g$ is the Killing form of the Lie algebra $\g$. Consider a nilpotent orbit $\mathbf O$ in $\g$ and fix an $\sl_2$-triple $(e,h,f)$ in $\g$ such that $f, e$ belong to $\mathbf O$ and $h$ belongs to $\h$.

According to \cite[Lemma 19]{brundan2007good}, there is an element $H$ in $\h$ such that the induced Lie algebra grading
$$ \g = \bigoplus_{\delta \in \Z} \g^{(H)}_\delta, \quad \text{where} \quad \g^{(H)}_\delta \defeq \{x \in \g ~ | ~ [H,x] = \delta x\}, $$
is a \textit{good grading} for~$f$, that is to say: the grading is a $\Z$-grading, $f$ belongs to~$\g^{(H)}_{-2}$ and the map $\ad(f) : \g^{(H)}_\delta \rightarrow \g^{(H)}_{\delta - 2}$ induced by the adjoint action of $f$ is injective for~$\delta \geqslant 1$ and surjective for~$\delta \leqslant 1$. Moreover, we can assume that $e$ belongs to $\g^{(H)}_2$ and $h$ belongs to $\g^{(H)}_0$. 

Fix such a good grading, and to simplify notation, denote $\g_\delta \defeq \g^{(H)}_\delta$. Set $\chi \defeq (f | \bullet)$ the linear form on $\g$ given by the scalar product by $f$. The homogeneous subspace $\g_1$ is equipped with the skewsymmetric form
$$ \omega(v,w) \defeq \chi([v,w]), \quad \text{for} \quad v, w \in \g_1. $$ 
Since the grading is good, this form is symplectic. Let $\l$ be an isotropic subspace of $\g_1$ and denote by $\l^{\perp, \omega}$ its orthogonal set with respect to $\omega$. The quotient $\l^{\perp, \omega} / \l$ is equipped with a skewsymmetric form induced by $\omega$ which is a symplectic form, also denoted by $\omega$.

Define the following nilpotent subalgebra of $\g$:
$$ \n_\l \defeq \l^{\perp, \omega} \oplus \g_{\geqslant 2}. $$
It is the Lie algebra of a unique unipotent subgroup $N_\l$ of $G$. 

By restriction, the linear form $\chi$ restricts to a linear form on $\n_\l$, denoted by $\overline{\chi}_\l$. Denote by $\Orb_\l \defeq \Ad^*(N_\l) \overline{\chi}_\l$ its coadjoint orbit in $\n_\l{}^*$, which is a smooth symplectic variety. The following lemma is stated without proof in \cite[Section 0.4]{desole2006finite}.

\begin{lemma} \label{lemma:symplectic-orbit}
	The map
	$$ \sigma_\l :  \l^{\perp, \omega} / \l \longrightarrow \Orb_\l, \quad (v \bmod \l) \longmapsto \overline{\chi}_\l + \ad^*(v)\overline{\chi}_\l, $$
	is well-defined and is a symplectic isomorphism. In particular, $\Orb_\l$ is an affine subspace of $\n_\l{}^*$.
\end{lemma}

\begin{proof}
    The action of the group $N_\l$ is given by exponentiation of its Lie algebra action, so
    $$  \Ad^*(N_\l) \overline{\chi}_\l = \Big\{ \sum_{k\geqslant 0} \frac{1}{k!} \ad^*(v)^k\overline{\chi}_\l ~ | ~ v \in \n_\l \Big\}. $$
    
    Because of the good grading, $ \ad^*(v)^k\overline{\chi}_\l = 0$ for $k \geqslant 2$. Indeed $\ad(v)^k f$ belongs to~$\g_{k -2} \subseteq \g_{\geqslant 0}$, which is orthogonal to $\n_\l \subseteq \g_{\geqslant 1}$ by the bilinear form $(\bullet|\bullet)$. Hence,
    $$ \Ad^*(N_\l) \overline{\chi}_\l = \{ \overline{\chi}_\l + \ad^*(v)\overline{\chi}_\l ~ | ~ v \in \n_\l\}. $$

    Because of the grading, one can check that $ \ad^*(v)^k\overline{\chi}_\l = 0$ for $v$ in $\g_{\geqslant 2}$. So we get a surjective affine map
    $$  \l^{\perp, \omega} \twoheadlongrightarrow \Orb_\l, \quad v \longmapsto \overline{\chi}_\l + \ad^*(v)\overline{\chi}_\l, $$
    and the fact that $\l$ is isotropic implies the isomorphism $ \l^{\perp, \omega} / \l \cong \Orb_\l$. 
\end{proof}

\subsection{Poisson structure on Slodowy slices} \label{subsection:slodowy-poisson-structure}

The group $N_{\l}$ acts on $\g^*$ by restriction of the coadjoint action, and the restriction map gives a moment map
$$ \pi_\l : \g^* \twoheadlongrightarrow \n_\l{}^*, \quad \xi \longmapsto \xi|_{\n_\l}. $$
Denote by $\Orb_\l^- \defeq - \Ad^*(N_\l) \overline{\chi}_\l$ the opposite coadjoint orbit. The fiber of $\Orb_\l^-$ under~$\pi_\l$ is equal to
$$ \pi_\l{}^{-1}(\Orb_\l^-) = - \chi + \ad^*(\mathfrak{k}) \chi \oplus \n_\l{}^\perp $$ where $\mathfrak{k}$ is a subspace of $\l^{\perp, \omega}$ such that $\l^{\perp, \omega} = \mathfrak{k} \oplus \l$.

The \textit{Slodowy slice} associated with the $\sl_2$-triple $(e,h,f)$ is the affine subspace of $\g^*$ defined as
$$ S_f \defeq - \chi + [\g, e]^\perp, \quad \text{where} \quad [\g, e]^\perp \defeq \{\xi \in \g^* ~  | ~  \xi([\g, e]) = 0\}. $$
The Slodowy slice $S_f$ is contained in this fiber $\pi_\l{}^{-1}(\Orb_\l^-)$.

\begin{remark}
    We work with $-\chi$ instead of $\chi$ to follow the conventions of \cite{kac2003quantum, kac2004quantum, desole2006finite}. Because any nilpotent element belongs to the same coadjoint orbit as its opposite, it does not change the objects up to some isomorphism.
\end{remark}

\begin{theorem}[{\cite[Lemma 2.1]{gan2002quantization}}]\label{theorem:slodowy-isomorphism}
    The coadjoint action of the group $N_\l$ induces an $N_\l$-isomorphism
    $$ \alpha_\l : N_\l \times S_f \overset{\sim}{\longrightarrow} \pi_\l{}^{-1}(\Orb_\l^-), \quad (g, \xi) \longmapsto \Ad^*(g)\xi, $$
    where the left-hand side is equipped with the action by left multiplication on $N_\l$. Hence, the quotient $\pi_\l{}^{-1}(\Orb_\l^-) /\!/ N_\l$ is isomorphic to the Slodowy slice as schemes:
    $$ S_f \cong \pi_\l{}^{-1}(\Orb_\l^-) /\!/ N_\l. $$
\end{theorem}

In the previous section, all the objects related to $\l$ were denoted with a ``$\l\,$'' subscript. For $\l = \{0\}$, we omit the subscript. In particular, $\n = \g_{\geqslant 1}$, $N = G_{\geqslant 1}$, there is an isomorphism $\g_1 \cong \Orb$ and~$\pi^{-1}(\Orb^-) = - \chi + (\g_{\geqslant 2})^\perp$.

\begin{corollary}[{\cite[Section 5.5]{gan2002quantization}}] \label{corollary:equivalent-poisson}
    For any isotropic subspace $\l$ in $\g_1$, there is an inclusion $\pi_\l{}^{-1}(\Orb_\l^-) \subseteq  \pi^{-1}(\Orb^-)$ which induces a Poisson isomorphism
    $$ \pi_\l{}^{-1}(\Orb_\l^-) /\!/ N_\l \cong  \pi^{-1}(\Orb^-) /\!/ N. $$
\end{corollary}

\subsection{Affine W-algebras related to a fixed isotropic subspace} \label{subsection:definition-w-algebras}

Consider the tensor product of vertex algebras $\V^k(\g) \otimes_\C \A(\l^{\perp, \omega} / \l)$. It is equipped with a vertex algebra map
$$ \Upsilon_\l : \V(\n_\l) \longrightarrow \V^k(\g) \otimes_\C \A(\l^{\perp, \omega} / \l) $$
defined for any element $x$ in $\n_\l$ by
$$ \Upsilon_\l(x) \defeq \begin{cases}
    &  x + \psi_{(x \bmod \l)}  \quad   \text{if} \quad  x \in \l^{\perp, \omega} \\
    &  x + \chi(x) \One  \quad   \text{if} \quad x \in \g_2 \\
    &  x  \quad   \text{otherwise}.
\end{cases} $$
The data of this chiral moment map allows to define the BRST cochain complex, denoted by $(\Cpx_{\l}, \mathcal{d}_\l)$, where
$$ \Cpx_{\l}^\bullet \defeq \V^k(\g) \otimes_\C \A(\l^{\perp, \omega} / \l) \otimes_\C \F^\bullet(\n_\l \oplus {\n_\l}^*). $$
Denote by $ \W^k(\g, f, H, \l) \defeq \Hgy^\bullet(\Cpx_{\l}, \mathcal{d}_\l) $ the vertex superalgebra constructed by this cohomology.

\begin{remark}
    This is the affine analogue of the BRST complex used in \cite{dandrea2006three} to construct finite W-algebra. All these constructions are well-known for finite W-algebras and provide isomorphic W-algebras \cite{gan2002quantization, brundan2007good,brundan2008highest-weight}. 
    
    The W-algebra~$\W^k(\g, f, H, \{0\})$ corresponding to the zero isotropic subspace is the one constructed by Kac, Roan, and Wakimoto in \cite{kac2003quantum}. In \cite{arakawa2015localization}, an $\hbar$-adic version of the affine W-algebra algebra, for a Lagrangian subspace $\l$, is constructed. To our knowledge, the W-algebra~$\W^k(\g, f, H, \l)$ for nonzero $\l$ appears only in \cite{arakawa2024arc}, for $\l$ Lagrangian. 
\end{remark}

Consider the Poisson variety $ \g^* \times (\l^{\perp, \omega} / \l)$. It has a natural diagonal action of~$N_\l$ given by the coadjoint action on $\g^*$ and the isomorphism of $\l^{\perp, \omega} / \l \cong \Orb_\l$. There is a moment map
$$ \mu_\l : \g^* \times (\l^{\perp, \omega} / \l) \longrightarrow \n_\l{}^*, \quad (\xi, v \bmod \l) \longmapsto \pi_{\l}(\xi) + \overline{\chi}_\l + \ad^*(v)\overline{\chi}_\l $$
It is clear that the natural projection $\g^* \times  (\l^{\perp, \omega} / \l) \twoheadrightarrow \g^*$ induces an $N_\l$-isomorphism
$$ \mu_\l^{-1}(0) \cong \pi_\l{}^{-1}(\Orb_\l^-). $$

If one takes the graded object corresponding to the Li filtration on the homomorphism $ \Upsilon_\l : \V(\n_\l) \rightarrow \V^k(\g) \otimes_\C \A(\l^{\perp, \omega} / \l)$ of vertex algebras, one gets the arcs of the moment map $\mu_\l : \g^* \times (\l^{\perp, \omega} / \l) \rightarrow {\n_\l}^*$. Denote by $(C_{\l, \infty}^\bullet, d_{\l, \infty})$ the associated Poisson vertex BRST complex, where:
$$ C_{\l, \infty}^\bullet \defeq \C[\Jet_\infty \g^*] \otimes_\C \C[\Jet_\infty (\l^{\perp, \omega} / \l)] \otimes_\C \Alt^\bullet_\infty(\n_\l \oplus \n_\l{}^*). $$

The following theorem generalises \cite[Theorem 4.1]{kac2004quantum} for $\l = \{0\}$ and \cite[Theorem 9.7]{arakawa2024arc} for $\l$ being Lagrangian.

\begin{theorem} \label{theorem:vanishing-w-algebra}
    The following cohomology vanishes outside degree $0$:
    $$ \Hgy^n(\Cpx_{\l}, \mathcal{d}_\l) \quad \text{for} \quad n \neq 0,  $$ 
    and there is a natural isomorphism 
    $$ \gr_{\Fil}\Hgy^0(\Cpx_{\l}^\bullet, \mathcal{d}_\l) \overset{\sim}{\longrightarrow} \Hgy^0(C_{\l, \infty}^\bullet, d_{\l, \infty}) $$
    where the filtration $\Fil$ on the cohomology $\Hgy^0(\Cpx_{\l}, \mathcal{d}_\l)$ is induced by the Li filtration on the complex $(\Cpx_{\l}, \mathcal{d}_\l)$. In particular, $ \W^k(\g, f, H, \l) = \Hgy^0(\Cpx_{\l}, \mathcal{d}_\l)$ is a purely even vertex algebra.
\end{theorem}

It is a straightforward generalisation~\cite[Theorem~9.7]{arakawa2024arc}, we only sketch the proof.

\begin{proof}[Sketch of proof]
    We first apply \Cref{theorem:quotient-quasi-isomorphic} to the BRST complex $(\Cpx_{\l}, \mathcal{d}_\l)$ where the vertex algebra $\V^k(\g) \otimes_\C \A(\l^{\perp, \omega} / \l)$ is graded by the Hamiltonian operator given by~$\Ham^\old \defeq \Ham^\g + L^\A_0$. So, we get an isomorphism $$\Hgy^n(\Cpx_{\l}^\bullet, \mathcal{d}_\l) \cong \Hgy^n(\Cpx_{\l, +}^\bullet, \mathcal{d}_{\l, +}),$$ where $(\Cpx_{\l, +}^\bullet, \mathcal{d}_{\l, +})$ is the quotient complex defined as in \Cref{subsection:quotient}. See \cite[Proposition 9.3]{arakawa2024arc} for details.

    Fix a basis $\{x_i\}_{i \in I(\n_\l)}$ of $\n_\l$ which is homogeneous for the good grading and denote by $\delta(x_i)$ the degree of $x_i$. Equip the complex $\Cpx_{\l}$ with the Hamiltonian operator
    $$ \Ham^\new \defeq \Ham^\g + \frac{1}{2}(\partial H)_{(1)} + L^\A_0 + L^\F(m_\bullet)_0, \quad \text{where} \quad m_\bullet = \frac{1}{2} \delta(x_\bullet). $$
    Thanks to the isomorphism $N_\l \times S_f \cong {\mu_\l}^{-1}(0)$ stated in \Cref{theorem:slodowy-isomorphism}, we can apply \Cref{theorem:vanishing-vertex} and conclude the proof.
\end{proof}

\subsection{Equivalence between the definitions} \label{subsection:equivalence-definitions}

For the isotropic subspace $\l = \{0\}$, the vertex algebra BRST complex is denoted by $(\Cpx^\bullet, \mathcal{d})$, with 
$$ \Cpx^\bullet = \V^k(\g) \otimes_\C \A(\g_1) \otimes_\C \F^\bullet(\n \oplus \n^*). $$
We want to compare its cohomology to the one of the BRST complex for a generic isotropic subspace $\l$ of $\g_1$, denoted by $(\Cpx_{\l}^\bullet, \mathcal{d}_{\l})$. 

Denote by $\A(\l^{\perp, \omega})$ the vertex subalgebra of the Weyl vertex algebra $\A(\g_1)$ which is strongly generated by the $\psi_{v}$ for $v$ in $\l^{\perp, \omega}$. The map $$ \Upsilon : \V(\g_{\geqslant 1}) \longrightarrow \V^k(\g) \otimes_\C \A(\g_1)$$ restricts to a vertex algebra homomorphism
$$ \Upsilon_\interm : \V(\n_\l) \longrightarrow \V^k(\g) \otimes_\C \A(\l^{\perp, \omega}). $$
As in \cite[Paragraph 3.2.5]{arakawa2015localization}, we can define an \emph{intermediary BRST complex} denoted by $(\Cpx_{\interm}^\bullet, \mathcal{d}_{\interm})$ where $\Cpx_{\interm}$ is the subalgebra of $\Cpx$ defined by
$$ \Cpx_{\interm}^\bullet \defeq \V^k(\g) \otimes_\C \A(\l^{\perp, \omega}) \otimes_\C \F^\bullet(\n_\l \oplus \n^*), $$
and the operator $\mathcal{d}_{\interm}$ is the restriction of the coboundary operator $\mathcal{d}$ to the subalgebra~$\Cpx_{\interm}$.

There is an inclusion homomorphism $\Pi_1 : (\Cpx_{\interm}^\bullet, \mathcal{d}_{\interm}) \hookrightarrow (\Cpx^\bullet, \mathcal{d})$ of cochain complexes and of vertex superalgebras. The projections $\l^{\perp, \omega} \twoheadrightarrow \l^{\perp, \omega} / \l$ and $\n^* \twoheadrightarrow {\n_\l}^*$ induce a surjective homomorphism $\Pi_2 : (\Cpx_{\interm}^\bullet, \mathcal{d}_{\interm}) \twoheadrightarrow (\Cpx_{\l}^\bullet, \mathcal{d}_{\l})$ of cochain complexes and vertex superalgebras.

Let us denote by $\Orb_{\interm} \defeq \Ad^*(N) \overline{\chi}_\l$ the $N$-orbit in $\n_\l{}^*$ of the restriction of $\chi$ to~$\n_{\l}$. The isomorphism $ \sigma :  \g_1 \overset{\sim}{\rightarrow} \Orb_{\max} $ induces an isomorphism 
$$ \sigma_{\interm} :  \g_1 / \l \longrightarrow \Orb_{\interm}, \quad (v \bmod \l) \longmapsto \overline{\chi}_\l + \ad^*(v)\overline{\chi}_\l. $$ 
On the coordinate ring side, we see that $ \C[\g_1 / \l] = \operatorname{Sym} \l^{\perp, \omega}$. The latter is a Poisson subalgebra of $\operatorname{Sym} \g_1$, whence $\g_1 / \l \cong \Orb_\interm$ is a Poisson variety.

\begin{remark}
    In \cite{arakawa2015localization}, they do not need an intermediary orbit because their choice of intermediate complex is slightly different from ours. We modify their construction with the goal of proving \Cref{main:new-brst-complex}~(\ref{theorem:new-construction-w-algebra}), for which such an intermediary orbit is necessary. 
\end{remark}

The intermediary moment map is defined as
$$ \mu_{\interm} : \g^* \times (\g_1 / \l) \longrightarrow \n_\l{}^*, \quad (\xi, v) \longmapsto \pi_{\l}(\xi) + \overline{\chi}_\l + \ad^*(v)\overline{\chi}_\l. $$
It is an $N$-equivariant homomorphism because $N_\l$ is a normal subgroup of~$N$.

\begin{lemma} \label{lemma:intermediary-isomorphisms}
    The projection map and the embedding
    $$  \g^* \times \g_1  \twoheadlongrightarrow \g^* \times (\g_1 / \l)  \quad \text{and} \quad  \g^* \times (\l^{\perp, \omega} / \l) \hooklongrightarrow \g^* \times (\g_1 / \l)  $$
    induce Poisson isomorphisms between the corresponding Hamiltonian reductions:
    $$ \mu_{\l}{}^{-1}(0) /\!/ N_{\l} \cong \mu_{\interm}{}^{-1}(0) /\!/ N \cong \mu^{-1}(0) /\!/ N. $$
\end{lemma}

\begin{proof}
    The projection map $\g^* \times \g_1  \twoheadrightarrow \g^* \times (\g_1 / \l)$ induces an $N$-equivariant isomorphism~$ \mu^{-1}(0) \cong \mu_{\interm}{}^{-1}(0)$. It is clear because of the parametrization $\sigma_{\interm}$ of the intermediary orbit $\Orb_{\interm}$. Then the lemma follows from \Cref{corollary:equivalent-poisson}.
\end{proof}
 
\begin{theorem} \label{theorem:equivalent-vertex}
    These inclusion and projection maps of cochain complexes
    $$ \Pi_1 : (\Cpx_{\interm}^\bullet, \mathcal{d}_{\interm}) \hooklongrightarrow (\Cpx^\bullet, \mathcal{d}) \quad \text{and} \quad \Pi_2 : (\Cpx_{\interm}^\bullet, \mathcal{d}_{\interm}) \twoheadlongrightarrow (\Cpx_{\l}^\bullet, \mathcal{d}_{\l}) $$
    induce vertex algebra isomorphisms between their cohomologies:
    $$ \Hgy^0(\Cpx_{\l}^\bullet, \mathcal{d}_{\l}) \cong  \Hgy^0(\Cpx_{\interm}^\bullet, \mathcal{d}_{\interm}) \cong \Hgy^0(\Cpx^\bullet, \mathcal{d}). $$
\end{theorem}

\begin{proof}
    Thanks to \Cref{lemma:intermediary-isomorphisms}, we can apply \Cref{theorem:quotient-quasi-isomorphic,theorem:vanishing-vertex} to the cochain complex $(\Cpx_{\interm}^\bullet, \mathcal{d}_{\interm})$ seen as a subcomplex of $(\Cpx^\bullet, \mathcal{d})$ and we get:
    \begin{align*}
        \Hgy^n(\Cpx_{\interm}^\bullet, \mathcal{d}_{\interm}) & = 0 \quad \text{for} \quad n \neq 0, \\
        \gr_{\Fil}\Hgy^0(\Cpx_{\interm}^\bullet, \mathcal{d}_{\interm}) & \cong \Hgy^0(\gr_{\Li} \Cpx_{\interm}^\bullet, \gr_{\Li} \mathcal{d}_{\interm}).
    \end{align*}

    The cochain maps $\Pi_1$ and $\Pi_2$ induce vertex algebra maps after taking the cohomology, denoted by
    $$ \Psi_1 : \Hgy^0(\Cpx_{\interm}^\bullet, \mathcal{d}_{\interm}) \longrightarrow \Hgy^0(\Cpx^\bullet, \mathcal{d}) \quad \text{and} \quad \Psi_2 : \Hgy^0(\Cpx_{\interm}^\bullet, \mathcal{d}_{\interm}) \longrightarrow \Hgy^0(\Cpx_{\l}^\bullet, \mathcal{d}_{\l}). $$
    Consider the Poisson vartiety isomorphisms of \Cref{lemma:intermediary-isomorphisms}. They induce isomorphisms of Poisson vertex algebras after passing to arc spaces and coordinate rings. These isomorphisms coincide with the associated graded maps
    $$ \gr_{\Fil} \Psi_1 : \Hgy^0(\gr_{\Li} \Cpx_{\interm}^\bullet, \gr_{\Li} \mathcal{d}_{\interm}) \longrightarrow \Hgy^0(\gr_{\Li} \Cpx^\bullet, \gr_{\Li} \mathcal{d})$$
    and
    $$\gr_{\Fil} \Psi_2 : \Hgy^0(\gr_{\Li} \Cpx_{\interm}^\bullet, \gr_{\Li} \mathcal{d}_{\interm}) \longrightarrow \Hgy^0(\gr_{\Li} \Cpx_{\l}^\bullet, \gr_{\Li} \mathcal{d}_{\l}), $$
    hence $\gr_{\Fil} \Psi_1$ and $\gr_{\Fil} \Psi_2$ are isomorphisms.

    The Hamiltonian operator $\Ham^\new$ introduced in the proof of Theorem \ref{theorem:vanishing-w-algebra} induces a nonnegative grading on the three vertex algebras 
    \begin{equation}
        \Hgy^0(\Cpx_{\l}^\bullet, \mathcal{d}_{\l}), \quad \Hgy^0(\Cpx_{\interm}^\bullet, \mathcal{d}_{\interm}) \quad \text{and} \quad \Hgy^0(\Cpx^\bullet, \mathcal{d}), \label{equation:three-cohomologies}
    \end{equation}
    and the maps $\Psi_1$ and $\Psi_2$ commute with the Hamiltonian operators. The filtrations~$\Fil$ induced on the cohomologies by the Li filtrations on the BRST complexes are finite on each homogeneous component of these Hamiltonian gradings. It follows from the fact that the cohomologies \eqref{equation:three-cohomologies} are isomorphic to 
    $$ \Hgy^0(\Cpx_{\l, +}^\bullet, \mathcal{d}_{\l, +}), \quad \Hgy^0(\Cpx_{\interm, +}^\bullet, \mathcal{d}_{\interm, +}) \quad \text{and} \quad \Hgy^0(\Cpx^\bullet_+, \mathcal{d}_+), $$
    and it follows from the proof of \Cref{lemma:li-quotient} that the filtrations $\Fil$ on the homogeneous components of these cohomologies are finite.
    
    Because the maps $\gr_{\Fil} \Psi_1$ and~$\gr_{\Fil} \Psi_2$ are isomorphisms and because the filtrations are finite, the maps $\Psi_1$ and $\Psi_2$ restrict to isomorphisms between the homogeneous components, so they are indeed isomorphisms.
\end{proof}

Denote by $\W^k(\g, f, H)$ the unique vertex algebra obtained by computing the cohomology $\Hgy^0(\Cpx_{\l}^\bullet, \mathcal{d}_{\l})$ for any isotropic subspace $\l$ of $\g_1$. It is called the \emph{affine W-algebra} associated with the pair $(f, H)$ at level $k$.

\begin{remark}
    The construction of Slodowy slices in \cite{gan2002quantization} suggests another BRST complex to define $\W^k(\g, f)$. Take $\m_{\l} \defeq \l \oplus \g_{\geqslant 2}$ and set $ \widetilde\Cpx_\l^\bullet \defeq \V^k(\g) \otimes_\C \F(\n_\l \oplus \m_\l{}^*)$. As above, one can show a vertex algebra isomorphism $\Hgy^0\big(\widetilde\Cpx_\l^\bullet\big) \cong \W^k(\g, f)$.
\end{remark}

\begin{theorem} \label{theorem:equivalent-constructions}
    For any nilpotent elements $f, f'$ in the orbit $\mathbf{O}$ and for any semisimple elements $H, H'$ defining good gradings for these nilpotent elements, there is a vertex algebra isomorphism between the associated W-algebras:
    $$ \W^k(\g, f, H) \cong \W^k(\g, f', H'). $$
\end{theorem}

\begin{proof}
    Using \Cref{theorem:equivalent-vertex}, the arguments used in \cite[Sections 4 and 5]{brundan2007good} to prove the analogue statement for finite W-algebras by the way of adjacent good gradings can be applied to affine W-algebras. 
\end{proof}

As a consequence, the affine W-algebras $\W^k(\g, f, H)$ constructed above only depend on the nilpotent orbit $\mathbf{O}$ which contains $f$. They are called \emph{the affine W-algebra associated with the orbit $\mathbf{O}$ at level $k$} and denoted by~$\W^k(\g, f)$.

\subsection{About the structure of the affine W-algebra} \label{subsection:structure-w-algebra}

Let us consider the construction of the W-algebra for the trivial isotropic subspace $\l = \{0\}$:
$$ \W^k(\g, f) = \Hgy^0(\Cpx^\bullet, \mathcal{d}), \quad \text{where} \quad  \Cpx^\bullet \defeq \V^k(\g) \otimes_\C \A(\g_1) \otimes_\C \F^\bullet(\n \oplus \n^*) $$
and $\n$ stands for the nilpotent Lie algebra $\g_{\geqslant 1}$. This is the construction originally done by Kac, Roan and Wakimoto in \cite{kac2003quantum, kac2004quantum}. 

Let us denote by $\{x_i\}_{i \in I(\g)}$ a basis of the Lie algebra $\g$ containing a basis~$\{x_i\}_{i \in I(\n)}$ of the subagebra $\n$ and assume that these bases are homogeneous the good grading of $\g$. Denote by $\delta(x_i)$ the degree of $x_i$. For any element $x$ in the Lie algebra $\g$, introduce the element
$$ J^x \defeq x + \sum_{j,k \in I(\n)} c_j^k(x) \NO{\phi_k \phi^*_j} \quad \text{in} \quad \Cpx_0, $$
where the $c_j^k(x)$ are the structure constants defined by the relations
\begin{equation}
   [x, x_j] = \sum_{k \in I(\g)} c_j^k(x) x_k. \label{equation:structure-coefficients}
\end{equation}

The BRST complex $ \Cpx^\bullet $ contains two interesting vertex subalgebras. The first one, denoted by $\Cpx_{+}$, is strongly generated by $J^x$ for $x$ in the negative parabolic subalgebra $\g_{\leqslant 0}$, the subalgebra $\A(\g_1)$ and the positively charged generators~$\phi^*_\xi$ for~$\xi$ in $\n^*$. The second one, denoted by $\Cpx_{-}$, is strongly generated by the negatively charged generators $\phi_x$ for $x$ in $\n$ and their images $\mathcal{d}(\phi_x)$ by the coboundary operator. Both are graded by the charge and stable by $\mathcal{d}$. There is the tensor product decomposition $\Cpx^\bullet = \Cpx_+^\bullet \otimes_\C \Cpx_-^\bullet$. Note that the subcomplex $\Cpx_+^\bullet$ is nonnegatively graded for the charge, so $$\Hgy^0(\Cpx_+^\bullet, \mathcal{d}) = \Ker(\mathcal{d}) \cap \Cpx_+^0.$$

\begin{proposition}[{\cite[(4.5) and Theorem 4.1]{kac2004quantum}}] \label{proposition:w-embeds-in-complex} 
    The vertex algebra inclusion~$\Cpx_+^\bullet \hookrightarrow \Cpx^\bullet $ induces an isomorphism between the cohomologies:
    $$ \Ker(\mathcal{d}) \cap \Cpx_+^0 \cong \W^k(\g, f). $$
\end{proposition}

We equip now the complex $\Cpx$ with the Hamiltonian operator
\begin{equation} \Ham \defeq \Ham^\g + \frac{1}{2}(\partial H)_{(1)} + L^\A_0 + L^\F(m_\bullet)_0, \quad \text{where} \quad m_\bullet = \frac{1}{2} \delta(x_\bullet). \label{equation:hamiltonian-w-algebra} \end{equation}
It commutes with coboundary operator $\mathcal{d}$ and then it induces a grading on the W-algebra $\W^k(\g, f)$. If the element $x$ in $\g$ is homogeneous of degree $\delta(x)$ for the good grading, then the conformal degree of $J^x$ is $ \Delta(J^x) = 1 - \frac{1}{2} \delta(x)$.

\begin{theorem}[{\cite[Theorem 4.1]{kac2004quantum}}] \label{theorem:basis-W-algebra}
    Let $\{x_i\}_{i=1}^{\mathcal{l}}$ be a basis of $\g^f$ which is homogeneous for the good grading. For each index $1 \leqslant i \leqslant \mathcal{l}$, one can construct an element $ J^{\{i\}}$ in $\Ker(\mathcal{d}) \cap \Cpx_+^0 \cong \W^k(\g, f) $ such that the following properties hold.
    \begin{enumerate}
        \item The family $\{J^{\{i\}}\}_{i=1}^{\mathcal{l}}$ freely generates the W-algebra.

        \item The Hamiltonian degree of each $J^{\{i\}}$ is $1-\frac{1}{2} \delta(x_i)$.

        \item The elements $J^{\{i\}}$ are of the form $J^{\{i\}} = J^{x_i} + T_i$, where $T_i$ is a linear combination of normally ordered products of elements of the form
        $$ \partial^n(J^{x_j}) \quad \text{where} \quad n \in \Z_{\geqslant 0} \quad \text{and} \quad \delta(x_i) < \delta(x_j) \leqslant 0, $$
        or the form
        $$ \partial^n(\psi_v) \quad \text{where} \quad n \in \Z_{\geqslant 0} \quad \text{and} \quad v \in \g_1, $$
        so that the Hamiltoniam degree of each product equals $1-\frac{1}{2} \delta(x_i)$.
    \end{enumerate}
\end{theorem}

\subsection{Geometric interpretation of the Kac--Roan--Wakimoto embedding} \label{subsection:geometric-krw}

Let $G^\natural$ be the connected subgroup of $G$ whose Lie algebra is 
$$ \g^\natural  \defeq \g^f \cap \g_0. $$ 
According to \cite[Lemma 2.4]{premet2007enveloping}, the group $G^\natural$ acts in a nontrivial way on the Slodowy slice $S_f \cong (-\chi + (\g_{\geqslant 2})^\perp) /\!/ N$. There exists an analogue of this action on the affine W-algebra $\W^k(\g, f)$. Let~$\{v_i\}_{i=1}^{2 s}$ be a basis of the symplectic space $\g_1$ and denote by $\{v^i\}_{i=1}^{2 s}$ the symplectic dual basis, these vectors correspond to fields~$\psi_i$ and $\psi^i$ in the Weyl vertex algebra~$\A(\g_1)$.

\begin{proposition}[{\cite[Theorem 2.4]{kac2003quantum}}] \label{proposition:krw-embedding}
    For any $x$ in $\g^\natural$, the element of $\Cpx^0_+$ defined by the formula
    $$J^{\{x\}} \defeq J^{x} + \frac{1}{2}\sum_{i=1}^{2 s} \NO{\psi^i \psi_{[u_i, x]}} $$
    belongs to the W-algebra $\W^k(\g, f)$.

    Moreover, there is an invariant symmetric bilinear form $\tau_k$ on the Lie algebra $\g^\natural$ such that the mapping $ x \mapsto J^{\{x\}} $ induces an injective vertex algebra embedding
    $$ \Theta^\natural : \V^{\tau_k}(\g^\natural) \hooklongrightarrow \W^k(\g, f). $$
\end{proposition}

In \cite{kac2003quantum, kac2004quantum}, the formula defining $\Theta^\natural$ comes from the explicit computation of the low degree generators in Theorem \ref{theorem:basis-W-algebra}. Premet does the same for finite W-algebras \cite{premet2007enveloping}. We found out that this formula has in fact a very natural geometric interpretation, that will be useful from the reduction by stages perspective. 

The group $G^\natural$ normalizes $N = G_{\geqslant 1}$. Define $ N^\natural \defeq N \rtimes G^\natural$ the semi-direct product of these groups. Denote by $\n^\natural$ the associated Lie algebra, by $\overline{\chi}^\natural$ the restriction of $\chi$ to $\n^\natural$ and by~$\Orb^\natural$ its coadjoint orbit in $(\n^\natural)^*$. Denote by $\Orb$ the coadjoint orbit of $\overline{\chi}$ in $\n^*$, where~$\overline{\chi}$ denotes the restriction of $\chi$ to $\n$.

\begin{lemma}
    There is a symplectic isomorphism given by the map
    $$ \sigma^\natural : \g_1 \longrightarrow \Orb^\natural, \quad v \longmapsto \overline{\chi}^\natural+ \ad^*(v) \overline{\chi}^\natural + \frac{1}{2} \ad^*(v)^2 \overline{\chi}^\natural. $$
    Moreover, the restriction map $(\n^\natural)^* \twoheadrightarrow \n^*$ induces a symplectic isomorphism between the coadjoint orbits~$\Orb^\natural$ and $\Orb$.
\end{lemma}

\begin{proof}
    The argument is similar to \Cref{lemma:symplectic-orbit}, see \Cref{proposition:orbit-description} for details.
\end{proof}

The group $N^\natural$ acts by the diagonal action on $\g^* \times \g_1  \cong \g^* \times \Orb^\natural$ and there is a moment map
$$ \mu^{\natural} : \g^* \times \g_1\longrightarrow (\n^\natural)^*, \quad (\xi, v) \longmapsto \xi|_{\n^\natural} + \overline{\chi}^\natural+ \ad^*(v) \overline{\chi}^\natural + \frac{1}{2} \ad^*(v)^2 \overline{\chi}^\natural. $$
The restriction of this action to the normal subgroup $N$ corresponds to the usual moment map
$$ \mu : \g^* \times \g_1 \longrightarrow \n^*, \quad (\xi, v) \longmapsto \pi(\xi) + \overline{\chi} + \ad^*(v)\overline{\chi}. $$

Because of the semi-direct product decomposition of $N^\natural$, there is an induced action of the quotient group $G^\natural \cong N^\natural / N$ on the Hamiltonian reduction $ \mu^{-1}(0) /\!/ N $ and the moment map $\mu^{\natural} : \g^* \times \g_1 \rightarrow (\n^\natural)^*$ induces a moment map
$$ \theta^{\natural} : \mu^{-1}(0) /\!/ N \longrightarrow (\g^\natural)^*. $$

\begin{proposition} \label{proposition:krw-embedding-geometry}
    The moment map $\theta^{\natural} : \mu^{-1}(0) /\!/ N \rightarrow (\g^\natural)^*$ is the geometric analogue of the vertex algebra embedding $\Theta^\natural : \V^{\tau_k}(\g^\natural) \hookrightarrow \W^k(\g, f)$ in the sense that the following square of Poisson vertex algebra homomorphism commutes:
    $$ \begin{tikzcd}
        \gr_\Li \V^{\tau_k}(\g^\natural) \arrow[r, "\sim"] \arrow[d, "\gr_\Li \Theta^\natural"'] & \C[\Jet_\infty (\g^\natural)^*] \arrow[d, "(\Jet_\infty \theta^{\natural})^*"] \\
        \gr_\Li \W^k(\g, f) \arrow[r, "\sim"'] & \C[\Jet_\infty (\mu^{-1}(0) /\!/ N )].
    \end{tikzcd} $$
    The horizontal isomorphisms are provided by \Cref{subsubsection:Kac-Moody} and \Cref{theorem:vanishing-w-algebra}.
\end{proposition}

\begin{proof}
    The comorphism $(\mu^{\natural})^*$ of the moment map $\mu^{\natural}$ is given for $x$ in $\g^\natural$ by:
    $$ (\mu^{\natural})^*(x) = x + \frac{1}{2}\sum_{i=1}^{2 s} \psi^i \psi_{[v_i, x]}, $$
    which lies in $\C[\g^*] \otimes_\C \C[\g_1]$. Because $\theta^{\natural}$ is induced by $\mu^{\natural}$, the proposition follows.
\end{proof}

\section{Reduction by stages for W-algebras}  \label{section:reduction-stages}

We state \Cref{main:reduction-by-stages-w-algebras}~(\ref{theorem:reduction-by-stages-w-algebras}), that is to say reduction by stages for affine W-algebras under the conditions \eqref{conditions} described in \Cref{subsection:sufficient-conditions}. For the proof, we need a new construction of the affine W-algebra $\W^k(\g, f_2)$, see \Cref{main:new-brst-complex}~(\ref{theorem:new-construction-w-algebra}). In \Cref{subsection:examples}, we provide some examples of reduction by stages.

\subsection{Sufficient condition on the good gradings} \label{subsection:sufficient-conditions} 

For $i=1,2$, let $f_i$ be a nilpotent element in $\g$ and $H_i$ be an element in~$\h$ such that the grading
$$ \g = \bigoplus_{\delta \in\Z}\g^{(i)}_\delta, \quad \g^{(i)}_\delta \defeq \g^{(H_i)}_\delta, $$
is a good grading for $f_i$. The homogeneous subspace $\g^{(i)}_1$ is equipped with the symplectic form $\omega_i(u,v)  \defeq (f_i | [u,v])$ for $u, v$ in $\g^{(i)}_1$. Introduce the Lie subalgebra~$\g^{\natural, 1} \defeq \g^{(1)}_0 \cap \g^{f_1}$.

Since $[H_1, H_2]=0$, one gets a bigrading on $\g$:
$$ \g = \bigoplus_{\delta_1, \delta_2 \in \Z} \g_{\delta_1, \delta_2}, \quad \text{where} \quad \g_{\delta_1, \delta_2} \defeq \g^{(1)}_{\delta_1} \cap \g^{(2)}_{\delta_2}. $$
Set $f_0 \defeq f_2 - f_1$. Consider the following conditions:
\begin{equation}\tag{$\bigstar$}\label{conditions}
  \begin{split}
        & \g^{(1)}_{\geqslant 2}  \subseteq \g^{(2)}_{\geqslant 1} \subseteq \g^{(1)}_{\geqslant 0}, \quad  \g^{(1)}_1 \subseteq  \bigoplus_{\delta=0}^2 \g_{1, \delta}, \quad \g^{(2)}_1 \subseteq \bigoplus_{\delta=0}^2 \g_{\delta, 1}, \\
        & f_0 \in \g_{0, -2}.
    \end{split}
\end{equation}

\begin{proposition} \label{proposition:construction-groups}
    Assume the conditions~\eqref{conditions}. Then there exist an isotropic subspace $\l_1$ of $\g^{(1)}_1$ and an isotropic subspace $\l_2$ of $\g^{(2)}_1$ such that both are $H_1$ and $H_2$-stable, and the nilpotent algebras
    $$ \n_1 \defeq {\l_1}^{\perp, \omega_1} \oplus \g^{(1)}_{\geqslant 2} \quad \text{and} \quad \n_2 \defeq {\l_2}^{\perp, \omega_2} \oplus \g^{(2)}_{\geqslant 2} $$
    satisfy the following properties:
    \begin{enumerate}
        \item \label{proposition:construction-groups,item:ideal} the algebra $\n_1$ is an ideal of $\n_2$,

        \item \label{proposition:construction-groups,item:m0} there exists a subalgebra of $\n_2$, denoted by $\n_0$, such that it is contained in the Lie subalgebra $\g^{\natural, 1}$ and there is a decomposition $\n_2 = \n_1 \oplus \n_0$. 
    \end{enumerate}
\end{proposition}

In the rest of this section, assume the conditions \eqref{conditions}.

\begin{lemma}
	The forms $\omega_1$ and $\omega_2$ coincide on $\g_{1,1} = \g^{(1)}_1 \cap \g^{(2)}_1$ and the subspace $\g_{1,1}$ is a symplectic subspace of $\g^{(1)}_1$ and of $\g^{(2)}_1$.
\end{lemma}

\begin{proof}
	Take $x, y \in\g_{1,1}$. Then, $[x,y] \in \g_{2,2}$ is orthogonal to $f_0 \in \g_{0,-2}$. Hence
	$$ \omega_1(x,y)= (f_1 | [x,y]) = (f_2 | [x,y]) = \omega_2(x,y). $$
    Let us prove that $\omega_1$ is nondegenerate on $\g_{1,1}$. Because the $H_1$-grading is good, the adjoint action of $f_1$, that belongs to $\g_{-2,-2}$, induces an isomorphism~$\g_{1,1} \cong \g_{-1,-1} $. Moreover, the symmetric invariant bilinear form $(\bullet | \bullet)$ is perfect pairing between the subspaces $\g_{1,1} $ and $ \g_{-1,-1} $, hence $\omega_1$ is nondegenerate.
\end{proof}

\begin{lemma} \label{lemma:construction-L1}
	The subspace  $\l_1 \defeq \g_{1,2}$ is isotropic in $\g^{(1)}_1$ and its orthogonal is given by ${\l_1}^{\perp, \omega_1} = \g_{1,1} \oplus \g_{1,2}$.
\end{lemma}

\begin{proof}
    Because  of the hypotheses of \Cref{proposition:construction-groups}, the vector space $\g^{(1)}_1$ decomposes as $\g^{(1)}_1 = \g_{1,0} \oplus \g_{1,1} \oplus \g_{1,2}$. The adjoint action of $f_1$ induces an isomorphism~$\g_{1,2} \cong \g_{-1,0}$ and the bilinear form~$(\bullet | \bullet)$ is a perfect pairing between the subspaces $\g_{-1,0}$ and $\g_{1,0}$, so $\g_{1,0}$ and~$\g_{1,2}$ are perfectly paired by $\omega_1$ and have the same dimension. Hence, $\g_{1,0}  \oplus \g_{1,2}$ is a symplectic subspace of $\g^{(1)}_1$ and $\g_{1,2}$ is a Lagrangian subspace of this symplectic subspace.
    
    By comparing the bidegrees, one can also see that $\omega_1(\g_{1,1} , \g_{1,0} \oplus \g_{1,2}) = 0 $. Hence, the subspace $\g_{1,2}$ is isotropic and its orthogonal is $\g_{1,1} \oplus \g_{1,2}$.
\end{proof}

\begin{lemma}\label{lemma:construction-L2}
	\begin{enumerate}
		\item \label{lemma:construction-L2,item:subspace} The subspace $\g_{0,1} \oplus \g_{2,1}$ is a symplectic subspace of $\g^{(2)}_1$. 
		
		\item \label{lemma:construction-L2,item:lagrangian} There exists a subspace $\mathfrak{a} $ in $ \g_{0,1} \cap \g^{f_1} $ such that $\mathfrak{a} \oplus \g_{2,1}$ is a Lagrangian subspace of $\g_{0,1} \oplus \g_{2,1}$. 
		
		\item \label{lemma:construction-L2,item:L2} The subspace $\l_2 \defeq \mathfrak{a} \oplus \g_{2,1}$ is isotropic in $\g^{(2)}_1$ and its orthogonal is given by~${\l_2}^{\perp, \omega_2} = \mathfrak{a} \oplus \g_{1,1} \oplus \g_{2,1}$.
	\end{enumerate}
\end{lemma}

\begin{proof}
	Because $f_0$ belongs to $\g_{0,-2}$, one has the inclusion $[f_2, \g_{1,1}] \subseteq  \g_{1,-1} \oplus \g_{-1,-1}$	and this space is orthogonal to $\g_{0,1} \oplus \g_{2,1}$ by the pairing $(\bullet | \bullet)$. Hence, 
	$$ \omega_2(\g_{0,1} \oplus \g_{2,1}, \g_{1,1}) = 0. $$
	The symplectic form $\omega_2$ has to be nondegenerate on $\g_{0,1} \oplus \g_{2,1}$, whence it is a symplectic subspace and \eqref{lemma:construction-L2,item:subspace} is proved.
	
	In the same way, one has $[f_2, \g_{2,1}] \subseteq \g_{0,-1} \oplus \g_{2,-1}$, so the subspace $\g_{2,1}$ is isotropic in $\g_{0,1} \oplus \g_{2,1}$.  It can be extended to a Lagrangian subspace of the form $\mathfrak{a} \oplus \g_{2,1}$, with $\mathfrak{a} \subseteq \g_{0,1}$. 
	
    We claim that $\mathfrak{a}$ is included in $\g_{0,1}$ belongs to $\g^{f_1}$. Any element in $\mathfrak{a}$ is of the form $x = u + v$ where $u \in \g^{f_1}$ and $v$ belongs to a complement of $\g^{f_1}$ in $\g$ which is~$H_1$ and $H_2$-stable. Then,
	$$ [f_2, x] = [f_0, u] + [f_1, v] + [f_0, v], $$
	where $[f_0, u] \in \g_{0,-1}$, $[f_1, v] \in \g_{-2,-1}$ and $[f_0, v] \in \g_{0,-1}$. Since $x $ belongs to $\mathfrak{a} \oplus \g_{2,1}$, which is isotropic, it implies that for all $y \in \g_{2,1}$,
	$$ \omega_2(x, y) = ([f_1, v]| y) = 0. $$
	Hence, $[f_1, v]$ is zero because the subspace $\g_{-2,-1}$ is paired with $\g_{2,1}$. Finally, we deduce that $v$ has to be zero because it lies in a complement of $\g^{f_1}$. This proves the existence of $\mathfrak{a}$ as in \eqref{lemma:construction-L2,item:lagrangian}.
	
	By construction, it is clear that $\mathfrak{a} \oplus \g_{2,1}$ is isotropic in $\g^{(2)}_1$ and its orthogonal is the subspace $\mathfrak{a} \oplus \g_{1,1} \oplus \g_{2,1}$, so \eqref{lemma:construction-L2,item:L2} follows.
\end{proof}

\begin{proof}[Proof of \Cref{proposition:construction-groups}]
	Consider the following subspaces of the Lie algebra $\g$:
    $$ \n_0 \defeq  \mathfrak{a} \oplus (\g^{(1)}_0 \cap \g^{(2)}_{\geqslant 2}), \quad \n_1 \defeq {\l_1}^{\perp, \omega_1} \oplus \g^{(1)}_{\geqslant 2}, \quad \n_2 \defeq  {\l_2}^{\perp, \omega_2} \oplus \g^{(2)}_{\geqslant 2}. $$
    We have the decomposition
    $$\g^{(2)}_{\geqslant 1} = (\g^{(2)}_1   \cap \g^{(1)}_0) \oplus (\g^{(2)}_1 \cap \g^{(1)}_1) \oplus (\g^{(2)}_1 \cap \g^{(1)}_2) \oplus (\g^{(2)}_{\geqslant 2} \cap \g^{(1)}_0) \oplus (\g^{(2)}_{2} \cap \g^{(1)}_1) \oplus (\g^{(2)}_{\geqslant 2} \cap \g^{(1)}_{\geqslant 2}),  $$
	whence we get the decomposition
    $$ \n_2 = \mathfrak{a} \oplus (\g^{(2)}_1 \cap \g^{(1)}_1) \oplus (\g^{(2)}_1 \cap \g^{(1)}_2) \oplus (\g^{(2)}_{\geqslant 2} \cap \g^{(1)}_0) \oplus (\g^{(2)}_{2} \cap \g^{(1)}_1) \oplus (\g^{(2)}_{\geqslant 2} \cap \g^{(1)}_{\geqslant 2}). $$
	We have the decomposition
    $$ \g^{(1)}_{\geqslant 1} =  (\g^{(1)}_1   \cap \g^{(2)}_0) \oplus (\g^{(1)}_1 \cap \g^{(2)}_1) \oplus (\g^{(1)}_1 \cap \g^{(2)}_2) \oplus (\g^{(1)}_{2} \cap \g^{(2)}_1) \oplus (\g^{(1)}_{\geqslant 2} \cap \g^{(2)}_{\geqslant 2}), $$
	whence we have the decomposition
    $$ \n_1 = (\g^{(1)}_1 \cap \g^{(2)}_1) \oplus (\g^{(1)}_1 \cap \g^{(2)}_2) \oplus (\g^{(1)}_{2} \cap \g^{(2)}_1) \oplus (\g^{(1)}_{\geqslant 2} \cap \g^{(2)}_{\geqslant 2}). $$

    To prove \eqref{proposition:construction-groups,item:ideal}, we just need to prove the inclusion $ [{\l_1}^{\perp, \omega_1}, \n_0] \subseteq {\l_1}^{\perp, \omega_1}$	because the $H_1$-grading is a Lie grading. But this inclusion is clear since
	$$ {\l_1}^{\perp, \omega_1} = \g^{(1)}_1 \cap \g^{(2)}_{\geqslant 1} \quad  \text{and} \quad \n_0 \subseteq \g^{(1)}_0 \cap \g^{(2)}_{\geqslant 1}.$$

	After comparing these decompositions, is is clear that $\n_2 = \n_1 \oplus \n_0$. To prove~\eqref{proposition:construction-groups,item:m0}, we just need to prove that $\n_0 \subseteq \g^{f_1}$. It is clear because of the inclusions
    $$ \mathfrak{a} \subseteq \g^{f_1} \quad \text{and} \quad  [f_1, \g^{(1)}_0 \cap \g^{(2)}_{\geqslant 2}] \subseteq \g^{(1)}_{-2} \cap \g^{(2)}_{\geqslant 0} = 0, $$
    where the last equality is a consequence of the inclusion $\g^{(1)}_{-2} \subseteq \g^{(2)}_{\leqslant -1}$.
\end{proof}

\subsection{Reduction by stages for Slodowy slices} \label{subsection:reduction-stages-slodowy}

If the conditions \eqref{conditions} hold, we get the nilpotent subalgebras $\n_2$, $\n_1$ and $\n_0$ such that the  corresponding unipotent subgroups of the reductive group $G$, denoted by $N_2$, $N_1$ and $N_0$, satisfy the semi-direct product decomposition $ N_2 = N_1 \rtimes N_0$. For $i=1,2$, the group $N_i$ acts on $\g^*$ by the restriction of the coadjoint action and this actions has a moment map given by the restriction map
$$ \pi_i : \g^* \longrightarrow {\n_i}^*, \quad \xi \longmapsto \xi|_{\n_i}. $$
Denote by $\overline{\chi_i}$ the restriction of the linear form $\chi_i = (f_i|\bullet)$ to the subalgebra $\n_i$ and denote by $\Orb_i \defeq \Ad^*(N_i)\overline{\chi_i}$ its coadjoint orbit.

The following theorem is a new formulation of \cite[Main Theorems 1 and~2]{genra2024reduction} under conditions \eqref{conditions}. Note that $\l_1$ and $\l_2$ may not be Lagrangian, in contrast to~\cite{genra2024reduction}.

\begin{theorem} \label{theorem:geometric-reduction-stages}
    Assume the conditions~\eqref{conditions}. Then Proposition~\ref{proposition:construction-groups} holds and we can use the objects introduced just before this theorem. We have the following consequences.
    \begin{enumerate}
        \item The coadjoint action of $N_0$ descends to the quotient~$\pi_1{}^{-1}(\Orb_1^-) /\!/ N_1$ with an induced moment map
        $$ \pi_0 : \pi_1{}^{-1}(\Orb_1^-) /\!/ N_1 \longrightarrow \n_0{}^*, \quad [\xi] \longmapsto  \xi|_{\n_0}. $$
        The linear form $\chi_2$ restricts to a character of $\n_0$, denoted by $\overline{\chi_0}$.

        \item there is an inclusion $\pi_2{}^{-1}(\Orb_2^-) \subseteq \pi_1{}^{-1}(\Orb_1^-) $ which induces a Poisson isomorphism:
        $$ \pi_2{}^{-1}(\Orb_2^-) /\!/ N_2 \cong {\pi_0}^{-1}(- \overline{\chi_0}) /\!/ N_0. $$  
    \end{enumerate}
\end{theorem}

The proof of Theorem \ref{theorem:geometric-reduction-stages} is almost identical to the construction done in \cite[Section 2]{genra2024reduction}. We only recall the important points. For $i=1,2$, embed the nilpotent element $f_i$ in an $\sl_2$-triple $(e_i, h_i, f_i)$ such that $e_i$ belongs to $\g^{(i)}_{2}$ and $h_i$ belongs to~$\g^{(i)}_0$, and let $ S_i \defeq - \chi_i + [\g, e_i]^\perp $ be the corresponding Slodowy slice. Moreover, accoding to \cite[Lemma 2.2.6]{genra2024reduction}, one can assume that $e_1$ belongs to $\g_{2,2}$ and $h_1$ belongs to $\g_{0,0}$.

According to Theoren \ref{theorem:slodowy-isomorphism}, the map
\begin{equation} \alpha_i : N_i \times S_i \longrightarrow \pi_i{}^{-1}(\Orb_i^-), \quad (g, \xi) \longmapsto \Ad^*(g) \xi \label{equation:gan-ginzburg-iso} \end{equation}
is an algebraic isomorphism. By the proof of Proposition \ref{proposition:construction-groups}, we have a symplectic isomorphism $ \g_{1,1} \cong \l_i / {\l_i}^{\perp, \omega_i} $ induced by the inclusion $\g_{1,1} \subseteq \l_i$. The following lemma follows from Lemma \ref{lemma:symplectic-orbit}.

\begin{lemma} \label{lemma:orbits-stages}
    Let $i$ be $1$ or $2$. We have the symplectic isomorphism
    $$ \sigma_i : \g_{1,1} \longrightarrow \Orb_i, \quad v \longmapsto \overline{\chi_i} + \ad^*(v) \overline{\chi_i}. $$
    Hence, the fiber of the orbit by the moment map is
    $$ \pi_i{}^{-1}(\Orb_i^-) = - \chi_i + \ad^*(\g_{1,1}) \chi_i \oplus \n_i{}^\perp.  $$
\end{lemma}

\begin{proof}[Sketch of proof of Theorem \ref{theorem:geometric-reduction-stages}]
    The natural projection $\n_2{}^* \twoheadrightarrow \n_1{}^*$ induces a symplectic isomorphism $\Orb_2 \cong \Orb_1$, and there is an inclusion of $ \pi_2{}^{-1}(\Orb_2^-)$ in $\pi_1{}^{-1}(\Orb_1^-)$, as direct consequences of \Cref{lemma:orbits-stages}.

    Because of the semidirect product decomposition, there is an induced action of~$N_0$ on the quotient 
    $\pi_1{}^{-1}(\Orb_1^-) /\!/ N_1$ and a moment map
    $$ \pi_0 : \pi_1{}^{-1}(\Orb_1^-) /\!/ N_1 \longrightarrow {\n_0}^*, \quad [\xi] \longmapsto \xi|_{\n_0}. $$

    By the same argument as \cite[Claim 2.4.3]{genra2024reduction}, the Gan-Ginzburg isomorphism
    $$ \alpha_1 : N_1 \times S_1 \longrightarrow \pi_1{}^{-1}(\Orb_1^-) \quad (g, \xi) \longmapsto \Ad^*(g) \xi,  $$
    restricts to an isomorphism
    \begin{equation}
        N_1 \times (- \chi_2 + [\g, e_1]^\perp \cap {\n_0}^\perp) \cong \pi_2{}^{-1}(\Orb_2^-). \label{equation:intermediate-stage}
    \end{equation}
    Using \eqref{equation:intermediate-stage}, we can conclude the proof as in \cite[Section 2.4]{genra2024reduction}.    
\end{proof}

\begin{corollary}[{\cite[Main Theorems 1 and 2]{genra2024reduction}}] \label{corollary:isomorphism-second-stage}
    There is an $N_0$-equivariant isomorphisms
    $$ N_0 \times S_2 \cong \pi_2{}^{-1}(\Orb_2^-) /\!/ N_1 \cong \pi_0{}^{-1}(-\overline{\chi}_0), $$
    where the left-hand side is equipped with the left mutliplication action on $N_0$.
\end{corollary}

\begin{remark}[{\cite[Proposition 2.3.1]{genra2024reduction}}]
    Let us denote by $\mathbf{O}_i$ the adjoint orbit of~$f_i$ and by $\overline{\mathbf{O}_i}$ its Zariski closure in $\g$. Under the hypotheses of Theorem \ref{theorem:geometric-reduction-stages}, then the inclusion $\overline{\mathbf{O}_1} \subseteq \overline{\mathbf{O}_2}$ holds.
\end{remark}

\subsection{Reduction by stages for affine W-algebras} \label{subsection:reduction-stages-w-algebra}

Introduce the nilpotent  Lie algebras
$$ \widetilde{\n}_1 \defeq \g^{(1)}_{\geqslant 1} \quad \text{and} \quad \widetilde{\n}_2 \defeq \g^{(1)}_{\geqslant 1} \oplus \n_0, $$
which satisfy the semi-direct product decomposition:
$$ \widetilde{\n}_2 = \widetilde{\n}_1 \oplus \n_0 \quad \text{and} \quad [\widetilde{\n}_1, \n_0] \subseteq \widetilde{\n}_1. $$
Pick a basis $\{x_i\}_{i \in I(\widetilde\n_2)}$ of $\widetilde{\n}_2$ which is the union of a basis $\{x_i\}_{i \in I(\widetilde\n_1)}$ of $\widetilde\n_1$ and a basis $\{x_i\}_{i \in I(\n_0)}$ of $\n_0$. 

\begin{example}\label{example:type-A3}
    Consider $\g \defeq \sl_4$ and denote by $E_{i,j}$ the elementary matrices. Let~$f_1 \defeq E_{4,1}$ be a minimal nilpotent element and $f_2 \defeq E_{4,1} + E_{3,2}$ be a rectangular nilpotent element. Take $H_1 \defeq E_{1,1} - E_{4,4}$ and $H_2 \defeq E_{1,1} + E_{2,2} - E_{3,3} - E_{4,4}$. Then, these data satisfy the conditions \eqref{conditions} and
    \begin{align*}
        \n_{1} = \left\{\left(\begin{smallmatrix} 1 & 0 & *& * \\ & 1 & 0 & * \\ & & 1 & 0 \\ & & & 1 \end{smallmatrix}\right)\right\}, & \quad \widetilde{\n}_1 = \g^{(1)}_{\geqslant 1} = \left\{\left(\begin{smallmatrix} 1 & * & *& * \\ & 1 & 0 & * \\ & & 1 & * \\ & & & 1 \end{smallmatrix}\right)\right\}, \\
        \g^{(2)}_{\geqslant 1} = \n_{2} = \left\{\left(\begin{smallmatrix} 1 & 0 & *& * \\ & 1 & * & * \\ & & 1 & 0 \\ & & & 1 \end{smallmatrix}\right)\right\} , & \quad \widetilde{\n}_2 = \left\{\left(\begin{smallmatrix} 1 & * & *& * \\ & 1 & * & * \\ & & 1 & * \\ & & & 1 \end{smallmatrix}\right)\right\}.
    \end{align*}
\end{example}

The Lie algebra $\widetilde{\n}_1$ corresponds to the choice of the zero isotropic subspace of $\g^{(1)}_1$ in the construction of the affine W-algebra $\W^k(\g, f_1)$ described in Subsection~\ref{subsection:definition-w-algebras}:
$$ \W^k(g, f_1) = \Hgy^0\big(\widetilde\Cpx^\bullet_1, \widetilde{\mathcal{d}}_1\big), \quad \text{where} \quad  \widetilde\Cpx^\bullet_1 \defeq \V^k(\g) \otimes_\C \A(\g^{(1)}_1) \otimes_\C \F^\bullet(\widetilde{\n}_1 \oplus \widetilde{\n}_1{}^*). $$
Moreover, as recalled in Proposition \ref{proposition:w-embeds-in-complex}, the affine W-algebra can be constructed as a vertex subalgebra of the BRST complex:
\begin{equation}
    \W^k(\g, f_1) \subseteq \Ker(\mathcal{d}_1) \cap \widetilde\Cpx^0_1. \label{equation:w-embeds-in-complex}
\end{equation}

Let~$\{v_i\}_{i=1}^{2 s}$ be a basis of the symplectic space $\g^{(1)}_1$ and denote by $\{v^i\}_{i=1}^{2 s}$ the symplectic dual basis. These vectors correspond to fields $\psi_i$ and $\psi^i$ in the Weyl vertex algebra~$\A(\g^{(1)}_1)$. Denote by $\{x_i\}_{i \in I(\widetilde\n_1)}$ a basis of the Lie algebra $\widetilde{\n}_1$ and by~$\phi_i, \phi_i^*$ the corresponding strong generators in the Clifford vertex algebra denoted by $\F^\bullet(\widetilde{\n}_1 \oplus \widetilde{\n}_1{}^*)$. Denote by $c_j^k(x)$ the structure coefficient as defined in the formula~\eqref{equation:structure-coefficients}.

\begin{lemma}
    There is an embedding of vertex algebras $ \V(\n_0) \hookrightarrow \W^k(\g, f_1) $
    given for any $x$ in $\n_0$ by the formula 
    $$  x \longmapsto  J^{\{x\}} = x + \frac{1}{2}\sum_{i=1}^{2 s} \NO{\psi^i \psi_{[v_i, x]}} + \sum_{j,k \in I(\widetilde\n_1)} c_j^k(x) \NO{\phi_k \phi^*_j}, $$
    where the right-hand side element is \emph{a priori} defined in $\widetilde\Cpx^0_1$.
\end{lemma}

\begin{proof}
    According to \Cref{proposition:krw-embedding}, there exists an invariant symmetric bilinear form~$\tau^{(1)}_k$ on the Lie algebra $\g^{\natural, 1} =  \g^{(1)}_0 \cap \g^{f_1}$ such that there is an explicit vertex algebra embedding $\V^{\tau^{(1)}_k}(\g^{\natural, 1}) \hookrightarrow \W^k(\g, f_1)$. The bilinear form~$\tau^{(1)}_k$ is identically zero on the nilpotent subalgebra $\n_0$ of $\g^{\natural, 1}$, the lemma follows.
\end{proof}

Because the restriction $\overline{\chi}_0$ of the linear form $\chi_2$ to the Lie subalgebra $\n_2$ is a character, we get a vertex algebra map
$$ \Upsilon_0 : \V(\n_0) \longrightarrow \W^k(\g, f_1)  $$
defined for a free generator $x$ in $\n_0$ by
$$ \Upsilon_0(x) \defeq \begin{cases}
     &   J^{\{x\}} + \chi_2(x) \One  \quad   \text{if} \quad x \in \g_{0,2} \\
     &   J^{\{x\}}  \quad   \text{otherwise}.
\end{cases} $$
We state \Cref{main:reduction-by-stages-w-algebras}.

\begin{theorem} \label{theorem:reduction-by-stages-w-algebras}
    Assume the conditions~\eqref{conditions}. Following \Cref{subsection:brst-vertex}, to the chiral comoment map $\Upsilon_0 : \V(\n_0) \rightarrow \W^k(\g, f_1)$ corresponds a BRST cochain complex~$(\Cpx_0^\bullet, \mathcal{d}_0)$, where $\Cpx_0 \defeq \W^k(\g, f_1) \otimes_\C \F^\bullet(\n_0 \oplus \n_0{}^*)$.    

    Then the cohomology of this complex is concentrated in degree $0$ and isomorphic to the affine W-algebra $\W^k(\g, f_2)$:
    $$ \Hgy^\bullet(\Cpx_0^\bullet, \mathcal{d}_0) \cong \delta_{\bullet = 0} \, \W^k(\g, f_2). $$
\end{theorem}

To prove Theorem \ref{theorem:reduction-by-stages-w-algebras}, we need to find a natural map
$$ \Hgy^0(\Cpx_0^\bullet, \mathcal{d}_0) \longrightarrow \W^k(\g, f_2). $$
To do so, we need to introduce a new construction of the affine-W-algebra $\W^k(\g, f_2)$ which relies on the nilpotent algebra $ \widetilde{\n}_2 = \g^{(1)}_{\geqslant 1} \oplus \n_0$,
which is natural because of our choice of the nilpotent Lie algebra $\widetilde{\n}_1 = \g^{(1)}_{\geqslant 1}$ to construct the affine W-algebra~$\W^k(\g, f_1)$. 

\begin{remark}
    In general, the Lie algebra $\widetilde{\n}_2$ is not equal to $\g^{(2)}_{\geqslant 1}$. In fact, it does not come from a good grading for $f_2$, see \Cref{example:type-A3}.
\end{remark}

Introduce the chiral comoment map
$$ \widetilde{\Upsilon}_2 : \V(\widetilde{\n}_2) \longrightarrow \V^k(\g) \otimes_\C \A(\g^{(1)}_1) $$
defined on strong generators $x \in \widetilde{\n}_2$ by:
$$\widetilde{\Upsilon}_2(x) \defeq \begin{cases}
     &  x + \frac{1}{2}\sum_{i=1}^{2 s} \NO{\psi^i \psi_{[v_i, x]}} + \chi_2(x)  \quad   \text{if} \quad x \in \n_0 \\
     & x + \psi_x  \quad   \text{if} \quad  x \in \g^{(1)}_1 \\
     & x + \chi_2(x) \One  \quad   \text{if} \quad x \in \g^{(1)}_2 \\
     & x  \quad   \text{otherwise}.
\end{cases} $$
Denote by $\big(\widetilde{\Cpx}^\bullet_2, \widetilde{\mathcal{d}}_2\big)$ the associated BRST complex, where
$$ \widetilde\Cpx^\bullet_2 \defeq \V^k(\g) \otimes_\C \A(\g^{(1)}_1) \otimes_\C \F^\bullet(\widetilde{\n}_2 \oplus \widetilde{\n}_2{}^*). $$
Recall that the differential $\widetilde{\mathcal{d}}_2$ is the $0$-th mode of the element
$$ \widetilde{Q}_2 \defeq \sum_{i \in I(\widetilde{\n}_2)} \NO{\widetilde{\Upsilon}_2(x_i) \phi_i^*} - \frac{1}{2} \sum_{i,j \in I(\widetilde{\n}_2)} \NO{\phi_{[x_i, x_j]} \phi_i^* \phi_j^*}. $$

As a charge-graded vector space, the cochain complex has the tensor product decomposition
$$ \widetilde{\Cpx}^\bullet_2 = \widetilde{\Cpx}^\bullet_1 \otimes_\C \F^\bullet(\n_0 \oplus \n_0{}^*).  $$
Because of the embedding \eqref{equation:w-embeds-in-complex}, one has a vertex superalgebra embedding 
$$ \Cpx^\bullet_0 \subseteq \widetilde{\Cpx}^\bullet_2. $$
Recall that the differential $\mathcal{d}_0$ on the cochain complex $\Cpx^\bullet_0$ is given by the 0-th mode of the element
$$ Q_0 \defeq \sum_{i \in I(\n_0)} \NO{\Upsilon_0(x_i) \phi_i^*} - \frac{1}{2} \sum_{i,j \in I(\n_0)} \NO{\phi_{[x_i, x_j]} \phi_i^* \phi_j^*}. $$

\begin{lemma}
    \begin{enumerate}
        \item The decomposition $\widetilde{Q}_2 = \widetilde{Q}_1 + Q_0$ holds in $\widetilde{\Cpx}^\bullet_2$. Hence, the decomposition $\widetilde{\mathcal{d}}_2 = \widetilde{\mathcal{d}}_1 + \mathcal{d}_0$ holds as operators on $\widetilde{\Cpx}^\bullet_2$.

        \item The $\lambda$-superbracket $\big[\widetilde{Q}_1{}_{\lambda} Q_0\big]$ is zero. Hence, the superbracket $[\widetilde{\mathcal{d}}_1, \mathcal{d}_0]$ of operators is zero.
    \end{enumerate}
\end{lemma}

\begin{proof}
    By definition, $ \NO{\phi_{[x_i, x_j]} \phi_i^*\phi_j^*} = \NO{\phi_{[x_i, x_j]} (\NO{\phi_i^*\phi_j^*})}$.
    Because their $\lambda$-bracket is zero, one can permute the last two terms:
     $$ \NO{\phi_{[x_i, x_j]} \phi_i^*\phi_j^*} = -\NO{\phi_{[x_i, x_j]} \phi_j^*\phi_i^*}. $$
    Then it is clear that 
    \begin{align*}
        \widetilde{Q}_2 = \widetilde{Q}_1 + \sum_{i \in I(\n_0)} \NO{\Big(x + \frac{1}{2}\sum_{i=1}^{2 s} \NO{\psi^i \psi_{[v_i, x]}} + \chi_2(x)\Big) \phi_i} + \sum_{\substack{i \in I(\n_0) \\ j \in I(\widetilde\n_1)}} \NO{\phi_{[x_i, x_j]} \phi_j^*\phi_i^*} \\
        - \frac{1}{2} \sum_{i,j \in I(\n_0)} \NO{\phi_{[x_i, x_j]} \phi_i^* \phi_j^*}.
    \end{align*}

    We recall the associator formula \cite[(1.40)]{desole2006finite}:
    \begin{multline*}
        \NO{(\NO{\phi_{[x_i, x_j]} \phi_j^*})\phi_i^*} - \NO{\phi_{[x_i, x_j]} (\NO{\phi_j^*\phi_i^*})} \\ = \sum_{n \geqslant 0}  \frac{1}{(n+1)!} \big( \NO{\partial^{n+1} \phi_{[x_i, x_j]} (\phi_j^*{}_{(n)} \phi_i^*)} - \NO{\partial^{n+1}\phi_j^*(\phi_{[x_i, x_j]} {}_{(n)} \phi_i^*)}\big).
    \end{multline*}
    If $i \in I(\n_0)$ and $j \in I(\widetilde{\n}_1)$, then one has for all nonnegative integer $n$:
    $$ \phi_j^*{}_{(n)} \phi_i^* = 0 \quad \text{and} \quad \phi_{[x_i, x_j]} {}_{(n)} \phi_i^* = 0, $$
    the second equality is due to the fact that $[x_i, x_j]$ belongs to $\widetilde\n_1$. Hence if~$i \in I(\n_0)$ and $j \in I(\widetilde{\n}_1)$, then
    $$ \NO{(\NO{\phi_{[x_i, x_j]} \phi_j^*})\phi_i^*} - \NO{\phi_{[x_i, x_j]} (\NO{\phi_j^*\phi_i^*})} = 0. $$
    We deduce that
    $$ \widetilde{Q}_2 = \widetilde{Q}_1 + \sum_{i \in I(\n_0)} \NO{(J^{\{x_i\}} + \chi_2(x)) \phi_i} - \frac{1}{2} \sum_{i,j \in I(\n_0)} \NO{\phi_{[x_i, x_j]} \phi_i^* \phi_j^*} = \widetilde{Q}_1 + Q_0. $$

    According to \cite[Theorem 2.4(a)]{kac2003quantum}, one has $\big[\widetilde{Q}_1{}_\lambda J^{\{x\}}\big] = 0$ for all $x$ in~$\n_0$. It follows that $\big[\widetilde{Q}_1{}_{\lambda} Q_0\big] = 0$.
\end{proof}

As a consequence of this lemma, the emedding \eqref{equation:w-embeds-in-complex} induces an embedding 
$$ \Ker(\mathcal{d}_0) \cap \Cpx^n_0 \subseteq \Ker(\widetilde{\mathcal{d}}_2) \cap \widetilde{\Cpx}^n_2 \quad \text{for} \quad n \in \Z. $$
Taking the cohomology, we get a natural map
\begin{equation}
    \Theta : \Hgy^\bullet(\Cpx_0^\bullet, \mathcal{d}_0) \longrightarrow \Hgy^\bullet\big(\widetilde\Cpx^\bullet_2, \widetilde{\mathcal{d}}_2\big). \label{equation:natural-map}
\end{equation}

To prove Theorem \ref{theorem:reduction-by-stages-w-algebras}, it is enough to show that this map is an isomorphism and to prove that the codomain is isomorphic to $\W^k(\g, f_2)$.

\begin{remark}
    If $\g^{(1)}_1 = \{0\}$, then $\n_1 = \widetilde{\n}_1$ and $\n_2 = \widetilde{\n}_2$. The construction is simplified because \Cref{main:new-brst-complex} is obvious in this case: one can use the ideas of~\cite{madsen1997secondary} and only needs homological algebra tools. To generalise the construction to the case $\g^{(1)}_1 \neq \{0\}$, we need to use some geometry.
\end{remark}

\subsection{New construction of the second Slodowy slice} \label{subsection:new-definition-slodowy}

\begin{lemma}
    The inclusion $[\widetilde{\n}_2, \widetilde{\n}_2] \subseteq \n_2$ holds.
\end{lemma}

\begin{proof}
    We use the decomposition $\widetilde{\n}_2 = \widetilde{\n}_1 \oplus \n_0$. We have the inclusions
    $$ [\widetilde{\n}_1, \widetilde{\n}_1]  \subseteq \g^{(1)}_{\geqslant 2} \subseteq \n_2 \quad \text{and} \quad [\n_0, \n_0] \subseteq \n_0 \subseteq \n_2. $$
    It remains to show the inclusion $[\widetilde{\n}_1, \n_0] \subseteq \n_2$.
    
    Recall the decompositions $ \widetilde{\n}_1 = \g_{1,0} \oplus (\g^{(1)}_{\geqslant 1} \cap \g^{(2)}_{\geqslant 1})$ and $\n_0 = \mathfrak{a} \oplus (\g^{(1)}_0 \cap \g^{(2)}_{\geqslant 2})$. The inclusion $[\widetilde{\n}_1, \n_0] \subseteq \n_2$ follows from the following inclusions:
    \begin{align*}
        [\widetilde{\n}_1, \g^{(1)}_0 \cap \g^{(2)}_{\geqslant 2}] & \subseteq [\g^{(2)}_{\geqslant 0}, \g^{(2)}_{\geqslant 2}] \subseteq [\g_{1,0}, \g_{0,1}] \subseteq \n_2, \\
        [\g^{(1)}_{\geqslant 1} \cap \g^{(2)}_{\geqslant 1}, \mathfrak{a}] & \subseteq [\g^{(1)}_{\geqslant 1} \cap \g^{(2)}_{\geqslant 1}, \g_{0,1}] \subseteq \g^{(2)}_{\geqslant 2} \subseteq \n_2,\\
        [\g_{1,0}, \mathfrak{a}] & \subseteq [\g_{1,0}, \g_{0,1}] \subseteq \g_{1,1} \subseteq \n_2.
    \end{align*}
\end{proof}

Take $i=1,2$. The nilpotent Lie algebra $\widetilde{\n}_i$ is the Lie algebra of a unipotent subgroup $\widetilde{N}_i$ of~$G$, which acts on $\g^*$ by the coadjoint action, and this action is Hamiltonian with the moment map given by the restriction map
$$ \widetilde{\pi}_i : \g^* \longrightarrow \widetilde{\n}_i{}^*, \quad \xi \longmapsto \xi|_{\widetilde{\n}_i}. $$
Denote by $\widetilde{\Orb}_i \defeq \Ad^*(\widetilde{N}_i) \widetilde{\chi}_i $ the coadjoint orbit of the linear form $\widetilde{\chi}_i \defeq \chi_i|_{\widetilde{\n}_i}$.

By the Gan--Ginzburg construction applied to $\widetilde{\n}_1 = \g^{(1)}_{\geqslant 1}$ (Subsection \ref{subsection:slodowy-poisson-structure}), there is a symplectic isomorphism
$$ \widetilde{\sigma}_1 : \g^{(1)}_1 \longrightarrow  \widetilde{\Orb}_1, \quad v \longmapsto \widetilde{\chi}_1 + \ad^*(v) \widetilde{\chi}_1. $$
Hence, the fiber of the orbit by the moment map is
$$ \widetilde{\pi}_1{}^{-1}(\widetilde{\Orb}_1^-) = \chi_1  + \ad^*(\g^{(1)}_1) \chi_1 \oplus (\g^{(1)}_{\geqslant 1})^\perp.  $$

There is a natural map given by the inclusion $\pi_1{}^{-1}(\Orb_1^-) \subseteq \widetilde{\pi}_1{}^{-1}\big(\widetilde{\Orb}_1^-\big)$. This inclusion is $N_0$-equivariant, each subspace is indeed $N_0$-stable because of the semidirect product decompositions. By Corollary \ref{corollary:equivalent-poisson}, this inclusion induces an isomorphism
$$ \pi_1{}^{-1}(\Orb_1^-) /\!/N_1 \cong \widetilde{\pi}_1{}^{-1}\big(\widetilde{\Orb}_1^-\big) /\!/ \widetilde{N}_1. $$
Moreover, the $N_0$-actions descend to both quotients and the isomorphism is $N_0$-equivariant.

Thanks to \Cref{example:type-A3}, we see that the Lie algebra $\widetilde{\n}_2$ is not covered by the Gan--Ginzburg construction in general. But analogous properties hold.

\begin{proposition} \label{proposition:orbit-description} 
    There is a symplectic isomorphism given by the map
    $$ \widetilde{\sigma}_2 : \g^{(1)}_1 \longrightarrow \widetilde{\Orb}_2, \quad v \longmapsto \widetilde{\chi}_2 + \ad^*(v) \widetilde{\chi}_2 + \frac{1}{2} \ad^*(v)^2 \widetilde{\chi}_2. $$
    Moreover, the restriction map $\widetilde{\n}_2{}^* \twoheadrightarrow \widetilde{\n}_1{}^*$ induces a symplectic isomorphism between the coadjoint orbits $\widetilde{\Orb}_2$ and $\widetilde{\Orb}_1$. Whence, the following diagram commutes:
    $$ \begin{tikzcd}[row sep=tiny]
        & \widetilde{\Orb}_2 \arrow[dd, "\sim"'{sloped}] \arrow[r, hook]& \widetilde{\n}_2{}^* \arrow[dd, two heads]  \\
        \g^{(1)}_1 \arrow[ru, "\widetilde{\sigma}_2", "\sim"'{sloped}] \arrow[rd, "\widetilde{\sigma}_1"', "\sim"{sloped}] & & \\
        &  \widetilde{\Orb}_1 \arrow[r, hook] & \widetilde{\n}_1{}^*.
    \end{tikzcd} $$
\end{proposition}

\begin{proof}
    Notice the equality $ \Ad^*(N_0) \widetilde{\chi}_2 = \{\widetilde{\chi}_2\}$ which follows from the decomposition $f_2 = f_1 + f_0$, the facts that $f_1$ is a fixed point for the action of $N_0$ and that the restriction of $\widetilde{\chi}_2$ to $\n_0$ is the character $\overline{\chi}_0$. Because there is the semidirect product decomposition $\widetilde{N}_2 = \widetilde{N}_1 \rtimes N_0$, one has
    $$ \widetilde{\Orb}_2 = \Ad^*(\widetilde{N}_1) \Ad^*(N_0) \widetilde{\chi}_2 = \Ad^*(\widetilde{N}_1)\widetilde{\chi}_2. $$
    
    The action of the group $\widetilde{N}_1$ is the exponential of the action of its Lie algebra $\widetilde{\n}_1$:
    $$ \Ad^*(\widetilde{N}_1)\widetilde{\chi}_2 = \Big\{\sum_{k \geqslant 0} \frac{1}{k!} \ad^*(v)^k \widetilde{\chi}_2 ~ | ~ v \in \widetilde{\n}_1 \Big\}. $$

    Take $k \geqslant 3$ and $v$ in $\widetilde{\n}_1$. Since $f_2$ is in $\g^{(1)}_{\geqslant -2}$, $\ad(v)^k(f_2)$ belongs to $\g^{(1)}_{\geqslant k - 2} \subseteq \g^{(1)}_{\geqslant 1}$. The inclusion $\g^{(1)}_{\geqslant 1} \subseteq \widetilde{\n}_2{}^\perp$ implies the equality $(\ad(v)^k(f_2) | \widetilde{\n}_2) = 0$. Hence,
    $$ \ad^*(v)^k \widetilde{\chi}_2 = 0. $$ 
    
    Let us check that if $v$ belongs to $\g^{(1)}_{\geqslant 2}$ and $k \geqslant 1$, then $\ad^*(v)^k \widetilde{\chi}_2 = 0$. The element $\ad(v)^k f_0$ belongs to $\g^{(1)}_{\geqslant 2 k} \subseteq \g^{(1)}_{\geqslant 2}$ that is orthogonal to $\widetilde{\n}_2 \subseteq \g^{(1)}_{\geqslant 0}$. The element $\ad(v)^k f_1$ belongs to $\g^{(1)}_{\geqslant 2 (k-1)}$, which is included in $\g^{(1)}_{\geqslant 2}$ for $k \geqslant 2$. 
    
    It remains to look at $\ad(v)^k f_1$ when $k=1$. We use the decomposition:
    $$ \widetilde{\n}_2 = \g_{\geqslant 2}^{(2)} \oplus \mathfrak{a} \oplus \g_{1,1} \oplus \g_{2,1} \oplus \g_{1,0}. $$
    Because $f_1$ is in $\g_{-2,-2}$ and $\g^{(1)}_{\geqslant 2} \subseteq \g^{(2)}_{\geqslant 1}$, the element $\ad(v) f_1$ belongs to $\g^{(1)}_{\geqslant 0} \cap \g^{(2)}_{\geqslant  -1}$. The only way~$\ad(v) f_1$ may not be orthogonal to $\widetilde{\n}_2$ is when $v$ is in $\g_{2,1}$. But, by construction of $\mathfrak{a}$, one has $\omega_2(v, \mathfrak{a} \oplus \g_{2,1}) = 0$, so $\ad^*(v) f_1$ is orthogonal to $\mathfrak{a}$ and then to $\widetilde{\n}_2$.
    
    This proves that
    $$ \widetilde{\sigma}_2 : \g^{(1)}_1 \longrightarrow \widetilde{\n}_2{}^*, \quad v \longmapsto \widetilde{\chi}_2 + \ad^*(v) \widetilde{\chi}_2 + \frac{1}{2} \ad^*(v)^2 \widetilde{\chi}_2. $$
    is a well-defined surjection. It is clear that the composition $\pi \circ \widetilde\sigma_2 = \widetilde\sigma_1$ because the projection $\widetilde{\n}_2{}^* \twoheadrightarrow \widetilde{\n}_1{}^*$ is $\widetilde N_1$-equivariant. This proves that $\widetilde\sigma_2$ is an isomorphism and that the projection induces an isomorphism~$\widetilde{\Orb}_2 \cong \widetilde{\Orb}_1$.
\end{proof}

\begin{corollary}\label{corollary:orbit-description}
    The orbit $\widetilde{\Orb}_2$ is described by the formula
    $$ \widetilde{\Orb}_2 = \Big\{\sum_{k=0}^2 \frac{1}{k!} \ad^*(v)^k \widetilde{\chi}_2 ~ | ~ v \in \g_{1,0} \oplus \g_{1,1}\Big\} \oplus \ad^*(\g_{1,2}) \widetilde{\chi}_2. $$
    Hence, the fiber of the orbit $\widetilde{\Orb}_2$ by the moment map is
    $$ \widetilde{\pi}_2{}^{-1}\big(\widetilde{\Orb}_2^-\big) = \Big\{- \sum_{k=0}^2 \frac{1}{k!} \ad^*(v)^k \chi_2 ~ | ~ v \in \g_{1,0} \oplus \g_{1,1}\Big\} \oplus \ad^*(\g_{1,2}) \chi_2 \oplus \widetilde{\n}_2{}^\perp. $$
\end{corollary}

\begin{proof}
    Take $v$ in $\g_{1,0} \oplus \g_{1,1}$ and $w$ in $\g_{1,2}$. One has
    $$ \ad^*(v + w)^2 \widetilde{\chi}_2 = \ad^*(v)^2 \widetilde{\chi}_2 + \ad^*(w)^2 \widetilde{\chi}_2 + \ad^*(v) \ad^*(w) \widetilde{\chi}_2 + \ad^*(w) \ad^*(v) \widetilde{\chi}_2. $$
    It is clear that the elements $\ad(w)^2 f_2$, $\ad(v) \ad(w) f_2$ and $\ad(w) \ad(v) f_2$ belong to~$\g^{(1)}_{\geqslant 0} \cap \g^{(2)}_{\geqslant 0}$ and the latter subspace is contained in the orthogonal subspace~$\widetilde{\n}_2{}^{\perp, (\bullet|\bullet)}$, hence we get the equality $\ad^*(v + w)^2 \widetilde{\chi}_2 = \ad^*(v)^2 \widetilde{\chi}_2$.    Then, the corollary follows immediately from Proposition \ref{proposition:orbit-description}.
\end{proof}

We compute the corresponding Hamiltonian reduction by using the fact that the Hamiltonian action of $N_2$ is well-understood and reduction by stages. 

\begin{lemma} \label{lemma:inclusion}
    The inclusion ${\pi_2}^{-1}(\Orb_{2}^-) \subseteq {\widetilde{\pi}_2}^{-1}\big(\widetilde{\Orb}_2^-\big)$ holds.
\end{lemma}

\begin{proof}
    By construction, there is the decomposition $\widetilde{\n}_2 = \n_2 \oplus \g_{1,0}$ of vector spaces. Then, the orthogonal of $\n_2$ for $(\bullet|\bullet)$ is $ {\n_2}^{\perp, (\bullet|\bullet)} = \widetilde{\n}_2{}^{\perp, (\bullet|\bullet)} \oplus \g_{-1,0}$. Using the fact that $\ad(f_1)$ induces an isomorphism $\g_{1,2} \cong \g_{-1,0}$, and the isomorphism~$\g \cong \g^*$ induced by the bilinear form $(\bullet|\bullet)$, we deduce the equality 
    $$ \widetilde{\n}_2{}^\perp \oplus \ad^*(\g_{1,2})\chi_1 = {\n_2}^\perp. $$

    We have the inclusion~$[f_2, \g_{1,2}] \subseteq \g_{-1,0} \oplus \g_{1,0}$ and the adjoint action of $f_2$ is injective on $\g_{1,2}$. Because of the inclusion $ \g_{1,0} \subseteq \widetilde{\n}_2$, we deduce that
    $$ {\n_2}^\perp = \widetilde{\n}_2{}^\perp \oplus \ad^*(\g_{1,2})\chi_1 = \widetilde{\n}_2{}^\perp \oplus \ad^*(\g_{1,2})\chi_2. $$

    One has the equality $$\ad^*(\g_{1,1}) \chi_2 = \ad^*(\g_{1,1}) \chi_1$$ because $[f_0, \g_{1,1}] \subseteq \g_{1,-1} = 0$. We get the explicit descriptions of both fibers:
    \begin{align*}
        \pi_2{}^{-1}(\Orb_2^-) & = - \chi_2 + \ad^*(\g_{1,1}) \chi_2 \oplus \n_2{}^\perp, \\
        \widetilde{\pi}_2{}^{-1}\big(\widetilde{\Orb}_2^-\big) & = \Big\{- \chi_{2} - \ad^*(v) \chi_2 - \ad^*(w) \chi_2 \\ & \quad \quad - \frac{1}{2} \ad^*(v + w)^2 \chi_2  ~ | ~ v \in \g_{1,0}, ~ w \in \g_{1,1} \Big\}\oplus \n_2{}^\perp.
    \end{align*}
    The inclusion is now clear.
\end{proof}

\begin{theorem}\label{theorem:new-construction-slodowy}
    The inclusion of Lemma \ref{lemma:inclusion} induces a map
    $$ {\pi_2}^{-1}(\Orb_{2}^-) /\!/ N_2 \longrightarrow  \widetilde{\pi}_2{}^{-1}\big(\widetilde{\Orb}_2^-\big)  /\!/ \widetilde{N}_2 $$
    which is a Poisson isomorphism.
\end{theorem}

For the proof, we need the following analogue of \eqref{equation:intermediate-stage}, see \cite[Claim~2.4.3]{genra2024reduction}.

\begin{lemma}\label{lemma:intermediate-stage-bis}
    The isomorphism given by Theorem \ref{theorem:slodowy-isomorphism},
    $$ \widetilde{\alpha}_1 : \widetilde{N}_1 \times S_1 \longrightarrow {\widetilde{\pi}_1}^{-1}\big(\widetilde{\Orb}_1^-\big), \quad (g, \xi) \longmapsto \Ad^*(g) \xi,  $$
    restricts to an isomorphism
    $$ \widetilde{N}_1 \times (- \chi_2 + [\g, e_1]^\perp \cap \n_0{}^\perp) \cong \widetilde{\pi}_2{}^{-1}\big(\widetilde{\Orb}_2^-\big). $$
\end{lemma}

\begin{proof}
    There is an inclusion $- \chi_2 + [\g, e_1]^\perp \cap \n_0{}^\perp \subseteq {\widetilde{\pi}_2}^{-1}\big(\widetilde{\Orb}_2^-\big)$ and the right-hand side variety is $\widetilde{N}_1$-stable, hence there is a closed embedding
    $$ \widetilde{\alpha}_1\big(\widetilde{N}_1 \times (- \chi_2 + [\g, e_1]^\perp \cap \n_0{}^\perp)\big) \subseteq {\widetilde{\pi}_2}^{-1}\big(\widetilde{\Orb}_2^-\big). $$
    Both varieties are irreducible, if their dimensions coincide then they are equal and the lemma is proved. The right-hand side has dimension
    \begin{align*}
        \dim (\widetilde{\pi}_2{}^{-1}\big(\widetilde{\Orb}_2^-\big)) & = \dim \g^{(1)}_1 + \dim (\widetilde{\n}_2{})^\perp = \dim \g^{(1)}_1 + \dim \g - \dim \widetilde{\n}_2  \\
        & = \dim \g^{(1)}_1 + \dim \g - \dim \widetilde{\n}_1 - \dim \n_0 
    \end{align*}
    Because $\widetilde{\alpha}_1$ is an isomorphism, we have $ \dim \widetilde{\n}_1 + \dim \g^{e_1} = \dim \g - \dim \widetilde{\n}_1 + \dim \g^{(1)}_1$.
    Hence
    $$ \dim (\widetilde{\pi}_2{}^{-1}\big(\widetilde{\Orb}_2^-\big)) =  \dim \g^{e_1} - \dim \n_0 + \dim \widetilde{\n}_1. $$
    Moreover, we proved in \cite[Claim 2.4.3]{genra2024reduction} that
    $$ \dim \g^{e_1} - \dim \n_0 = \dim [\g, e_1]^\perp \cap \n_0{}^\perp, $$
    whence 
    $$  \dim \g^{e_1} - \dim \n_0 + \dim \widetilde{\n}_1 = \dim \widetilde{\n}_1 + \dim [\g, e_1]^\perp \cap \n_0{}^\perp. $$
    This is the dimension of the left-hand side $\widetilde{\alpha}_1\big(\widetilde{N}_1 \times (-\chi_2 + [\g, e_1]^\perp \cap \n_0{}^\perp)\big)$, hence
    $$ \dim \widetilde{\alpha}_1\big(\widetilde{N}_1 \times (-\chi_2 + [\g, e_1]^\perp \cap \n_0{}^\perp)\big) = \dim (\widetilde{\pi}_2{}^{-1}\big(\widetilde{\Orb}_2\big)). $$
\end{proof}

\begin{proof}[Proof of Theorem \ref{theorem:new-construction-slodowy}]
    Consider the affine GIT quotients:
    \begin{align*}
        {\pi_{2}}^{-1}(\Orb_{2}^-)  /\!/ N_{2} & = \Spec \C[{\pi_{2}}^{-1}(\Orb_{2}^-)]^{N_{2}}, \\
        {\widetilde{\pi}_2}^{-1}\big(\widetilde{\Orb}_2^-\big)  /\!/ \widetilde{N}_2 & = \Spec \C\big[{\widetilde{\pi}_2}^{-1}\big(\widetilde{\Orb}_2^-\big)\big]^{\widetilde{N}_2}.
    \end{align*}
    We have a $\C$-algebra map
    $$ \C\big[\widetilde{\pi}_2{}^{-1}\big(\widetilde{\Orb}_2^-\big)\big]^{\widetilde{N}_2} \longrightarrow \C[{\pi_{2}}^{-1}(\Orb_{2}^-)]^{N_2}, $$
    which clearly respect the Poisson structures on both sides.
    
    The semi-direct product decompositions $N_2 = N_1 \rtimes N_0$ and~$\widetilde{N}_2 = \widetilde{N}_1 \rtimes N_0$ imply:
    \begin{align*}
\C\big[{\widetilde{\pi}_2}^{-1}\big(\widetilde{\Orb}_2^-\big)\big]^{\widetilde{N}_2} & = \Big(\C\big[{\widetilde{\pi}_2}^{-1}\big(\widetilde{\Orb}_2^-\big)\big]^{\widetilde{N}_1}\Big)^{N_0}, \\
         \C[{\pi_{2}}^{-1}(\Orb_{2}^-)]^{N_2}& = \big(\C[{\pi_{2}}^{-1}(\Orb_{2}^-)]^{N_1}\big)^{N_0}.
    \end{align*}

    We apply Isomorphism \eqref{equation:intermediate-stage} and Lemma \ref{lemma:intermediate-stage-bis}: the algebra map
    $$ \C\big[{\widetilde{\pi}_2}^{-1}\big(\widetilde{\Orb}_2^-\big)\big]^{\widetilde{N}_1} \longrightarrow \C[{\pi_{2}}^{-1}(\Orb_{2}^-)]^{N_1}. $$
    is a $N_0$-isomorphism. After taking the $N_0$-invariants, we get the desired isomorphism.
\end{proof}

\begin{corollary} \label{corollary:vanishing-lie-cohomology}
    The Lie algebra cohomology of $\C\big[{\Jet_\infty \widetilde{\pi}_2}^{-1}\big(\widetilde{\Orb}_2^-\big)\big]$ as $\widetilde{\n}_2[t]$-module is given by the superalgebra isomorphism
    $$ \Hgy^\bullet\Big(\widetilde{\n}_2[t], \C\big[{\Jet_\infty \widetilde{\pi}_2}^{-1}\big(\widetilde{\Orb}_2^-\big)\big]\Big) \cong \delta_{\bullet = 0} \, \C[\Jet_\infty S_2], $$
    where is the isomorphism is induced by the action map
    $$ \widetilde{N}_2 \times S_2 \longrightarrow \widetilde{\pi}_2{}^{-1} \big(\widetilde{\Orb}_2^-\big), \quad (g, \xi) \longmapsto \Ad^*(g) \xi. $$
\end{corollary}

\begin{proof}
    For $m$ a nonnegative integer, consider the $m$-jets. We can compute~$ \Hgy^\bullet\Big(\widetilde{\n}_2[t] / (t^{m+1}), \C\big[{\Jet_m \widetilde{\pi}_2}^{-1}\big(\widetilde{\Orb}_2^-\big)\big]\Big) $ by using the Hochschild--Serre spectral sequence $\{E_r\}_{r = 0}^\infty$ corresponding to the decomposition $ \widetilde{N}_2 = \widetilde{N}_1 \rtimes N_0 $. 
    
    The $1$-st page is given by
    $$ E_1^{p, q} = \Hgy^q\Big(\widetilde{\n}_1[t] / (t^{m+1}), \C\big[{\Jet_m \widetilde{\pi}_2}^{-1}\big(\widetilde{\Orb}_2^-\big)\big]\Big) \otimes_\C \Alt^p\big(\big(\n_0[t] / (t^{m+1})\big)^*\big). $$
    The $2$-nd page is given by
    $$ E_2^{p, q} = \Hgy^p\bigg(\n_0[t] / (t^{m+1}), \Hgy^q\Big(\widetilde{\n}_1[t] / (t^{m+1}), \C\big[{\Jet_m \widetilde{\pi}_2}^{-1}\big(\widetilde{\Orb}_2^-\big)\big]\Big)\bigg). $$
    The infinity page is
    $$ E_\infty^{p, q} = \gr^p \Hgy^{p+ q}\Big(\widetilde{\n}_2[t] / (t^{m+1}), \C\big[{\Jet_m \widetilde{\pi}_2}^{-1}\big(\widetilde{\Orb}_2^-\big)\big]\Big), $$
    where the implicit filtration is the one induced by the double complex structure. The double complex is bounded so the spectral sequence is convergent.

    \Cref{lemma:intermediate-stage-bis} gives the following isomorphism of $\widetilde{\n}_1[t] / (t^{m+1})$-module:
    $$ \C\big[{\Jet_m \widetilde{\pi}_2}^{-1}\big(\widetilde{\Orb}_2^-\big)\big] \cong \C\big[\Jet_m \widetilde{N}_1\big] \otimes_\C \C[\Jet_m(-\chi_2 + [\g, e_1]^\perp \cap {\n_0}^\perp)]. $$
    The Lie algebra cohomology of $\C\big[\Jet_m \widetilde{N}_1\big]$ is the algebraic de Rham cohomology of~$\Jet_m \widetilde{N}_1$, and because this groups is unipotent, its cohomology is $\C$ in degree $0$ and zero otherwise.

    Hence the first page component $E_1^{p, q}$ is zero when $q \neq 0$ and
    $$ E_1^{p, 0} \cong \C\big[\Jet_m \pi_0{}^{-1}(-\overline{\chi}_0)\big]\otimes_\C \Alt^p\big(\big(\n_0[t] / (t^{m+1})\big)^*\big) $$
    because of the $N_0$-isomorphism $\widetilde{\pi}_2{}^{-1}\big(\widetilde{\Orb}_2^-\big) /\!/ \widetilde{N}_1 \cong \pi_0{}^{-1}(-\overline{\chi}_0)$. Applying \Cref{corollary:isomorphism-second-stage}, we can use the $N_0$-isomorphism $ N_0 \times S_2 \cong \pi_0{}^{-1}(-\overline{\chi}_0) $ and the same argument as previously to say that the second page component $E_2^{p, q}$ is zero except if $p = q = 0$, so the spectral sequence collapses. 

    Then we deduce the isomorphism
    $$ \Hgy^0\Big(\widetilde{\n}_2[t] / (t^{m+1}), \C\big[{\Jet_m \widetilde{\pi}_2}^{-1}\big(\widetilde{\Orb}_2^-\big)\big]\Big) \cong \C[\Jet_m S_2], $$
    and the other cohomology groups are trivial.

    Consider the Lie algebra cochain complex
    $$ C_m \defeq \C\big[{\Jet_m \widetilde{\pi}_2}^{-1}\big(\widetilde{\Orb}_2^-\big)\big] \otimes_\C \Alt^\bullet\big(\widetilde{\n}_2[t] / (t^{m+1})\big)^* $$
    for all $m \geqslant 0$. The sequence $\{C_m\}_{m \geqslant 0}$ form an inductive sequence and its colimit is the Lie algebra cochain complex of the $\widetilde{\n}_2[t]$-module $\C\big[{\Jet_\infty \widetilde{\pi}_2}^{-1}\big(\widetilde{\Orb}_2^-\big)\big]$:
    $$ C_\infty \defeq \C\big[{\Jet_\infty \widetilde{\pi}_2}^{-1}\big(\widetilde{\Orb}_2^-\big)\big] \otimes_\C \Alt^\bullet_\infty\big(\widetilde{\n}_2{}^*\big). $$

    Because taking the cohomology commutes with inductive colimits, one gets
    $$\Hgy^0\Big(\widetilde{\n}_2[t], \C\big[{\Jet_\infty \widetilde{\pi}_2}^{-1}\big(\widetilde{\Orb}_2^-\big)\big]\Big) \cong \colim_{m \geqslant 0} \C[\Jet_m S_2] \cong \C[\Jet_\infty S_2]. $$

    \begin{remark}
        We conjecture that the action map
        $$ \widetilde{N}_2 \times S_2 \longrightarrow \widetilde{\pi}_2{}^{-1} \big(\widetilde{\Orb}_2^-\big), \quad (g, \xi) \longmapsto \Ad^*(g) \xi $$
        is an isomorphism. Such isomorphism should be used to apply \Cref{main:vanishing-vertex}. In fact, \Cref{corollary:vanishing-lie-cohomology} will be enough for our purpose.
    \end{remark}

    % Recall the discussion of Paragraph \ref{paragraph:git-arc}: the action of the Lie algebra $\widetilde{\n}_2[t] / (t^{m+1})$ can be lifted to an action of the Lie algebra $\widetilde{\n}_2[t]$ such that the normal subalgebra $t^{m+1} \widetilde{\n}_2[t]$ acts trivially. We claim that
    % $$ \Hgy^\bullet\Big(\widetilde{\n}_2[t] / (t^{m+1}), \C\big[{\Jet_m \widetilde{\pi}_2}^{-1}\big(\widetilde{\Orb}_2\big)\big]\Big) \cong \Hgy^\bullet\Big(\widetilde{\n}_2[t], \C\big[{\Jet_m \widetilde{\pi}_2}^{-1}\big(\widetilde{\Orb}_2\big)\big]\Big). $$

    % Hence, for every nonnegative $m$, one has
    % $$ \Hgy^0\Big(\widetilde{\n}_2[t], \C\big[{\Jet_m \widetilde{\pi}_2}^{-1}\big(\widetilde{\Orb}_2\big)\big]\Big) \cong \C[S_2] $$
    % and the other cohomology groups are trivial. Taking the colimit (which commutes with the cohomology), the corollary follows.
\end{proof}

\subsection{New construction of the second affine W-algebra} \label{subsection:new-definition-w-algebra}

The Lie algebra $\n_2$ corresponds to the choice of the isotropic subspace of $\l_2$ in the construction of the affine W-algebra $\W^k(\g, f_2)$ described in Subsection~\ref{subsection:definition-w-algebras}. Because of the isomorphism~$\g_{1,1} \cong \l_2{}^{\perp, \omega_2} / \l_2$, we have
$$ \W^k(\g, f_2) = \Hgy^0(\Cpx^\bullet_2, \mathcal{d}_2), \quad \text{where} \quad \Cpx^\bullet_2 \defeq \V^k(\g) \otimes_\C \A(\g_{1,1}) \otimes_\C \F^\bullet(\n_2 \oplus \n_2{}^*). $$

Geometrically, the complex $(\Cpx^\bullet_2, \mathcal{d}_2)$ corresponds to the moment map
$$ \mu_2 : \g^* \times \g_{1,1}  \longrightarrow \n_2{}^*, \quad (\xi, v)  \longmapsto \pi_2(\xi) + \overline{\chi}_{2} + \ad^*(v)\overline{\chi}_{2},  $$
where the acting group is $N_2$. One has the $N_2$-equivariant isomorphism
$$ \mu_2{}^{-1}(0) \cong {\pi_{2}}^{-1}(\Orb_{2}^-). $$ 

The complex $\big(\widetilde\Cpx^\bullet_2, \widetilde{\mathcal{d}}_2\big)$ corresponds to the moment map
$$\widetilde{\mu}_2 : \g^* \times \g^{(1)}_{1} \longrightarrow \widetilde{\n}_2{}^*, \quad (\xi, v)  \longmapsto \pi_2(\xi) + \widetilde{\chi}_{2} + \ad^*(v)\widetilde{\chi}_{2} +  \frac{1}{2}\ad^*(v)\widetilde{\chi}_{2},$$
where the acting group is $\widetilde{N}_2$.  One has the $\widetilde{N}_2$-equivariant isomorphism
$$ \widetilde{\mu}_2{}^{-1}(0) \cong \widetilde{\pi}_{2}{}^{-1}\big(\widetilde{\Orb}_{2}^-\big). $$
Denote by $\big(\widetilde{C}_{2, \infty}^\bullet, \widetilde{d}_{2, \infty}\big)$ the Poisson vertex BRST complex associated with the arc space of this moment map.

\begin{proposition} \label{proposition:vanishing-new-complex}
    The following cohomology vanishes outside degree $0$: 
    $$ \Hgy^n\big(\widetilde\Cpx^\bullet_2, \widetilde{\mathcal{d}}_2\big) = 0 \quad \text{for} \quad n \neq 0, $$ 
    and there is a natural isomorphism 
    $$ \gr_{\Fil}\Hgy^0\big(\widetilde\Cpx^\bullet_2, \widetilde{\mathcal{d}}_2\big) \overset{\sim}{\longrightarrow} \Hgy^0\big(\widetilde{C}_{2, \infty}^\bullet, \widetilde{d}_{2, \infty}\big), $$ 
    where the filtration $\Fil$ on the cohomology $\Hgy^0\big(\widetilde\Cpx^\bullet_2, \widetilde{\mathcal{d}}_2\big)$ is induced by the Li filtration on the complex $\big(\widetilde\Cpx^\bullet_2, \widetilde{\mathcal{d}}_2\big)$.
\end{proposition}

\begin{proof}
    Consider the Hamiltonian operator $ \widetilde{\Ham}^\old_2 \defeq \Ham^\g + L^\A_0$ on the vertex algebra~$\V^k(\g) \otimes_\C \A(\g^{(1)}_1)$, it induces a nonnegative grading. The homogeneous subspaces are all finite dimensional and for any $x \in \widetilde{\n}_2$, the image $\widetilde{\Upsilon}_2(a)$ is the sum of terms of degrees less than one. The standard comoment map is given by:
    $$ \widetilde{\Upsilon}_{2, \mathrm{st}}(x) \defeq \begin{cases}
        & x + \frac{1}{2}\sum_{i=1}^{2 s} \NO{\psi^i \psi_{[v_i, x]}}  \quad   \text{if} \quad x \in \n_0, \\
        & x  \quad   \text{if} \quad  x \in \widetilde{\n}_1.
    \end{cases} $$
    The induced action of $\mathcal{U}(\widetilde{\n}_2[t^{-1}]t^{-1})$ is free because $\V^k(\g) \otimes_\C \A(\g^{(1)})$ is freely generated by any basis of $\g$ and any basis of $\g^{(1)}$. Hence \Cref{theorem:quotient-quasi-isomorphic} applies and the cohomology of the complex $\big(\widetilde\Cpx^\bullet_2, \widetilde{\mathcal{d}}_2\big)$ is isomorphic to the cohomology of the corresponding quotient complex $\big(\widetilde\Cpx^\bullet_{2, +}, \widetilde{\mathcal{d}}_{2,+}\big)$.

    Let $\l$ be a Lagrangian subspace of $\g_{1,1}$ and $\l^{\mathrm{c}}$ be a Lagrangian complement such that $\l$ and $\l^{\mathrm{c}}$ are perfectly paired by the symplectic form. Then $\g_{1,0} \oplus \l$ and~$\l^{\mathrm{c}} \oplus \g_{1,2}$ are perfectly paired by the symplectic form $\omega_1$. 
    
    Denote by $\{v_i\}_{i=1}^{2s}$ a basis of $\g^{(1)}_{1}$ such that $v_i = v^{2s - i + 1}$ for all $1 \leqslant i \leqslant s$,
    \begin{align*}
        \Span_\C \{v_i\}_{i=1}^{s_0} = \g_{1,0}, \quad & \Span_\C \{v_i\}_{i=s_0 + 1}^{s} = \l, \\
        \Span_\C \{v_i\}_{i=2s -s_0}^{s + 1} = \l^{\mathrm c}, \quad & \Span_\C \{v_i\}_{i=2s - s_0+1}^{2s} = \g_{1,2}.
    \end{align*}
    For $1 \leqslant i \leqslant s_0$, set $a_i \defeq -\frac{1}{2}$ and for $s_0 + 1 \leqslant i \leqslant s$, set $a_i \defeq 0$.
    
    Choose a basis $\{x_i\}_{i \in I(\widetilde{\n}_2)}$  of $\widetilde{\n}_2$ which is homogeneous for the $H_2$-grading and denote by $\delta(x_i)$ the degree of $x_i$. Consider the Hamiltonian operator
    $$ \widetilde{\Ham}^\new_2 = \Ham^\g + \frac{1}{2}(\partial {H_2})_{(1)} + L^\A(a_\bullet)_0 + L^\F(m_\bullet)_0 \quad \text{where} \quad m_\bullet = \frac{1}{2}\delta(x_\bullet). $$
    The degrees of the strong generators are
    \begin{align*}
        &\Delta(x) = 1-\frac{\delta}{2} \quad \text{for} \quad x \in \g^{(2)}_\delta \quad \text{and} \quad \delta \in \Z, \\
        &\Delta(\psi_v) = 1 -\frac{\delta}{2} \quad \text{for}\quad v \in \g _{1,\delta} \quad \text{and} \quad \delta \in \{0,1,2\}, \\
        &\Delta(\phi_i) = 1 -\frac{\delta(x_i)}{2} \quad \text{and} \quad  \Delta(\phi_i^*) =\frac{\delta(x_i)}{2}  \quad \text{for} \quad i \in I(\widetilde{\n}_2).
    \end{align*}
    The element $\widetilde{Q}_2$ is homogeneous of degree $1$ for $\widetilde{\Ham}^\new_2$ and the quotient complex~$\widetilde\Cpx_{2, +}$ has an induced grading which is nonnegative. We can apply Theorem \ref{theorem:vanishing-vertex} by replacing the assumption \eqref{theorem:vanishing-vertex,item:section} by the vanishing result of Corollary \ref{corollary:vanishing-lie-cohomology}, which is enough, and we are done.    
    \end{proof}

We want to build an intermediary complex to compare $(\Cpx^\bullet_2, \mathcal{d}_2)$ and $\big(\widetilde\Cpx^\bullet_2, \widetilde{\mathcal{d}}_2\big)$, by analogy with Theorem \ref{theorem:equivalent-constructions}. To do so, we need the following lemma.

\begin{lemma} \label{lemma:intermediary-chiral-comoment}
    The image of $\V(\n_2)$ by the chiral comoment map $\widetilde{\Upsilon}_2$ is contained in the vertex subalgebra~$\V^k(\g) \otimes_\C \A(\g_{1,1} \oplus \g_{1,2})$ of $\V^k(\g) \otimes_\C \A(\g^{(1)}_1)$. 
\end{lemma}

\begin{proof}
    We just need to check it for a strong generator $x$ in $\n_2 = \n_1 \oplus \n_0$. If $x$ is in $\n_1$, it is clear that $\widetilde{\Upsilon}_2(x)$ belongs to $\V^k(\g) \otimes_\C \A(\l_1{}^{\perp, \omega_1})$. For $x$ in $\n_0$, we need to check that the quadratic term $ \sum_{i=1}^{2 s} \NO{\psi^i \psi_{[v_i, x]}}$ belongs to~$\A(\g_{1,1} \oplus \g_{1,2})$. 
    
    According to the proof of \Cref{lemma:construction-L1}, there is a basis $\{v_i\}_{i=1}^{2s}$ of $\g^{(1)}_1$ such that 
    $$ \Span_\C \{v_i\}_{i=1}^{s_0} = \g_{1,0}, \quad \{v_i\}_{i=s_0 + 1}^{2s - s_0} = \g_{1,1}, \quad \quad \Span_\C \{v_i\}_{i=2s - s_0+1}^{2s} = \g_{1,2}$$
    and $v_i = v^{2s - i + 1}$ for all~$1 \leqslant i \leqslant s_0$. In particular, $\l_1{}^{\perp, \omega_1} = \g_{1,1} \oplus \g_{1,2}$ is spanned by $\{v^i\}_{i=1}^{2s - s_0}$.

    The conditions \eqref{conditions} imply the inclusion $[\g_{1,2}, \m_0] \subseteq \g^{(1)}_{1} \cap \g^{(2)}_{\geqslant 3} = 0$, hence
    $$ \sum_{i=1}^{2 s} \NO{\psi^i \psi_{[v_i, x]}} = \sum_{i=1}^{2s - s_0} \NO{\psi^i \psi_{[v_i, x]}} \in \A(\g_{1,1} \oplus \g_{1,2}). $$
\end{proof}

\begin{remark}
    The last lemma is coherent with the description of the orbit $\widetilde{\Orb}_2$ given in Proposition~\ref{corollary:orbit-description}.
\end{remark}

By restriction, we have a chiral comoment map
$$ \Upsilon_{2, \interm}: \V(\n_2) \longrightarrow \V^k(\g) \otimes_\C \A(\g_{1,1} \oplus \g_{1,2}). $$
As a consequence of the previous lemmas, we can define an intermediary BRST complex, denoted by $(\Cpx_{2,\interm}^\bullet, \mathcal{d}_{2,\interm})$, where $\Cpx_{2,\interm}$ is the subalgebra of $\widetilde{\Cpx}_2$ defined by
$$ \Cpx_{2,\interm}^\bullet \defeq \V^k(\g) \otimes_\C \A(\g_{1,1} \oplus \g_{1,2}) \otimes_\C \F^\bullet(\n_2 \oplus \widetilde{\n}_2{}^*), $$
and the operator $\mathcal{d}_{2, \interm}$ is the restriction of the coboundary operator $\mathcal{d}_2$ to the subalgebra~$\Cpx_{2,\interm}$.

By construction, one has there is a vertex algebra embedding map
$$ \Pi_1 : (\Cpx_{2, \interm}^\bullet, \mathcal{d}_{2, \interm}) \hooklongrightarrow \big(\widetilde{\Cpx}_2^\bullet, \widetilde{\mathcal{d}}_2\big) $$
and it commutes with the differentials because of Lemma \ref{lemma:intermediary-chiral-comoment}.

\begin{lemma}
    The projection maps $\g_{1,1} \oplus \g_{1,2} \twoheadrightarrow \g_{1,1}$ and $\widetilde{\n}_2{}^* \twoheadrightarrow \n_2{}^*$ induce an surjective map
    $$  \Pi_2 : (\Cpx_{2, \interm}^\bullet, \mathcal{d}_{2, \interm}) \twoheadlongrightarrow (\Cpx_2^\bullet, \mathcal{d}_2) $$
    of cochain complexes and vertex superalgebras.
\end{lemma}

\begin{proof}
    We use the same notations as in the proof of Lemma \ref{lemma:intermediary-chiral-comoment}. The projection is clearly a vertex algebra map. For $x$ in $\n_0$, the quadratic term~$\sum_{i=1}^{2s - s_0} \NO{\psi^i \psi_{[v_i, x]}}$ in $\Upsilon(x)$ is killed by the projection. 
    
    Indeed, if $1 \leqslant i \leqslant s_0$, then $\psi^i$ belongs to $\g_{1,2}$, and if $s_0 + 1 \leqslant i \leqslant 2s - s_0$, then $v_i \in \g_{1,1}$ and $[v_i, \n_0]\subseteq \g^{(1)}_1 \cap \g^{(1)}_{\geqslant 2} = \g_{1,2}$. 
    
    Hence the following diagram of vertex algebra homomorphisms commutes:
    $$ \begin{tikzcd}[row sep = small]
        & \V^k(\g) \otimes_\C \A(\g_{1,1} \oplus \g_{1,2}) \arrow[dd, two heads] \\
        \V(\n_2) \arrow[ru, "\Upsilon_{2, \interm}"] \arrow[rd, "\Upsilon_2"'] & \\
        & \V^k(\g) \otimes_\C \A(\g_{1,1}).
    \end{tikzcd} $$
    It implies that the projection $\Pi_2$ commutes with the differential.
\end{proof}

Let us denote by $\Orb_{2, \interm} \defeq \Ad^*\big(\widetilde{N}_2\big)\overline{\chi}_2$ the $\widetilde{N}_2$-orbit in $\n_2{}^*$ of the restriction of the linear form $\chi_2$ to $\n_2$. 

\begin{lemma} \label{lemma:intermediary-orbit}
    The orbit $\Orb_{2, \interm}$ is described by the formula
    $$ \Orb_{2, \interm} = \Big\{\sum_{k=0}^2 \frac{1}{k!} \ad^*(v)^k \widetilde{\chi}_2 ~ | ~ v \in \g_{1,0} \oplus \g_{1,1}\Big\}. $$
\end{lemma}

\begin{proof}
    The projection map $\widetilde{\n}_2{}^* \twoheadrightarrow \n_2{}^*$ is $\widetilde{N}_2$-equivariant and the orbit $\Orb_{2, \interm}$ is the image of $\widetilde{\Orb}_2$. Using the explicit description of the latter orbit in Corollary \ref{corollary:orbit-description} and the inclusion $\ad^*(\g_{1,2})\chi_2 \subseteq {\n_2}^\perp$, the lemma follows.
\end{proof}

The intermediary moment map is defined as
\begin{align*}
    \mu_{2, \interm} : \g^* \times (\g_{1,0} \oplus \g_{1,1}) & \longrightarrow \n_2{}^*, \\
    (\xi, v) & \longmapsto \pi_2(\xi) + \overline{\chi}_{2} + \ad^*(v)\overline{\chi}_{2} + \frac{1}{2} \ad^*(v)^2 \overline{\chi}_{2},
\end{align*}
it is a $\widetilde{N}_2$-equivariant homomorphism because $N_2$ is a normal subgroup of $\widetilde{N}_2$.

\begin{lemma} \label{lemma:intermediary-isomorphisms-f2}
    The projection map and the embedding
    $$  \g^* \times \g^{(1)}_1  \twoheadlongrightarrow \g^* \times (\g_{1,0} \oplus \g_{1,1})  \quad \text{and} \quad  \g^* \times (\g_{1,1}) \hooklongrightarrow \g^* \times (\g_{1,0} \oplus \g_{1,1})  $$
    induce Poisson isomorphisms between the corresponding Hamiltonian reductions:
    $$ \mu_2{}^{-1}(0) /\!/ N_2 \cong \mu_{2,\interm}{}^{-1}(0) /\!/ \widetilde{N}_2 \cong \widetilde{\mu}_2{}^{-1}(0) /\!/ \widetilde{N}_2. $$
\end{lemma}

\begin{proof}
    The projection map $\g^* \times \g^{(1)}_1  \twoheadrightarrow \g^* \times (\g_{1,0} \oplus \g_{1,1})$ induces an $\widetilde{N}_2$-equivariant isomorphism~$\widetilde{\mu}_2{}^{-1}(0) \cong \mu_{2,\interm}{}^{-1}(0)$. It follows from the parametrization of the intermediary orbit $\Orb_{2,\interm}$ given in \Cref{lemma:intermediary-orbit}. Then the lemma follows from \Cref{theorem:new-construction-slodowy}.
\end{proof}

\begin{theorem} \label{theorem:new-construction-w-algebra}
    These inclusion and projection maps of cochain complexes
    $$ \Pi_1 : (\Cpx_{2, \interm}^\bullet, \mathcal{d}_{2, \interm}) \hooklongrightarrow \big(\widetilde\Cpx^\bullet_2, \widetilde{\mathcal{d}}_2\big) \quad \text{and} \quad \Pi_2 : (\Cpx_{2, \interm}^\bullet, \mathcal{d}_{2,\interm}) \twoheadlongrightarrow (\Cpx_2^\bullet, \mathcal{d}_2) $$
    induce vertex algebra isomorphisms between their cohomologies:
    $$ \Hgy^0(\Cpx_2^\bullet, \mathcal{d}_2) \cong  \Hgy^0(\Cpx_{2, \interm}^\bullet, \mathcal{d}_{2,\interm}) \cong \Hgy^0\big(\widetilde\Cpx^\bullet_2, \widetilde{\mathcal{d}}_2\big). $$
\end{theorem}

\begin{proof}
    The proof relies on the same arguments as the proof of \Cref{theorem:equivalent-vertex}.
\end{proof}

\subsection{Proof of reduction by stages for W-algebras} \label{subsection:end-proof-reduction-stages}

Recall the BRST cochain complex $(\Cpx_0, \mathcal{d}_0)$ introduced in Theorem \ref{theorem:reduction-by-stages-w-algebras}. Geometrically, this complex corresponds to the moment map
$$ \mu_0 : \widetilde{\pi}_1{}^{-1}\big(\widetilde{\Orb}_1\big) /\!/ \widetilde{N}_1  \longrightarrow \n_0{}^*, \quad [\xi] \longmapsto \widetilde{\pi}_2(\xi) + \overline{\chi}_{0}, $$
where the acting group is $N_0$.  One has the $N_2$-equivariant isomorphism
$$ \mu_0{}^{-1}(0) \cong \pi_{0}{}^{-1}(-\overline{\chi}_0). $$
Denote by $(C_{0, \infty}^\bullet, d_{0, \infty})$ the Poisson vertex BRST complex associated with the arc space of this moment map.

\begin{proposition} \label{proposition:vanishing-second-stage}
    The following cohomology vanishes outside degree $0$:
    $$\Hgy^n(\Cpx_0, \mathcal{d}_0) \quad \text{for} \quad n \neq 0, $$
    and there is a natural isomorphism
    $$ \gr_{\Fil}\Hgy^0(\Cpx_0, \mathcal{d}_0) \overset{\sim}{\longrightarrow} \Hgy^0(C_{0, \infty}^\bullet, d_{0, \infty}), $$ 
    where the filtration $\Fil$ on the cohomology $\Hgy^0(\Cpx_0, \mathcal{d}_0)$ is induced by the Li filtration on the complex $(\Cpx_0, \mathcal{d}_0)$.
\end{proposition}

\begin{proof}
    We equip $\W^k(\g, f_1)$ with the Hamiltonian operator~$\Ham_0^\old$~induced by the operator defined in equation \eqref{equation:hamiltonian-w-algebra}. Because of Theorem \ref{theorem:basis-W-algebra}, we can provide a strong basis $\{J^{\{i\}}\}_{i=0}^\ell$ of the W-algebra corresponding to a $H_1$-homogeneous basis $\{x_i\}_{i=0}^\ell$ of $\g^{f_1}$. Each element $J^i$ is $\Ham_0^\old$-homogeneous of degree $\Delta(J^{\{i\}}) = 1 - \frac{1}{2} \delta(x_i)$, where~$\delta(x_i)$ is the degree of $x_i$. Because $\delta(x_i)$ is nonpositive, $\Delta(J^{\{i\}})$ is positive. This implies that each homogeneous component of the W-algebra $\W^k(\g, f_1)$ is finite dimensional. 
    
    For $x$ in $\n_0$, the image $\Upsilon_0(x)$ is the sum of terms of degrees less than one. The standard comoment map is given by $\Upsilon_{0, \mathrm{st}}(x) = J^{\{x\}}.$ The induced action of $\mathcal{U}(\n_0[t^{-1}]t^{-1})$ is free because the affine W-algebra is freely generated. Hence, \Cref{theorem:quotient-quasi-isomorphic} holds and the cohomology of the complex $(\Cpx_0, \mathcal{d}_0)$ is isomorphic to the cohomology of the corresponding quotient complex $(\Cpx_{0,+}, \mathcal{d}_{0,+})$. 

    Let $\{x_i\}_{i \in I(\n_0)}$ be a basis of $\n_0$ which is homogeneous for the $H_2$-grading and denote by $\delta(x_i)$ the degree of $x_i$. Let $H_0 \defeq H_2 - H_1$ which belongs to the Lie algebra $\g^{\natural, 1}$. Define the following Hamiltonian operator on the complex $\Cpx_0$: 
    $$ \Ham_0^\new = \Ham_0^\old + \frac{1}{2} \big(\partial J^{\{H_0\}}\big)_{(1)} + L^\F(m_\bullet), \quad \text{where} \quad m_\bullet = \frac{1}{2} \delta(x_\bullet). $$
    The degrees of the strong generators are
    \begin{align*}
        &\Delta (J^{\{x\}}) = 1-\frac{\delta}{2} \quad \text{for} \quad x \in \n_0 \cap\g^{(2)}_\delta \quad \text{and} \quad \delta \in \Z, \\
        &\Delta(\phi_i) = 1 -\frac{\delta(x_i)}{2} \quad \text{and} \quad \Delta (\phi_i^*)=\frac{\delta(x_i)}{2},\quad \text{for} \quad i \in I(\widetilde{\n}_0).
    \end{align*}
    The element $Q_0$ is homogeneous of degree $1$ and the quotient complex $\Cpx_{0, +}$ has an induced grading which is nonnegative. Moreover, the moment map is smooth (see the explanation after \cite[Claim~2.4.2]{genra2024reduction}) and we have the isomorphism given by Corollary \ref{corollary:isomorphism-second-stage}. So we can apply Theorem \ref{theorem:vanishing-vertex} and the proof is done.
\end{proof}

\begin{proposition}[{\cite[Theorem 4.17]{arakawa2015associated}}] \label{proposition:li-filtration-w-algebra}
    The filtration induced on the W-algebra $\W^k(\g, f_1)$ by the Li filtration on the BRST cochain complex $\widetilde{\Cpx}_1$ coincides with its Li filtration: for all $p$ in $\Z$, $\FilLi^p \W^k(\g, f_1) = \Hgy^0\big(\FilLi^p \widetilde{\Cpx}_1^\bullet, \mathcal{d}_1\big)$.
\end{proposition}

Now we can prove our main theorem, that is to say reduction by stages.

\begin{proof}[Proof of Theorem \ref{theorem:reduction-by-stages-w-algebras}]
    Recall the map $\Theta : \Hgy^0(\Cpx_0^\bullet, \mathcal{d}_0) \rightarrow \Hgy^0\big(\widetilde\Cpx^\bullet_2, \widetilde{\mathcal{d}}_2\big)$ which is defined in~\eqref{equation:natural-map}. This map respects the Hamiltonian gradings defined on both sides, which are nonnegative. We consider the filtrations on the cohomology induced by the Li filtrations on the cochain complexes. These filtrations are finite on each homogeneous part of the gradings, so the map $\Theta$ will be an isomorphism if the associated graded map $\gr \Theta$ is. 
    
    On the left-hand side, using \Cref{proposition:vanishing-second-stage,proposition:li-filtration-w-algebra}, there is the isomorphism
    $$ \gr_{\Fil} \Hgy^0(\Cpx_0^\bullet, \mathcal{d}_0) \cong \C[\Jet_{\infty}(\pi_0{}^{-1}(-\overline{\chi}_0) /\!/ N_0)]. $$
    On the right-hand side, according to \Cref{proposition:vanishing-new-complex}, there is the isomorphism
    $$ \gr_{\Fil} \Hgy^0\big(\widetilde\Cpx^\bullet_2, \widetilde{\mathcal{d}}_2\big) \cong \C\big[\Jet_\infty\big(\widetilde{\pi}_{2}{}^{-1}\big(\widetilde{\Orb}_{2}^-\big) /\!/ \widetilde{N}_2\big) \big]. $$
    Hence, $\gr \Theta$ coincides with the Poisson vertex isomorphism induced by the reduction by stages for Slodowy slices: $\widetilde{\pi}_{2}{}^{-1}\big(\widetilde{\Orb}_{2}^-\big) /\!/ \widetilde{N}_2 \cong \pi_0{}^{-1}(-\overline{\chi}_0) /\!/ N_0$.
\end{proof}

\subsection{Examples of reduction by stages} \label{subsection:examples}

\subsubsection{Examples in type $\mathrm A$}

Let $n$ be an integer greater than $2$ and let $\g \defeq \sl_n$ be the Lie algebra of traceless square matrices of size $n$. Denote by $E_{i,j}$ the elementary matrices. We recall that the nilpotent orbits of $\sl_n$ are in bijection with partitions of $n$ \cite[Chapter 5]{collingwood1993nilpotent}. A partition of $n$ is a nonincreasing sequences of positive integers~$(a_1, \dots, a_r)$ such that $\sum_{i=1}^r a_i = n$. 

They can be represented by Young diagrams, that we draw in the $xy$-plane so that the rows are parallel to the $x$-axis. We assume that the box size is twice the unit length in this coordinate system. Good gradings of nilpotent elements can be classified using \emph{pyramids}, see \cite{elashvili2005classification, brundan2007good} for details. These pyramids are built by shifting the rows in a Young diagram and fixing a bijective labeling of the boxes by numbers between $1$ and $n$.

Consider two partitions of $n$ the form 
$$ a_\bullet^{(1)} \defeq (a_1, \dots, a_r, a, 1^b) \quad \text{and} \quad a_\bullet^{(2)} \defeq (a_1, \dots, a_r, a + 1, 1^{b-1}). $$ 
These partitions start with the same element $a_1, \dots, a_r$ and finish with two different hook-type partitions, $(a, 1^b)$ and $(a + 1, 1^{b-1})$.

We associate to theses partitions the pyramids $P_1$ and $P_2$ as follows. Both pyramids start with the same first $r$ rows of length $(a_1, \dots, a_r)$, that we draw left-aligned and with the same labels. The hook-type part $(a, 1^b)$ of $P_1$ is also left aligned, but there is a half-box shift between the rows $(a_1, \dots, a_r)$ and the hook-type part. Fix an arbitrary labeling of these boxes. The hook-type part of $P_2$ is built in two steps:
\begin{enumerate}
    \item shift the hook-type part $(a, 1^b)$ of $P_1$ by half a box to the right, 
    
    \item move the boxes corresponding to the part $(1^b)$ by one box to the left and one box down.
\end{enumerate} 

\begin{example} \label{example:partitions}
    Consider the following partitions of $n = 9$:
    $$ a_\bullet^{(1)} = (4, 4, 2, 1^3) \quad \text{and} \quad a_\bullet^{(2)} = (4, 4, 3, 1^2). $$
    The associated pyramids are drawn in \Cref{figure:pyramids}. 
\end{example}

\setlength{\unitlength}{0.7pt}

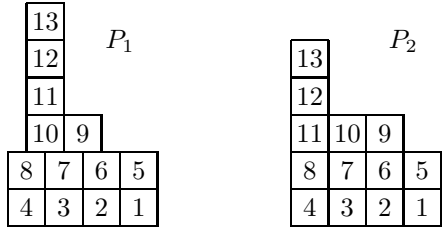
\begin{figure}[t]
    \caption{Examples of pyramids}
    \label{figure:pyramids}
    $$ \begin{picture}(80,120)
        % horizontal
	\put(0,0){\line(1,0){80}}
        \put(0,20){\line(1,0){80}}
        \put(0,40){\line(1,0){80}}
        \put(10,60){\line(1,0){40}}
        \put(10,80){\line(1,0){20}}
        \put(10,100){\line(1,0){20}}
        \put(10,120){\line(1,0){20}}
        % vertical lines
        \put(0,0){\line(0,1){40}}
        \put(10,40){\line(0,1){80}}
        \put(20,0){\line(0,1){40}}
        \put(30,40){\line(0,1){80}}
        \put(40,0){\line(0,1){40}}
        \put(50,40){\line(0,1){20}}
        \put(60,0){\line(0,1){40}}
        \put(80,0){\line(0,1){40}}

        \put(60,100){\makebox(0,0){{$P_1$}}}
	\put(70,10){\makebox(0,0){{$1$}}}
        \put(50,10){\makebox(0,0){{$2$}}}
        \put(30,10){\makebox(0,0){{$3$}}}
        \put(10,10){\makebox(0,0){{$4$}}}
        \put(70,30){\makebox(0,0){{$5$}}}
        \put(50,30){\makebox(0,0){{$6$}}}
        \put(30,30){\makebox(0,0){{$7$}}}
        \put(10,30){\makebox(0,0){{$8$}}}
        \put(40,50){\makebox(0,0){{$9$}}}
        \put(20,50){\makebox(0,0){{$10$}}}
        \put(20,70){\makebox(0,0){{$11$}}}
        \put(20,90){\makebox(0,0){{$12$}}}
        \put(20,110){\makebox(0,0){{$13$}}}
    \end{picture} \hspace{5em}  \begin{picture}(80,120)
        % horizontal
	\put(0,0){\line(1,0){80}}
        \put(0,20){\line(1,0){80}}
        \put(0,40){\line(1,0){80}}
        \put(0,60){\line(1,0){60}}
        \put(0,80){\line(1,0){20}}
        \put(0,100){\line(1,0){20}}
        % vertical lines
        \put(0,0){\line(0,1){100}}
        \put(20,0){\line(0,1){100}}
        \put(40,0){\line(0,1){60}}
        \put(60,0){\line(0,1){60}}
        \put(80,0){\line(0,1){40}}

        \put(60,100){\makebox(0,0){{$P_2$}}}
	\put(70,10){\makebox(0,0){{$1$}}}
        \put(50,10){\makebox(0,0){{$2$}}}
        \put(30,10){\makebox(0,0){{$3$}}}
        \put(10,10){\makebox(0,0){{$4$}}}
        \put(70,30){\makebox(0,0){{$5$}}}
        \put(50,30){\makebox(0,0){{$6$}}}
        \put(30,30){\makebox(0,0){{$7$}}}
        \put(10,30){\makebox(0,0){{$8$}}}
        \put(50,50){\makebox(0,0){{$9$}}}
        \put(30,50){\makebox(0,0){{$10$}}}
        \put(10,50){\makebox(0,0){{$11$}}}
        \put(10,70){\makebox(0,0){{$12$}}}
        \put(10,90){\makebox(0,0){{$13$}}}
    \end{picture}$$
\end{figure}

\begin{proposition}
    For $i=1,2$, let $f_i$ be a nilpotent matrix and $H_i$ be the diagonal element in $\sl_n$ corresponding to the pyramid $P_i$ described above. Then, the pairs~$(f_1, H_1)$ and $(f_2, H_2)$ satisfy the condition \eqref{conditions} and \Cref{main:reduction-by-stages-w-algebras} holds.
\end{proposition}

\begin{proof}
    From the left to the right, and from the bottom to the top, denote by $j_1, \dots, j_{a}, j_{a+1}, \dots, j_{a+b}$ the labels in the hook-type part (in \Cref{example:partitions}, these labels are $9,10,11,12,13$).

    Take $1 \leqslant j, k  \leqslant n$. Let $x^{(i)}_j$, respectively $x^{(i)}_k$, be the abscissa of the center of the box labeled with $j$, respectively with $k$. Then
    $$ [H_i, E_{j,k}] = (x^{(i)}_j - x^{(i)}_k) E_{j,k}. $$

    The nilpotent element $f_i$ is the sum of elementary matrices $E_{j,k}$ such that the boxes labeled by $j$ and $k$ are in the same row and $x^{(i)}_j - x^{(i)}_k = -2$. Then $ f_0 = f_2 - f_1$ is equal to $E_{j_{a+1}, j_a}$ and since $x^{(i)}_{j_{a+1}} = x^{(i)}_{j_a}$, one has $[H_1, f_0] = 0$.
    
    Let us prove that for any labels $j, k$, 
    \begin{align*}
        & \text{if} \quad x^{(1)}_{j} - x^{(1)}_k \geqslant 2 \quad \text{then} \quad x^{(2)}_{j} - x^{(2)}_k \geqslant 1, \\
        & \text{if} \quad x^{(2)}_{j} - x^{(2)}_k \geqslant 1 \quad \text{then} \quad x^{(1)}_{j} - x^{(1)}_k \geqslant 0, \\
        & \text{if} \quad |x^{(1)}_{j} - x^{(1)}_k | =  1 \quad \text{then} \quad |x^{(1)}_{j} - x^{(1)}_k| \in \{0,1,2\}.
    \end{align*}
    If $j$ labels one box in the hook-type part and $k$ labels a box in the first $r$ rows, it is because the relative position between a box in the hook-type part and a box in the first $r$ rows only vary of half a box when $P_1$ is turned into $P_2$. The other cases are clear. 
    
    Hence the conditions \eqref{conditions} are all satisfied.
\end{proof}

\begin{remark}
    This family of examples contains cases which are already known, see \Cref{table:examples}. However, some reduction by stages from \cite{fasquel2026connecting, fasquel2025virasoro} are not covered by this family. For example, $\g = \sl_4$, $f_1$ being a rectangular nilpotent element (partition~$(2,2)$) and $f_2$ being a subregular one. 
\end{remark}

\subsubsection{Other examples} The computations done in {\cite[Section 4]{genra2024reduction}} prove the following proposition.

\begin{proposition}
    \Cref{main:reduction-by-stages-w-algebras} holds in the following cases:
    \begin{enumerate}
        \item the Lie algebra $\g$ is of type $\mathrm B$, $f_1$ is a subregular nilpotent element and $f_2$ is regular,

        \item the Lie algebra $\g$ is of type $\mathrm{C}_r$ ($r \geqslant 3$), $f_1$ is a nilpotent element of partition~$(2^2, 1^{2r - 4})$ and $f_2$ is regular,

        \item the Lie algebra $\g$ is of type $\mathrm{G}_2$ ($r \geqslant 3$), $f_1$ is a nilpotent element of Bala-Carter label $\widetilde{\mathrm{A}}_1$ and $f_2$ is subregular.
    \end{enumerate}
\end{proposition}

\printbibliography
	
\end{document}